\documentclass{article}
\usepackage[utf8]{inputenc}
\usepackage[margin=1in]{geometry} 
\usepackage{amsmath,amsthm,amssymb,amsfonts}
\usepackage{tikz}
\usetikzlibrary{shapes.geometric, arrows}
\usepackage{verbatim}
\usepackage{dsfont}
\usepackage{cite}
\usepackage{algpseudocode}
\usepackage{algorithm}
\usepackage{subfigure}
\usepackage[colorlinks=true]{hyperref}
\hypersetup{urlcolor=blue, citecolor=red, linkcolor=blue}
\usepackage[capitalise,noabbrev,nameinlink]{cleveref}
\usepackage[bottom]{footmisc}
\usepackage{multirow}
\usepackage[group-digits=false]{siunitx}
\usepackage{mwe}
\usepackage{mathtools}
\usepackage[inline]{enumitem}
\usepackage[]{mdframed}
\usepackage[skins,breakable]{tcolorbox}

\theoremstyle{plain}
\newtheorem{thm}{Theorem}[section] 
\newtheorem*{theorem*}{Main Result}
\newtheorem{lem}[thm]{Lemma}
\newtheorem{prop}[thm]{Proposition}

\newtheorem{cor}[thm]{Corollary}

\theoremstyle{definition}
\newtheorem{defn}[thm]{Definition} 
 
\newtheorem{assump}[thm]{Assumption}
\newtheorem{remark}[thm]{Remark}

\numberwithin{equation}{section}

\DeclareMathOperator*{\maximize}{\text{maximize}}
\DeclareMathOperator*{\argmin}{arg\,min}

\tikzstyle{dnc_step} = [rectangle, rounded corners, minimum width=3cm, minimum height=.7cm,text centered, text width = 4.5cm, draw=white]
\tikzstyle{arrow} = [draw, -latex',rounded corners]
\tikzstyle{point} = [circle, minimum width=.5cm, minimum height=.5cm, text centered, draw=black, fill=gray!30]
\tikzstyle{line} = [draw, -latex']

\usetikzlibrary{shapes,calc,arrows}

\begin{document}
\title{Fast and Inverse-Free Algorithms for Deflating Subspaces}
\author{James Demmel\footnote{Department of EECS (Computer Science Division) and Department of Mathematics, University of California Berkeley} \and Ioana Dumitriu\footnote{Department of Mathematics, University of California San Diego}
\and Ryan Schneider\footnote{Department of Mathematics, University of California Berkeley (ryan.schneider@berkeley.edu)} }
\date{}

\maketitle

\begin{abstract}
This paper explores a key question in numerical linear algebra:\ how can we compute projectors onto the deflating subspaces of a regular matrix pencil $(A,B)$, in particular without using matrix inversion or defaulting to an expensive Schur decomposition? We focus specifically on \textit{spectral} projectors, whose associated deflating subspaces correspond to sets of eigenvalues/eigenvectors. In this work, we present a high-level approach to this computational problem, which combines rational function approximation with an inverse-free arithmetic of Benner and Byers [Numerische Mathematik 2006]. The result is a numerical framework that captures existing inverse-free methods, generates an array of new options, and provides straightforward tools for pursuing efficiency on structured problems (e.g., definite pencils). To exhibit the efficacy of this framework, we consider a handful of methods in detail, including Implicit Repeated Squaring and iterations based on the matrix sign function. In an appendix, we demonstrate that recent, randomized divide-and-conquer eigensolvers -- which are built on fast methods for individual projectors -- can be adapted to produce the generalized Schur form of any matrix pencil in nearly matrix multiplication time.
\end{abstract}

\small
\hspace{2mm} \textbf{Keywords:} Matrix pencil, spectral projectors, deflating subspaces, repeated squaring, matrix sign 

\hspace{20.7mm} function, generalized Schur decomposition  \\
\vspace{2mm}
\indent \hspace{1mm} \textbf{MSC Class:} 65F15 65F60
\normalsize

\vspace{2mm}
\pagebreak

\tableofcontents

\pagebreak

\section{Introduction}\label{section: intro}

We consider the task of computing projectors onto the right/left deflating subspaces of an arbitrary, regular matrix pencil $(A,B) \in {\mathbb C}^{n \times n} \times {\mathbb C}^{n \times n}$. We follow the standard definition for these subspaces:\ ${\mathcal X},{\mathcal Y} \subseteq {\mathbb C}^n$ are right and left deflating subspaces of $(A,B)$, respectively, if $\text{dim}({\mathcal X}) = \text{dim}({\mathcal Y})$ and 
\begin{equation}\label{eqn: deflating_subspace_def}
    \text{span} \left\{ Ax, Bx : x \in {\mathcal X} \right\} = {\mathcal Y}.
\end{equation}
When ${\mathcal X}$ is spanned by a set of right eigenvectors of $(A,B)$, projectors onto ${\mathcal X}$ and ${\mathcal Y}$ are called \textit{spectral projectors} of $(A,B)$.\footnote{Any set of right eigenvectors is guaranteed to span a corresponding right deflating subspace, which is easy to see from the generalized Schur form of $(A,B)$ \cite{Stewart_Schur}. Despite the similarity in naming, left deflating subspaces of $(A,B)$ are not usually spanned by left eigenvectors.} If additionally $B = I$, or more generally if $A$ and $B$ commute, then these projectors are the same, corresponding to an invariant eigenspace of $A$. Deflating subspaces are only defined if $(A,B)$ is regular, which here means that the characteristic polynomial $p(\lambda) = \det(A - \lambda B)$ is \textit{not} identically zero. Since it will be important to have in mind throughout, note that the pencil $(A,B)$ and the matrix $B^{-1}A$ have the same eigenvalues and right eigenvectors, provided $B$ is invertible. \\
\indent In this introduction, we assume familiarity with standard definitions from linear algebra. For a summary of the relevant notation see \cref{section: notation}.

\subsection{Motivation}
\indent Spectral projectors and their associated deflating subspaces are essential in numerical linear algebra. While specific projectors/subspaces are of interest in certain applications --  e.g., corresponding to the largest or smallest $k$ eigenvalues as in \cite{hermitian_shattering} -- our primary motivation is their use in divide-and-conquer eigensolvers, which recursively diagonalize a matrix pencil (or individual matrix, setting $B=I$) as follows:

\begin{figure}[h]
   \centering
   \begin{tikzpicture}[node distance=1.2cm]
     \node (base) [dnc_step] {$(A,B)$};
     \node (sub1) [dnc_step, below of = base, xshift = 3cm] {$((U_L^{(2)})^HAU_R^{(2)}, (U_L^{(2)})^HBU_R^{(2)})$};
     \node (sub2) [dnc_step, below of = base, xshift = -3cm]
     {$((U_L^{(1)})^HAU_R^{(1)}, (U_L^{(1)})^HBU_R^{(1)})$};
    \path [arrow] (base) |- ($(base)+(0,-.4)$) -| (sub2);
    \path [arrow] (base) |- ($(base)+(0,-.4)$) -| (sub1);
   \end{tikzpicture}
\end{figure}

\indent Here, $U_R^{(1)},U_L^{(1)} \in {\mathbb C}^{n \times k}$ and $U_R^{(2)},U_L^{(2)} \in {\mathbb C}^{n \times (n-k)}$ contain orthonormal bases for pairs of right/left deflating subspaces corresponding to disjoint sets of eigenvalues. Note that $U_R^{(1)}$ and $U_R^{(2)}$ are not orthogonal to one another in general (and the same is true for $U_L^{(1)}$ and $U_L^{(2)}$). Typically, these matrices are obtained by computing rank-revealing factorizations of the associated spectral projectors. In this way, spectral projectors are the key ingredient of the block diagonalization procedure. \\ 
\indent While divide-and-conquer has existed in the literature for several decades \cite{BEAVERS1974143,bulgakov,MALYSHEV,Bai:CSD-94-793,Ballard2010MinimizingCF}, it has only recently been formulated in a way that provably succeeds (and achieves near-optimal complexity) on arbitrary inputs \cite{banks2020pseudospectral,arXiv}. As efforts to deploy these algorithms move forward, 
it is increasingly important to optimize methods for computing projectors, as the step of divide-and-conquer that relies on them tends to dominate in both computational time and precision requirements.\footnote{Indeed, divide-and-conquer may require as many as $O(n)$ projectors to diagonalize a single pencil.} With this in mind, we focus specifically on the following problem, which -- as is typical in divide-and-conquer -- assumes that we have already identified a region of the complex plane containing a piece of the spectrum of $(A,B)$.
\begin{mdframed}
\noindent \textbf{Problem Statement.} Given an arbitrary, regular matrix pencil $(A,B)$ and a set $S \subseteq {\mathbb C}$ containing some subset of eigenvalues of $(A,B)$, compute the projectors $P_R$ and $P_L$ onto the right/left deflating subspaces corresponding to $S$, where the former is spanned by (right) eigenvectors associated to eigenvalues in $S$ and the latter is defined according to \eqref{eqn: deflating_subspace_def}. 
\end{mdframed}

\indent Our primary contribution in this paper is a general framework for developing fast, inverse-free algorithms that can solve this problem. Here, ``fast" means that each method requires at most $O(\log(\frac{n}{\delta}) T_{\text{MM}}(n))$ operations to compute $P_R$ and $P_L$ to forward accuracy $\delta$ in the spectral norm -- assuming the problem is not too ill-conditioned -- for $T_{\text{MM}}(n)$ the complexity of $n \times n$ matrix multiplication. Several such algorithms already exist in the literature (and are summarized in the subsequent sections); our framework provides a means of understanding/adapting them and generating new options. \\
\indent In this context, conditioning depends primarily on the largest $\epsilon > 0 $ for which the $\epsilon$-pseudospectrum of $(A,B)$, defined appropriately (see \cref{defn:pencil_pseudospectrum}), is well-separated from the boundary of $S$. To achieve the ``fast" moniker this $\epsilon$ will need to be at least polynomial in $n^{-1}$ and $\delta$, as we discuss further in \cref{section: SIGN}. Importantly, this notion of ill-conditioning speaks to complexity and not feasibility/accuracy; the methods we present can still produce accurate approximations even if the aforementioned pseudospectral condition does not hold, though in that case they may exceed $O(\log(\frac{n}{\delta}) T_{\text{MM}}(n))$ complexity. A problem is ill-posed only when $(A,B)$ has an eigenvalue on the boundary of $S$.
\begin{remark}
    Absent these fast methods, computing $P_R$ and $P_L$ typically requires obtaining a full Schur decomposition of $(A,B)$ and reorganizing it so that the leading eigenvalues belong to $S$ -- i.e., \texttt{xGGES} in LAPACK.\footnote{See \cite{Kagstrom_poromaa} for a discussion of the reorganization procedure.} Indeed, computing a generalized Schur decomposition is equivalent to computing a sequence of nested deflating subspaces. With this in mind, and because Schur form itself has a variety of applications \cite{Demmel_Kagstrom1,Demmel_Kagstrom2,CANDECOMP_Schur,Lindley_Schur}, we discuss it in more detail in \cref{section: fast_schur}. There we demonstrate that the aforementioned divide-and-conquer approach can be adapted to produce a Schur decomposition (of any pencil in nearly matrix multiplication time) by leveraging the fast methods for individual deflating subspaces presented in the rest of the paper. 
\end{remark}

\subsection{High-Level Strategy} 
At a high level, our goal will be to transform $(A,B)$ into the pencil $(A_S,B_S)$, which has the same right eigenvectors as $(A,B)$ but whose eigenvalues are zero and one, depending on whether or not the corresponding eigenvalues of $(A,B)$ belong to $S$. If this can be done, computing projectors is straightforward. In particular $P_R = B_S^{-1}A_S$, while $P_L^H$ can be obtained by repeating this procedure with $(A^H,B^H)$ and $S^* = \left\{ \overline{z} : z \in S \right\}$.\footnote{This follows from the observation that left eigenvectors of $B^{-H}A^H$ associated to $S$, which are not necessarily left eigenvectors of either $(A,B)$ or $(A^H,B^H)$, span the appropriate left deflating subspace.} Since a rank-revealing QR factorization of $B_S^{-1}A_S$ can be computed implicitly, for example using the \textbf{GRURV} algorithm of Ballard et al.\ \cite{grurv}, accomplishing the transformation $(A,B) \rightarrow (A_S,B_S)$ is sufficient for computing $P_R$ and $P_L$ without inversion. \\
\indent In practice, obtaining $(A_S,B_S)$ from $(A,B)$ naturally reduces to (rationally) approximating the indicator function
\begin{equation} \label{eqn: indicator}
    \mathds{1}_S(z) = \begin{cases}
      1 & z \in S \\
      0 & z \in {\mathbb C} \backslash \overline{S} \\
      \text{undefined} & z \in \partial S
    \end{cases}
\end{equation}
where $\overline{S}$ denotes the closure of $S$. Evaluating such an approximation $r(z)$ at $B^{-1}A$ (without taking inverses or forming the product) will yield an approximation of $(A_S,B_S)$, as doing so preserves right eigenvectors while mapping eigenvalues $\lambda$ to $r(\lambda)$. That this can be done implicitly follows from work of Benner and Byers \cite{BB2001,Benner_Byers}, which develops a general, inverse-free arithmetic on matrix pencils. The central definition of their arithmetic is the \textit{matrix relation}
\begin{equation}\label{eqn: matrix_relation}
    (B \backslash A) = \left\{ (x,y) \in {\mathbb C}^n \times {\mathbb C}^n : Ax = By \right\},
\end{equation}
which we can think of as a  representation of $B^{-1}A$ that neither requires $B$ to be invertible nor risks instability if $B$ \textit{is} invertible but ill-conditioned.\footnote{This representation is not unique. We can left multiply $A$ and $B$ by any matrix $M$ whose null space only trivially overlaps with the range of $\begin{bmatrix} A & B \end{bmatrix}$ without changing the relation.} The corresponding arithmetic is defined by the following operations. 
\begin{defn}[Arithmetic for Matrix Relations]\label{defn: matrix_relation_op}
    The \textit{sum} and \textit{product} of two matrix relations $(B_1 \backslash A_1)$ and $(B_2 \backslash A_2)$ are subsets of ${\mathbb C}^n \times {\mathbb C}^n$ given by
    $$ \aligned 
    (B_2 \backslash A_2) + (B_1 \backslash A_1) &=  \left\{ (x,z) : \exists \; \; y_1,y_2 \; \; \text{s.t.} \; \begin{pmatrix} A_1 & -B_1 & 0 & 0 \\ A_2 & 0 & -B_2 & 0 \\ 0 & I & I & -I \end{pmatrix} \begin{pmatrix} x \\ y_1 \\ y_2 \\ z \end{pmatrix} =  0 \right\} \\
    (B_2 \backslash A_2) (B_1 \backslash A_1) &= \left\{ (x,z) : \exists \; \; y \; \; \text{s.t.} \; \begin{pmatrix} A_1 & -B_1 & 0 \\ 0 & A_2 & -B_2\end{pmatrix} \begin{pmatrix} x \\ y \\ z \end{pmatrix} = 0 \right\}.
    \endaligned $$
\end{defn}
\indent Observing that a block QR factorization
\begin{equation}\label{eqn: block_QR}
    \begin{pmatrix} A \\ B \end{pmatrix} = \begin{pmatrix} Q_{11} & Q_{12} \\ Q_{21} & Q_{22} \end{pmatrix} \begin{pmatrix} R \\ 0 \end{pmatrix} 
\end{equation}
implies $\text{null} \left( \begin{bmatrix} Q_{12}^H & Q_{22}^H \end{bmatrix} \right) = \text{range} \left( \begin{bmatrix} A \\ B \end{bmatrix} \right)$, yet another result of Benner and Byers suggests that the sum and product of any pair of matrix relations can be computed using only QR and matrix multiplication \cite[Theorems 2.3 and 2.7]{Benner_Byers}. 
\begin{thm}[Benner and Byers \cite{Benner_Byers}]\label{thm: matrix_relation_comp}
    Let $(B_1 \backslash A_1)$ and $(B_2 \backslash A_2)$ be two matrix relations with $A_1,A_2,B_1,B_2 \in {\mathbb C}^{n \times n}$. Suppose
    $$ \text{\normalfont null} \left( \begin{bmatrix} Q_1 & Q_2\end{bmatrix} \right) = \text{\normalfont range}\left(\begin{bmatrix} -B_1 \\ A_2 \end{bmatrix} \right)\; \text{and}\; \; \text{\normalfont null}\left( \begin{bmatrix} U_1 & U_2 \end{bmatrix} \right)) = \text{\normalfont range}\left( \begin{bmatrix} -B_1 \\ B_2 \end{bmatrix} \right) .$$
    Then 
    $$ (B_2 \backslash A_2)(B_1 \backslash A_1) = ((Q_2B_2) \backslash (Q_1A_1)) $$
    and 
    $$ (B_2 \backslash A_2) + (B_1 \backslash A_1) = ((U_2B_2) \backslash (U_1A_1 + U_2A_2)).$$
\end{thm}
\indent It can be shown that these operations extend to the spectra of $(A_1,B_1)$ and $(A_2,B_2)$ in a natural way; if $\lambda$ and $\mu$ are eigenvalues of these pencils associated to a shared (right) eigenvector $v$, then $(\lambda + \mu,v)$ and $(\lambda \mu, v)$ are eigenpairs of the pencils corresponding to $(B_2 \backslash A_2) + (B_1 \backslash A_1)$ and $(B_2 \backslash A_2)(B_1 \backslash A_1)$, respectively\footnote{Technically, this only holds if $(A_1,B_1)$ and $(A_2,B_2)$ are both regular and $\left\{ \lambda , \mu \right\} \neq \left\{ \infty, 0 \right\}$, though this will always be the case for our purposes.} \cite[Theorems 2.5 and 2.8]{Benner_Byers}. Hence, a polynomial can be applied to a regular pencil $(A,B)$ by evaluating it at the relation $(B \backslash A)$, which implicitly evaluates the polynomial at $B^{-1}A$ and maps eigenvalues accordingly. \\
\indent To extend this to rational functions we need only introduce some notion of a multiplicative inverse, which can be done as $(B \backslash A)^{-1} = (A \backslash B)$. While this is only a true inverse when $B^{-1}A$ exists and is invertible, it is sufficient for our purposes since the eigenvalues of $(B,A)$ are clearly the reciprocals of the eigenvalues of $(A,B)$.\footnote{Note that if $(A,B)$ has an eigenvalue at zero then $(B,A)$ has a corresponding eigenvalue at infinity.} \\
\indent Given any rational function $r(z) = p(z)/q(z)$ for polynomials $p$ and $q$, we can now evaluate $r(B \backslash A)$ using the arithmetic of Benner and Byers as 
\begin{equation}\label{eqn: evaluate_rational}
    r(B \backslash A) = (q(B \backslash A))^{-1} p(B \backslash A). 
\end{equation}
Here we make a somewhat arbitrary choice; we could evaluate $r(B \backslash A)$ as $p(B \backslash A)(q(B \backslash A))^{-1}$, though this simply yields a different representation of the same matrix relation. Sticking with \eqref{eqn: evaluate_rational}, we obtain the following high-level, inverse-free framework for computing $(A_S,B_S)$:
\begin{enumerate}
    \item Approximate the indicator function $\mathds{1}_S(z)$ with a rational function $r(z)$.
    \item Evaluate $r(B \backslash A)$ via \cref{thm: matrix_relation_comp} and \eqref{eqn: block_QR}.
    \item Set $(B_S \backslash A_S) = r(B \backslash A)$.
\end{enumerate}

\indent In this framework, each choice of $r(z)$ generates a different numerical method for solving our original problem. Choosing between them suggests the following.
\begin{mdframed}
    \textbf{The Indicator Approximation Problem:} Given $S \subseteq {\mathbb C}$, what is the ``best" rational function approximation to $\mathds{1}_S(z)$?
\end{mdframed}

\indent This question is more subtle than it initially appears. While in general we prefer approximations that are as close to $\mathds{1}_S(z)$ as possible -- in an appropriate norm and at least near the eigenvalues of $(A,B)$ if not on all of $S$ -- we must also account for the expensive nature of Benner and Byer's inverse-free arithmetic. That is, each addition/multiplication required to evaluate a given rational function requires a block $2n \times n$ QR factorization. Hence, a better approximation to $\mathds{1}_S(z)$ may not yield a more efficient method if it requires even one or two additional operations. There is one exception here:\ Möbius transformations; $r(z) = \frac{az+b}{cz+d}$ can be applied to $(B \backslash A)$ for free as $((cA+dB)\backslash (aA + bB))$, and it will therefore be advantageous to write $r(z)$ in terms of Möbius transformations whenever possible. \\
\indent Regardless of the choice of $r(z)$, the framework outlined above naturally promotes stability in floating-point arithmetic by avoiding matrix inversion. This is particularly valuable when working with a pencil $(A,B)$ in which $B$ and/or $A$ is singular or nearly singular. Moreover, it leaves the door open to implement any of the methods considered here in a communication-optimal fashion (in the vein of Ballard et al.\ \cite{Ballard_min_comm}). These concerns are especially relevant for divide-and-conquer eigensolvers, where the benefits of avoiding inversion have already been explored -- e.g., \cite{arXiv,Ballard2010MinimizingCF,My_thesis}.

\subsection{Contributions}
In the remainder of the paper, we discuss a handful of specific instances of our high-level framework, which are summarized in \cref{tab: iteration_comp}. Our goal throughout is to produce rigorous performance guarantees for each method. Since all four can be implemented iteratively -- either by implementing the composition of a fixed rational function or via iterative squaring -- relative efficiency will be determined by the operation count per iteration and the number of iterations required to reach a given level of accuracy. As we will see, executing $k$ iterations of each method requires $O(kT_{\text{MM}}(n))$ operations. \\
\indent The main contributions of the remaining sections can now be summarized as follows:
\begin{itemize}
    \item In \cref{section: IRS}, we consider the Implicit Repeated Squaring (\textbf{IRS}) routine of Malyshev \cite{MALYSHEV}, arguably the most widely used inverse-free method for computing projectors. We include in this section a new, general stability bound for the method in floating-point arithmetic, which is based on a standard finite-precision model that can accommodate fast matrix multiplication. 
    \item \cref{section: SIGN} presents inverse-free methods based on the matrix sign function of Beavers and Denman \cite{BEAVERS1974143}. There, we use the Indicator Approximation Problem to refine one of these methods -- the Halley iteration -- in an effort to improve efficiency on problems with real eigenvalues.\footnote{This includes the important definite generalized eigenvalue problem. See \cite[Section VI.3]{stewart1990matrix} for further background.} The result is a new (generalized) dynamically weighted Halley iteration.   
    \item In \cref{section: examples}, we provide a handful of numerical examples that test the methods explored in the preceding sections. 
     \item Finally, \cref{section: fast_schur} presents the aforementioned discussion of divide-and-conquer for Schur form, \cref{section: appendix} contains the proof of the finite-precision stability bound for \textbf{IRS}, \cref{appendix: Zolo} discusses the viability of generalizing the results of \cref{section: SIGN}, and \cref{appendix: examples} covers additional numerical examples.
\end{itemize}

\indent In total, this work demonstrates the efficacy of the previous section's high-level framework, both as a tool to understand and adapt existing methods and to produce new ones. To that end, we hope it prompts further development of specialized routines for this problem.

\renewcommand{\arraystretch}{1.4}
\begin{table}[t]
    \centering
    \begin{tabular}{lll}
        \hline
         Method &  $S$ & Approximation $r(z)$ \\
         \hline
         Implicit Repeated Squaring & $|z|<1$ & $(1+z^{2^k})^{-1}$ \\
         Newton Iteration & $\text{Re}(z) > 0$ & $\frac{1}{2} \left[1+f \circ \cdots \circ f(z) \right]$ with $f(z) = \frac{1}{2}(z+z^{-1})$\\
         Halley Iteration & $\text{Re}(z) > 0$ & $\frac{1}{2}\left[ 1+ f \circ \cdots \circ f(z)\right]$ with $f(z) = z \frac{z^2+3}{3z^2+1}$\\
         Dynamically Weighted Halley Iteration & $\text{Re}(z) > 0$ & $\frac{1}{2} \left[1+f_k \circ \cdots \circ f_1(z)\right]$ with $f_i(z) = z \frac{a_i z^2 + b_i}{c_iz^2 + d_i}$ \\
         \hline 
    \end{tabular}
    \caption{Methods of approximating $\mathds{1}_S(z)$ for different choices of $S$. The Newton and Halley iterations are based on the (complex) scalar sign function \eqref{eqn: scalar_sign}, as $\mathds{1}_S(z) = \frac{1}{2}(\text{sign}(z)+1)$ if $S$ is the right half plane.}
    \label{tab: iteration_comp}
\end{table}
\renewcommand{\arraystretch}{1}

\subsection{Notation and Conventions}\label{section: notation}
Throughout the paper, we use $(A,B) \in {\mathbb C}^{n \times n} \times {\mathbb C}^{n \times n}$ to denote a square matrix pencil, which is assumed to be regular. $A^H$ and $A^{-H}$ denote the Hermitian transpose and inverse Hermitian transpose, respectively, while $I$ is used to denote the identity matrix, with size implied by context. $\|\cdot\|_2$ is the spectral norm on matrices and the Euclidean norm on vectors, and $\kappa_2(\cdot)$ is the spectral norm condition number. $\Lambda(A,B)$ is used to denote the spectrum of $(A,B)$. Finally, $\sigma_i(A)$ is the $i$-th singular value of $A$, though when convenient $\sigma_{\min}(A)$ may be used to denote the smallest singular value of $A$ instead. 
\section{Implicit Repeated Squaring}\label{section: IRS}
We start with Implicit Repeated Squaring (\textbf{IRS}), a routine for repeatedly squaring a product $A^{-1}B$ without forming it. \textbf{IRS} originates in early efforts to implement inverse-free, divide-and-conquer eigensolvers, first in work\footnote{The paper \cite{Malyshev_original} was translated from Russian in two parts \cite{Malyshev1,Malyshev2}. Much of its content was subsequently presented in \cite{MALYSHEV}.} of Malyshev \cite{Malyshev1989,Malyshev_original} and later Bai, Demmel, and Gu \cite{Bai:CSD-94-793}. It eventually appeared as in \cref{alg:IRS} under the name \textbf{IRS} in work of Ballard, Demmel, and Dumitriu \cite{Ballard2010MinimizingCF}.

\begin{algorithm}
\caption[Implicit Repeated Squaring (\textbf{IRS})]{Implicit Repeated Squaring (\textbf{IRS})\\
\textbf{Input:} $A, B \in {\mathbb C}^{n \times n}$ and $p$ a positive integer. }\label{alg:IRS}
\begin{algorithmic}[1]
\State $A_0 = A$
\State $B_0 = B$
\For{$j = 0:p-1$}
    \vspace{1mm}
    \State $\begin{pmatrix} B_j \\ - A_j 
    \end{pmatrix}  = \begin{pmatrix} Q_{11} & Q_{12} \\
    Q_{21} & Q_{22} \end{pmatrix} \begin{pmatrix} R_{j} \\ 0  \end{pmatrix}$ 
    \vspace{1mm}
    \State $A_{j+1} = Q_{12}^H A_j$
    \State $B_{j+1} = Q_{22}^H B_j$
\EndFor
\State \Return $A_p, B_p$
\end{algorithmic}
\end{algorithm}

\indent \textbf{IRS} can be used to compute projectors by applying our framework with $r(z) = (1+z^{2^p})^{-1}$ and $S = \left\{z : |z| < 1 \right\}$. In these terms, the pseudocode of \cref{alg:IRS} can be viewed as a straightforward application of \cref{thm: matrix_relation_comp} to $(A_p \backslash B_p) = (A \backslash B)^{2^p}$, where squaring naturally drives eigenvalues to zero and infinity (assuming none are on the unit circle). Applying the Möbius transformation $(1 + z)^{-1}$, which sends zero to one and infinity to zero, the projector $P_R$ can be obtained from $((A_p+B_p) \backslash A_p)$ as 
\begin{equation}
    P_R \approx (A_p+B_p)^{-1}A_p.
    \label{eqn: IRS_right_proj}
\end{equation}
Repeating this process with $(A^H,B^H)$ yields the left projector $P_L \approx {\mathcal A}_p^{H}({\mathcal A}_p+{\mathcal B}_p)^{-H}$ for $({\mathcal A}_p \backslash {\mathcal B}_p) = (A^H \backslash B^H)^{2^p}$. Note here that \textbf{IRS} is applied to $(A \backslash B)$ rather than $(B \backslash A)$. This is done to maintain consistency with the presentation of \textbf{IRS} in \cite{Bai:CSD-94-793, Ballard2010MinimizingCF}, though it also means that $P_R$ and $P_L$ are spectral projectors of $(A,B)$ corresponding to $\left\{ z : |z| >1 \right\}$ rather than $S$. To avoid confusion, we label the projectors as $P_{R, |z|>1}$ and $P_{L, |z| > 1}$ to make clear the subset of the spectrum of $(A,B)$ they depend on. 

\subsection{Condition Number}
\indent Accuracy guarantees for \textbf{IRS} have been derived by both Malyshev \cite{MALYSHEV} and Bai, Demmel, and Gu \cite{Bai:CSD-94-793}. The latter bounds $\|(A_p+B_p)^{-1}A_p - P_{R, |z|>1}\|_2$ in terms of the reciprocal of a distance to the nearest ill-posed problem $d_{(A,B)}$, which is defined as follows (see \cite[Theorem 1]{Bai:CSD-94-793})



\begin{defn}\label{defn: ill-posed}
The distance from $(A,B)$ to the nearest ill-posed problem is 
$$d_{(A,B)} = \inf \left\{ \|E\|_2 + \|F\|_2 : (A + E) - z(B + F) \; \text{is singular for some} \; |z| =1 \right\}.  $$ 
\end{defn}

This quantity is specialized to the setting where \textbf{IRS} is employed to compute spectral projectors. Indeed, $d_{(A,B)}^{-1}$ is infinite if $(A,B)$ is singular or has an eigenvalue on the unit circle, in which case squaring cannot successfully produce a projector by driving eigenvalues to zero or infinity. As a result, $d_{(A,B)}^{-1}$, or more precisely the scale-invariant quantity $\|(A,B)\|_2/d(A,B)$, can be interpreted as a condition number for repeated squaring that is specialized to the spectral projector application. The same can be said of the condition number used in Malyshev's analysis \cite{MALYSHEV}, which is also infinite if $(A,B)$ has an eigenvalue on the unit circle. \\
\indent We seek something more general here. Our motivation lies in the potential for \textbf{IRS} to be applied to other problems in numerical linear algebra. Take for example the matrix exponential $e^A$. The most commonly used algorithm for computing $e^A$ is the scaling and squaring method \cite{scaling_and_squaring}, which evaluates the exponential as
\begin{equation}\label{eqn: scaling_and_squaring}
        e^A \approx \left[ q(A/2^p)^{-1} p(A/2^p) \right]^{2^p}
\end{equation}
for two polynomials $p$ and $q$. Clearly, \textbf{IRS} can be used to handle the squaring step of this approach without relying on inversion. Nevertheless, the aforementioned condition number is ill-suited to capture its performance, since eigenvalues of $(q(A/2^p), p(A/2^p))$ on the unit circle are irrelevant. \\
\indent With this in mind, we choose to work with the following condition number instead. While its definition may seem somewhat contrived, we justify our choice in the next subsection, which presents a general floating-point stability bound for \textbf{IRS} in terms of $\kappa_{\text{IRS}}(A,B,p)$. 
\begin{defn}\label{defn: irs_condition_num}
    Given $A,B \in {\mathbb C}^{n \times n}$ and $p \geq 1$, define the block matrix
    $$ D_{(A,B)}^p = \begin{pmatrix} B & & & \\
        -A & B & & \\
        & -A & \ddots &\\
        & & \ddots & B \\
        & & &  -A
        \end{pmatrix} \in {\mathbb C}^{2^pn \times (2^p-1)n}.$$
    The \textit{condition number} of \textbf{IRS} corresponding to the inputs $A,B$, and $p$ is 
    $$ \kappa_{\text{IRS}}(A,B,p) = \sigma_{\min}(D_{(A,B)}^p)^{-1} \left| \left| \binom{A}{B} \right| \right|_2. $$
\end{defn} 
It is not hard to show that $\kappa_{\text{IRS}}(A,B,p)$ is invariant to both swapping $A$ and $B$ and scaling $(A,B)$ and  also satisfies $\kappa_{\text{IRS}}(A,B,p) \geq 1$ for any inputs.\footnote{This follows from the observation $\sigma_{\min}(D_{(A,B)}^p) \leq \sigma_{\min}\binom{A}{B}$.} Moreover, we have the following lemma, which suggests that bounds involving $\kappa_{\text{IRS}}(A,B,p)$ should be sharper than existing results. 
\begin{lem}\label{lem: condition_numbers}
    Let $(A,B)$ be an $n \times n$ regular pencil. Then for any $p \geq 1$ we have
    $$\sigma_{\min}(D_{(A,B)}^p) \geq d_{(A,B)}. $$ 
\end{lem}
\begin{proof}
    Let $m = 2^p$ and define the $mn \times mn$ block matrix 
    \begin{equation}
        M_p(A,B) = \begin{pmatrix} 
        -A & & & -B \\
        B & -A & & \\
        & \ddots & \ddots & \\
        & & B & -A \end{pmatrix} .
    \end{equation}
    To first show $\sigma_{\min}(D_{(A,B)}^p) \geq \sigma_{\min}(M_p(A,B))$, let $x = [x_1 \; x_2 \; \cdots \;  x_{m-1}]^T \in {\mathbb C}^{(m-1)n}$ be the unit vector satisfying $\sigma_{\min}(D_{(A,B)}^p) = \|D_{(A,B)}^px\|_2$, where $x_i \in {\mathbb C}^{n}$ for each $i$. Padding $x$ with zeros to obtain another unit vector 
    \begin{equation}
        y = [x_{m-1} \; \; x_{m-2} \; \cdots \; x_1 \; \;  0] \in {\mathbb C}^{mn}
        \label{eqn: y_unit_vector}
    \end{equation}
    it is easy to see $\|M_p(A,B)y\|_2 = \|D_{(A,B)}^px\|_2$ and therefore 
    \begin{equation}
       \sigma_{\min}(M_p(A,B)) \leq \|M_p(A,B)y\|_2 = \|D_{(A,B)}^px\|_2 = \sigma_{\min}(D_{(A,B)}^p).
        \label{eqn: min_sv_bound}
    \end{equation}
   The first inequality now follows from an observation of Bai, Demmel, and Gu, who show that $M_p(A,B)$ is unitarily equivalent to the block matrix $\text{diag}(-A + e^{i\theta_1}B, \ldots, -A + e^{i\theta_m}B)$ for $e^{i \theta_1}, \ldots, e^{i \theta_m}$ the $m^{\text{th}}$ roots of $-1$. That is, 
   \begin{equation}
        \aligned
        \sigma_{\min}(M_p(A,B)) &= \min_{1 \leq j \leq m} \sigma_n(-A + e^{i\theta_j}B) \\
        &\geq \min_{\theta} \sigma_n(-A + e^{i\theta}B) = d_{(A,B)},
        \endaligned 
        \label{eqn: first condition inequ}
    \end{equation}
    which completes the proof. 
\end{proof}
As we might expect, $\kappa_{\text{IRS}}(A,B,p)$ is not necessarily infinite if $(A,B)$ has an eigenvalue on the unit circle. 
It is also the only condition number posed for \textbf{IRS} that has an explicit dependence on $p$, the number of steps of squaring. While $\kappa_{\text{IRS}}(A,B,p)$ increases with $p$, \cref{lem: condition_numbers} implies the $p$-independent upper bound
\begin{equation}
    \kappa_{\text{IRS}}(A,B,p) \leq d_{(A,B)}^{-1} \left| \left| \binom{A}{B} \right| \right|_2.
    \label{eqn: condition_upper_bound}
\end{equation}
Thinking of $p$ as an input to the procedure not only provides a tighter condition number but also allows us to quantify the stability of \textbf{IRS} in terms of the number of steps taken (and its dependence on $n$). 

\subsection{Floating-Point Stability Bound}
\indent With $\kappa_{\text{IRS}}(A,B,p)$ at our disposal, we now pursue a rigorous and general stability bound for \textbf{IRS} in finite-precision arithmetic. Here, we assume a floating-point (i.e., finite-precision) model of computation, where
\begin{equation}
    fl(x \circ y) = (x \circ y)(1 + \Delta), \; \; |\Delta| \leq {\bf u} 
    \label{eqn: floating point}
\end{equation}
for basic operations $\circ \in \left\{ +, -, \times, \div \right\}$ and a machine precision ${\bf u}(\epsilon, n)$, which is a function of the desired accuracy $\epsilon$ and the size of the problem $n$. This is a standard formulation for finite-precision arithmetic (see for example \cite{Higham}). Given ${\bf u}$, the number of bits of precision required to achieve \eqref{eqn: floating point} is $\log_2(1/{\bf u})$. \\
\indent To analyze \textbf{IRS} in this model, we further assume access to the following black-box algorithms for QR and matrix multiplication, where $\mu_{\text{MM}}(n)$ and $\mu_{\text{QR}}(n)$ are (low-degree) polynomials in $n$. 

\begin{assump}[Matrix Multiplication]
    There exists a $\mu_{\text{MM}}(n)$-stable $n \times n$ multiplication algorithm ${\bf MM}( \cdot, \cdot)$ satisfying 
    $$\|{\bf MM}(A,B) - AB \|_2 \leq \mu_{\text{MM}}(n) {\bf u}\|A\|_2 \|B\|_2 $$
    in $T_{\text{MM}}(n)$ arithmetic operations.
    \label{MM assumption}
\end{assump}

\begin{assump}[QR Factorization]
    There exists a $\mu_{\text{QR}}(n)$-stable full QR algorithm ${\bf QR}(\cdot)$ satisfying
    \begin{enumerate}[itemsep=0em]
        \item $[Q,R] = {\bf QR}(A)$ for $A,R \in {\mathbb C}^{2n \times n}$ and $Q \in {\mathbb C}^{2n \times 2n}$.
        \item $R$ is exactly upper triangular
        \item There exist $A' \in {\mathbb C}^{2n \times n}$ and unitary $Q' \in {\mathbb C}^{2n \times 2n}$ such that $A' = Q'R$,
            $$ \|Q'-Q\|_2 \leq \mu_{\text{QR}}(n) {\bf u}, \; \; \text{and} \; \; \|A' - A\|_2 \leq \mu_{\text{QR}}(n) {\bf u} \|A\|_2 $$
    \end{enumerate}
    in $T_{\text{QR}}(n)$ arithmetic operations.
    \label{QR assumption}
\end{assump}

\indent Once again, these black-box assumptions are standard \cite[Section 3.5 and Chapter 19]{Higham}. While we won't be too particular about $\mu_{\text{MM}}(n)$ and $\mu_{\text{QR}}(n)$, we do note that they are compatible with fast\footnote{i.e., sub-$O(n^3)$} matrix multiplication; that is, QR can be implemented stably (in a mixed sense) using fast matrix multiplication \cite{2007}, which itself can be formulated to satisfy the forward error bound given by \cref{MM assumption} \cite{fastMM}. Consequently our analysis applies to \textbf{IRS} implemented with a variety of fast matrix multiplication routines \cite{STRASSEN1969,COPPERSMITH1990251,fastermm}, including the current fastest known algorithm of Williams et al.\ \cite{williams2023new}, and we may additionally assume $T_{\text{QR}}(n) = O(T_{\text{MM}}(n))$. \\
\indent In terms of these black-box assumptions, each iteration of \textbf{IRS} consists of the following.
\begin{enumerate}
    \item $[Q,R] = \textbf{QR} \left( \begin{bmatrix} B_j \\ -A_j \end{bmatrix} \right)$ with $Q = \begin{pmatrix} Q_{11} & Q_{12} \\ Q_{21} & Q_{22} \end{pmatrix}$
    \item $A_{j+1} = \textbf{MM}(Q_{12}^H, A_j)$
    \item $B_{j+1} = \textbf{MM}(Q_{22}^H, B_j)$
\end{enumerate}
Allowing the error guarantees implied by our black-box assumptions to propagate through multiple iterations yields the following (weak) forward stability bound for the algorithm overall. \\
\indent The proof of \cref{thm: main_IRS_bound} -- which we defer to \cref{section: appendix} -- follows an argument originally developed by Malyshev \cite{Malyshev1}. Our result is more general in that it accommodates the rigorous floating-point assumptions summarized above (and therefore also finite-precision, fast matrix multiplication) and is stated in terms of $\kappa_{\text{IRS}}(A,B,p)$. The latter ensures that the bound is usable in both the spectral projector setting and more general applications.  

\begin{thm}\label{thm: main_IRS_bound}
    Given $A,B \in {\mathbb C}^{n \times n}$ and  $p \geq  1$, let $[\widetilde{A}_p, \widetilde{B}_p] = \text{\normalfont \bf IRS}(A, B, p)$ on a floating-point machine with precision ${\bf u}$. For $\epsilon \in (0,1)$, $\kappa_{\text{\normalfont IRS}}(A,B,p)$ as in \cref{defn: irs_condition_num}, and $\mu(n) = \max \left\{ \mu_{\text{\normalfont MM}}(n), \mu_{\text{\normalfont QR}}(n), \sqrt{n} \right\}$, suppose
    $${\bf u} \leq \frac{\epsilon}{324\mu(n) \kappa_{\text{\normalfont IRS}}(A,B,p) \max \left\{p^2+4p-5,1 \right\}}. $$ 
    Then there exist matrices $\mathring{A}_p, \mathring{B}_p \in {\mathbb C}^{n \times n}$ such that $\mathring{A}_p^{-1} \mathring{B}_p = (A^{-1}B)^{2^p}$ and
    $$ \| \widetilde{A}_p - \mathring{A}_p\|_2, \|\widetilde{B}_p - \mathring{B}_p\|_2 \leq \epsilon \left| \left| \binom{A}{B} \right| \right|_2.$$
\end{thm}

\indent \cref{thm: main_IRS_bound} implies that the number of bits of precision required for \textbf{IRS} to compute $\widetilde{A}_p$ and $\widetilde{B}_p$ to within $\epsilon \| \binom{A}{B}\|_2$ of a corresponding set of exact outputs is at most
\begin{equation}
    \log_2(1/{\bf u}) = O \left( \log_2(1/\epsilon) + \log_2(n) + \log_2(\kappa_{\text{IRS}}(A,B,p)) + \log_2(p) \right).
    \label{eqn: irs_bit_requirement}
\end{equation}
When used to compute projectors as part of the randomized versions of divide-and-conquer developed in \cite{banks2020pseudospectral,arXiv}, this precision requirement is provably lower, in general, than alternatives that require inversion (see the discussion in \cite[Chapter 6]{My_thesis}).

\section{Sign Function Methods}\label{section: SIGN}
We consider next a family of methods based around the matrix sign function. In terms of our high-level framework, $S$ is now the right half plane, in which case
\begin{equation}
    \mathds{1}_S(z) = \frac{1}{2} (\text{sign}(z) +1), 
    \label{eqn: indicator_for_rhp}
\end{equation}
for $\text{sign}(z)$ the scalar sign function
\begin{equation}
    \text{sign}(z) = \begin{cases}
       + 1 & \text{Re}(z) > 0  \\
      -1 & \text{Re}(z) < 0 \\
      \text{undefined} & \text{otherwise.}
    \end{cases}
    \label{eqn: scalar_sign}
\end{equation}
In this setting, approximations of $\mathds{1}_S(z)$ can be derived from approximations of $\text{sign}(z)$, and moreover computing $P_R$ and $P_L$ reduces to approximating (implicitly) the matrix sign function $\text{sign}(B^{-1}A)$, as defined by Roberts \cite{SGN}. Accordingly, eigenvalues are driven to $\pm 1$.
\begin{defn}\label{defn: matrix_sign}
    Let $A \in {\mathbb C}^{n \times n}$ have no eigenvalues on the imaginary axis and suppose 
    $$A = P \begin{pmatrix} J_+ & \\ & J_- \end{pmatrix} P^{-1}$$
    is its Jordan canonical form, with blocks $J_+$ and $J_-$ associated to eigenvalues in the right and left half planes, respectively. Then the matrix sign function of $A$ is
    $$ \text{sign}(A) = P \begin{pmatrix} I & \\ & -I \end{pmatrix} P^{-1}, $$
    where $I$ and $-I$ have the same dimensions as $J_+$ and $J_-$.
\end{defn}

\subsection{Newton Iteration}
The most commonly used method for approximating $\text{sign}(A)$ is a simple Newton iteration of Roberts \cite{SGN}. From the viewpoint of function approximation, this iteration computes $\text{sign}(z)$ via the rational function obtained by repeatedly composing $f(z) = \frac{1}{2}(z+z^{-1})$ with itself. 
\begin{defn}\label{defn: Newton_iteration}
    The \textit{Newton iteration} for computing $\text{sign}(A)$ is given by 
    $$ A_{j+1} = \frac{1}{2}(A_j + A_j^{-1}); \; \; \; \; A_0 = A. $$
\end{defn}
Recalling that the multiplicative ``inverse" of $(B \backslash A)$ is $(A \backslash B)$, the standard Newton iteration can be applied to matrix relations as follows:
\begin{equation}
    (B_{j+1} \backslash A_{j+1}) = \frac{1}{2} \left[(B_j \backslash A_j) + (A_j \backslash B_j) \right]; \; \; \; \; (B_0 \backslash A_0) = (B \backslash A).
    \label{eqn: matrix_relation_newton}
\end{equation}
\cref{alg:IF_Newton} executes $p$ steps of this iteration according to \cref{thm: matrix_relation_comp}.\footnote{This is not the first inverse-free implementation of the Newton iteration. In fact, Benner and Byers present their own version of \textbf{IF-Newton} \cite[Algorithm 1]{Benner_Byers}, which incorporates scaling to promote faster convergence.} Here, the factor of $\frac{1}{2}$ is applied by scaling $A_{j+1}$ by $\frac{1}{\sqrt{2}}$ and $B_{j+1}$ by $\sqrt{2}$, which is necessary to guarantee convergence of the individual matrices as $j \rightarrow \infty$ in exact arithmetic (see \cite[Theorem 3.6]{Benner_Byers}). As in the approach based on \textbf{IRS}, some post-processing is necessary to obtain $P_{R, \text{Re}(z)>0}$, where again the subscript clarifies the corresponding subset of the spectrum of $(A,B)$. In this case, $B_j^{-1}A_j$ approximates $\text{sign}(B^{-1}A)$ and therefore
\begin{equation}
    P_{R,\text{Re}(z)>0} \approx \frac{1}{2} (B_j^{-1}A_j + I) = \frac{1}{2}B_j^{-1}(A_j + B_j).
    \label{eqn: P_R_from_Newton}
\end{equation} 
Equivalently, $P_{R,\text{Re}(z)>0}$ corresponds to the matrix relation $(2B_j \backslash (A_j+B_j))$. \\
\indent As its name suggests, the Newton iteration can be viewed as an extension of classical Newton's method, which finds roots of $z^2-1$ according to $z_{j+1} = \frac{1}{2}(z_j+z_j^{-1})$. Indeed, the Newton iteration for $\text{sign}(A)$ applies this version of Newton's method to the eigenvalues of $A$, and quadratic convergence of classical Newton's method implies quadratic convergence for \eqref{eqn: matrix_relation_newton}. 

\renewcommand{\arraystretch}{1}
\begin{algorithm}[t]
\caption[Inverse-Free Newton Iteration (\textbf{IF-Newton})]{Inverse-Free Newton Iteration (\textbf{IF-Newton})\\
\textbf{Input:} $A, B \in {\mathbb C}^{n \times n}$, $p$ a number of iterations. \\
\textbf{Requires:} $(A,B)$ has no eigenvalues on $\text{Re}(z) = 0$.}\label{alg:IF_Newton}
\begin{algorithmic}[1]
\State $A_0 = A$
\State $B_0 = B$
\For{$j=0:p-1$}
    \vspace{1mm}
    \State $\begin{pmatrix} -A_j \\ B_j 
    \end{pmatrix}  = \begin{pmatrix} Q_{11} & Q_{12} \\
    Q_{21} & Q_{22} \end{pmatrix} \begin{pmatrix} R_j \\ 0  \end{pmatrix}$
    \vspace{1mm}
    \State $A_{j+1} = \frac{1}{\sqrt{2}} (Q_{12}^HB_j + Q_{22}^HA_j)$
    \State $B_{j+1} = \sqrt{2} Q_{22}^HB_j$
\EndFor
\State \Return $(A_p, B_p)$
\end{algorithmic}
\end{algorithm}

\subsection{(Weighted) Halley Iteration}
\indent The connection between \textbf{IF-Newton} and classical Newton's method suggests that alternative algorithms for computing $\text{sign}(A)$ can be obtained from other root finding iterations. Halley's method, for example, approximates roots of $z^2-1$ according to the third-order iteration
\begin{equation}
    z_{j+1} = z_j\frac{z_j^2+3}{3z_j^2+1}.
    \label{eqn: Halley_scalar_iteration}
\end{equation}
Consequently, it implies the following iteration for $\text{sign}(A)$. 
\begin{defn}\label{defn: Halley_iteration}
    The \textit{Halley iteration} for computing $\text{sign}(A)$ is given by 
    $$A_{j+1} = A_j(3A_j^2+1)^{-1}(A_j^2+3); \; \; \; \; A_0 = A.$$
\end{defn}
Recalling that any Möbius transformation can be applied to $(B \backslash A)$ for free, only two QR factorizations are required to run the Halley iteration on matrix relations if evaluated as 
\begin{equation}
    (B_{j+1} \backslash A_{j+1}) = (B_j \backslash A_j) h((B_j \backslash A_j)^2); \; \; \; \; (B_0 \backslash A_0) = (B \backslash A)
    \label{eqn: Halley_for_relation}
\end{equation}
for  $h(z) = \frac{3z+1}{z+3}$. As in the Newton iteration, the approximation of $\text{sign}(z)$ corresponding to \eqref{eqn: Halley_for_relation} can be obtained by repeated composition, this time with $f(z) = zh(z^2)$. Applying \cref{thm: matrix_relation_comp} yields \cref{alg:IF_Halley}, which executes $p$ steps of this Halley iteration on an arbitrary pencil $(A,B)$. As in \textbf{IF-Newton}, the outputs of this routine yield the projector $P_{R,\text{Re}(z)>0}$ according to \eqref{eqn: P_R_from_Newton}. 

\begin{algorithm}
\caption[Inverse-Free Halley Iteration (\textbf{IF-Halley})]{Inverse-Free Halley Iteration (\textbf{IF-Halley})\\
\textbf{Input:} $A, B \in {\mathbb C}^{n \times n}$, $p$ a number of iterations. \\
\textbf{Requires:} $(A,B)$ has no eigenvalues on $\text{Re}(z) = 0$.}\label{alg:IF_Halley}
\begin{algorithmic}[1]
\State $A_0 = A$
\State $B_0 = B$
\For{$i = 0:p-1$} 
    \vspace{1mm}
    \State $\begin{pmatrix} -B_i \\ A_i
    \end{pmatrix}  = \begin{pmatrix} Q_{11} & Q_{12} \\
    Q_{21} & Q_{22} \end{pmatrix} \begin{pmatrix} R_i \\ 0  \end{pmatrix}$
    \vspace{1mm}
    \State $C_i = Q_{12}^HA_i + 3 Q_{22}^HB_i$
    \State $D_i = 3Q_{12}^HA_i + Q_{22}^HB_i$
    \vspace{1mm}
    \State $\begin{pmatrix} -D_i \\ A_i \end{pmatrix} = \begin{pmatrix} U_{11} & U_{12} \\ U_{21} & U_{22} \end{pmatrix} \begin{pmatrix} \widehat{R}_i \\ 0 \end{pmatrix} $
    \vspace{1mm}
    \State $A_{i+1} = U_{12}^HC_i$
    \State $B_{i+1} = U_{22}^HB_i$
\EndFor
\State \Return $(A_p,B_p)$
\end{algorithmic}
\end{algorithm}

\indent If \textbf{IF-Halley} can exhibit third-order convergence, as we demonstrate more rigorously in the next section, is it necessarily a better choice than \textbf{IF-Newton}? In practice, the answer is no; in fact, it is likely to be less efficient in general. As a third-order method, we can expect \textbf{IF-Halley} to cut the number of iterations required by \textbf{IF-Newton} by roughly a factor of $\log_2(3)$ (in general $\log_2(m)$ for an order $m$ iteration). Since each iteration of \textbf{IF-Halley} is twice as expensive as one of \textbf{IF-Newton}, the latter is likely to be less costly overall, despite converging more slowly. \\
\indent  Nevertheless, we might hope that by modifying the Halley iteration we can overcome this drawback, guaranteeing faster convergence on at least some problems. To that end, we consider varying the Möbius transformation $h$, replacing \eqref{eqn: Halley_scalar_iteration} with 
\begin{equation}
    z_{j+1} = z_j \frac{a_jz_j^2 + b_j}{c_jz_j^2 + d_j}
    \label{eqn: weighted_scalar_Halley}
\end{equation}
for some $a_j,b_j,c_j,d_j \in {\mathbb C}$ satisfying $a_jd_j- b_jc_j \neq 0$, which are allowed to evolve with each iteration. Note that doing so will not change the complexity of the method. To streamline the choice of these coefficients, we make the following simplifications:
\begin{enumerate}
    \item Without loss of generality, we can set $d_j = 1$ for all $j$. 
    \item Since $\text{sign}(z)$ fixes $\pm 1$, we enforce that the iteration does as well. This will be guaranteed as long as $c_j = a_j+b_j-1$.
\end{enumerate}
From here we obtain a specialized Indicator Approximation Problem:\ if we know that the spectrum of $(A,B)$ is contained in particular subsets of $S$ and ${\mathbb C} \setminus S$, can we choose $a_j$ and $b_j$ to optimally approximate $\text{sign}(z)$ on them? \\
\indent When the spectrum of $(A,B)$ is real, and in particular lies in a union of intervals $[-1,-l_0) \cup (l_0, 1]$ for some $l_0 > 0$,\footnote{We assume here that zero is not an eigenvalue of $(A,B)$, as in that case the sign function is not defined.} a solution to this problem is known. By considering only real $a_j,b_j$ in this case, we guarantee that eigenvalues lie in a similar union $[-1,-l_j) \cup (l_j, 1]$ at each step, where
\begin{equation}
    l_{j+1} = \min_{l_j \leq x \leq 1} x \frac{a_jx^2 + b_j}{(a_j + b_j - 1)x^2 + 1}.
    \label{eqn: Halley_progress}
\end{equation}
To promote fast convergence, we want $l_j$ to be as close to one as possible at each iteration (see \cref{lem: shifted_dwh_bound}). We therefore obtain the following optimization problem.
\begin{equation}
    \maximize_{a_j,b_j} \; l_{j+1} \; \; \; \; \text{subject to} \;  \;  \; \; a_j,b_j > 0 \; \; \text{and} \;  \; a_j+b_j > 1.
    \label{eqn: Halley_optimization}
\end{equation}
Here, requiring $a_j,b_j,c_j > 0$ ensures that $(l_j, 1] \mapsto (l_{j+1}, 1]$ guarantees also $[-1,-l_j) \mapsto [-1, -l_{j+1})$. \\
\indent For any starting value $l_0$, \eqref{eqn: Halley_optimization} admits the following solution: 
\begin{equation}
    \begin{cases}
        \gamma_j = \sqrt[3]{\frac{4 (1 - l_j^2)}{l_j^4}}; \; \; \; \; \; b_j = \sqrt{1 + \gamma_j} + \frac{1}{2} \sqrt{8 - 4 \gamma_j + \frac{8(2 - l_j^2)}{l_j^2 \sqrt{1 + \gamma_j}}}; \\
        a_j = \frac{1}{4}(b_j - 1)^2; \; \; \; l_{j+1} = l_j \frac{a_j l_j^2 + b_j}{(a_j+b_j-1)l_j^2 + 1}
   \end{cases}
    \label{eqn: weighted_coefficients}
\end{equation}
For this choice of $a_j$ and $b_j$ we note that $(a_j,b_j) \rightarrow (1,3)$ as $l_j \rightarrow 1$, meaning the corresponding modified Halley iteration gradually approaches the standard version as it converges. \\
\indent The solution \eqref{eqn: weighted_coefficients} is due to  Nakatsukasa, Bai, and Gygi \cite{optimizing_halley}, who introduced a \textit{dynamically weighted} Halley iteration to compute the polar decomposition of a matrix. In their case, a variation of \cref{defn: Halley_iteration} converges to the (unitary) polar factor of $A$ by driving its singular values to one. As a result, they arrive at the same optimization problem, seeking a Möbius transformation that yields a particularly accurate approximation of  $\text{sign}(z)$ on a portion of the real axis (in their case just $(l_j,1]$). While the optimized coefficients were originally obtained via a direct and exhaustive search, a connection to the work of Zolotarev \cite{Zolo} was later made by Nakatsukasa and Freund \cite{zolo_dnc}, who demonstrated that the rational function corresponding to \eqref{eqn: weighted_coefficients} can be interpreted as an optimal approximation to the sign function on $[-1, -l_j] \cup [l_j, 1]$ -- i.e., it solves exactly our specialized Indicator Approximation Problem.\footnote{Optimal meaning the best (in the infinity norm) rational function approximation $p(x)/q(x)$ for $p(x)$ and $q(x)$ real polynomials of degree three and two, respectively.} \\
\indent Applying \eqref{eqn: weighted_coefficients} to \cref{alg:IF_Halley} produces \cref{alg:IF_DWH}, which -- borrowing  terminology from Nakatsukasa, Bai, and Gygi -- we call an inverse-free dynamically weighted Halley iteration (\textbf{IF-DWH}). Note that this routine requires not only that the spectrum of $(A,B)$ is contained in a symmetric union of intervals in $[-1,1]$ but that a lower bound $l_0$ on the minimum eigenvalue (in magnitude) is known.

\begin{algorithm}[h]
\caption[Inverse-Free Weighted Halley Iteration (\textbf{IF-DWH})]{Inverse-Free Weighted Halley Iteration (\textbf{IF-DWH})\\
\textbf{Input:} $A, B \in {\mathbb C}^{n \times n}$, $p$ a number of iterations, $l_0 > 0$. \\
\textbf{Requires:} Eigenvalues $\lambda$ of $(A,B)$ are real with $ l_0 < |\lambda|  \leq 1$.}\label{alg:IF_DWH}
\begin{algorithmic}[1]
\State $A_0 = A$
\State $B_0 = B$
\For{$j = 0:p-1$} 
    \State $\gamma_j = \left(4 (1 - l_j^2)/l_j^4 \right)^{1/3} $
    \State $b_j = \sqrt{1 + \gamma_j} + \frac{1}{2} \sqrt{8 - 4 \gamma_j + 8(2 - l_j^2)/(l_j^2 \sqrt{1 + \gamma_j})}$
    \State $a_j = \frac{1}{4} (b_j - 1)^2$
    \State $c_j = a_j + b_j - 1$
    \vspace{1mm}
    \State $\begin{pmatrix} -B_j \\ A_j
    \end{pmatrix}  = \begin{pmatrix} Q_{11} & Q_{12} \\
    Q_{21} & Q_{22} \end{pmatrix} \begin{pmatrix} R_j \\ 0  \end{pmatrix}$ \Comment{Apply Halley iteration}
    \vspace{1mm}
    \State $C_j = a_jQ_{12}^HA_j + b_jQ_{22}^HB_j$
    \State $D_j = c_iQ_{12}^HA_j + Q_{22}^HB_j$
    \vspace{1mm}
    \State $\begin{pmatrix} -D_j \\ A_j \end{pmatrix} = \begin{pmatrix} U_{11} & U_{12} \\ U_{21} & U_{22} \end{pmatrix} \begin{pmatrix} \widehat{R}_j \\ 0 \end{pmatrix} $
    \vspace{1mm}
    \State $A_{j+1} = U_{12}^HC_j$
    \State $B_{j+1} = U_{22}^HB_j$
    \State $l_{j+1} = l_j(a_jl_j^2 + b_j)/(c_jl_j^2 + 1)$ \Comment{Compute next value of $l$}
\EndFor
\State \Return $(A_p,B_p)$, optionally $l_p$
\end{algorithmic}
\end{algorithm}

\indent \textbf{IF-DWH} is relevant for any regular matrix pencil that has (nonzero) real eigenvalues. This includes the important definite generalized eigenvalue problem, where $A$ and $B$ are Hermitian and the \textit{Crawford number}
\begin{equation}\label{eqn: Crawford}
    \gamma(A,B) = \min_{\|x\|_2=1}|x^H(A + iB)x| = \min_{\|x\|_2 = 1}\sqrt{(x^HAx)^2 + (x^HBx)^2}
\end{equation}
is strictly positive. As a generalization of the Hermitian eigenvalue problem, definite pencils appear frequently in applications -- see e.g., \cite{SVM, FORD1974337}. In the companion paper \cite{structure_dnc}, we incorporate \textbf{IF-DWH} in a specialized divide-and-conquer algorithm for definite pencils. 

\subsection{Convergence Bounds}
\begin{figure}[t]
    \centering
    \includegraphics[width=.9\linewidth]{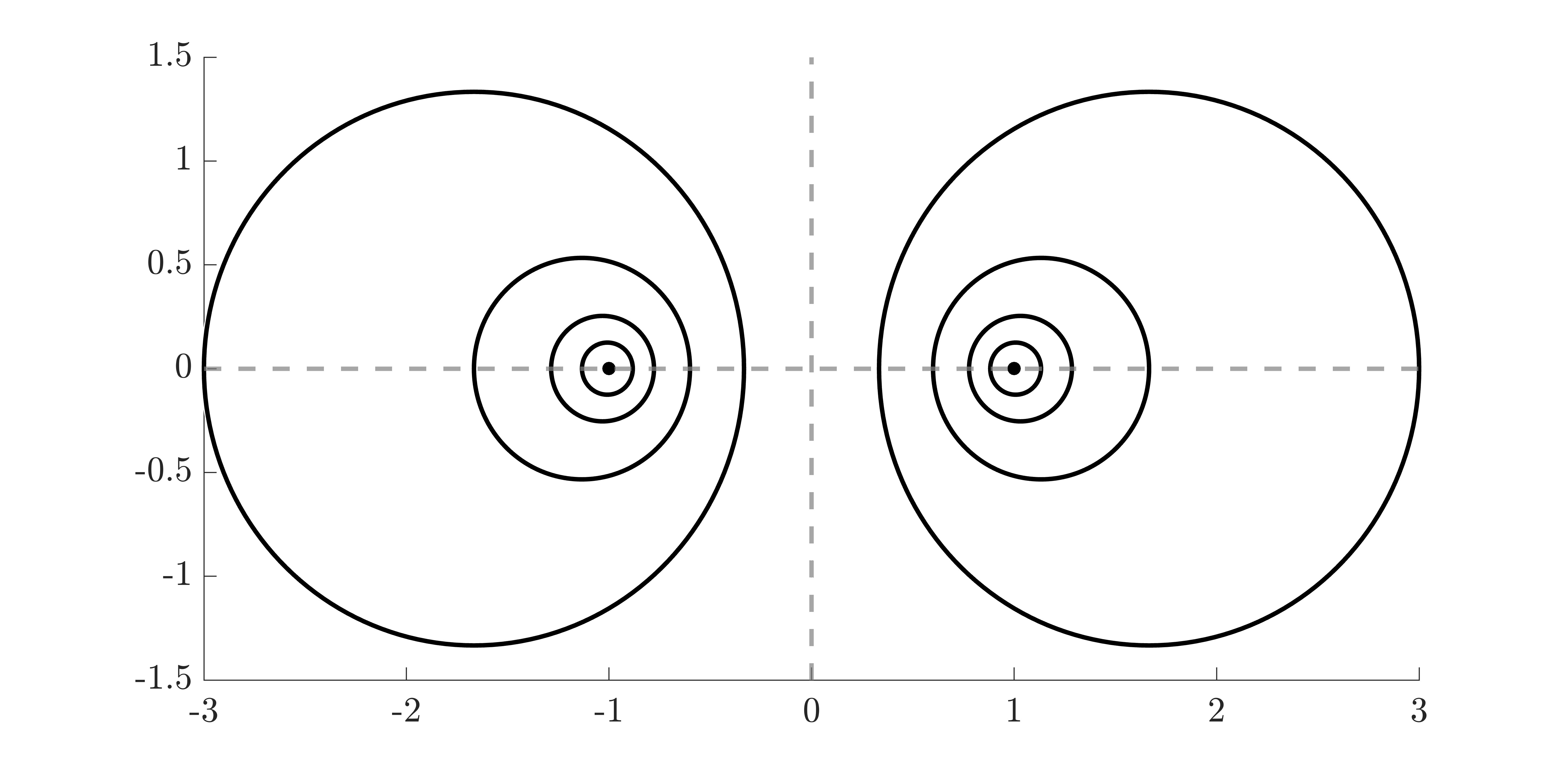}
    \caption{The circles of Apollonius corresponding to $\alpha = \frac{1}{2}, \frac{1}{4}, \frac{1}{8}$, and $\frac{1}{16}$.}
    \label{fig: circles}
\end{figure}

We close this section with a handful of convergence results. \cref{subsub: convergence_1} covers \textbf{IF-Newton}/\textbf{IF-Halley} and their generalizations, while \cref{subsub: convergence_2} presents a new analysis of the dynamically weighted Halley iteration. The latter improves on the work of Nakatsukasa, Bai, and Gygi \cite{optimizing_halley} and is the main contribution of this section. \\
\indent As in much of the literature for sign-function-based methods, the circles of Apollonius are the key theoretical tool here. 
\begin{defn}\label{defn: circles_of_Apollonius}
    For $\alpha \in (0,1)$ let
    $$ C_{\alpha}^+ = \left\{ z : \left| \frac{1-z}{1+z} \right| \leq \alpha \right\}, \; \; \; \; C_{\alpha}^- = \left\{ z :  \left| \frac{1+z}{1-z} \right|  \leq \alpha \right\} $$
    be sets in the right and left half planes, respectively. The boundaries $\partial C_{\alpha}^+$ and $\partial C_{\alpha}^-$ of these sets are the \textit{circles of Apollonius} corresponding to $\alpha$.
\end{defn}
\indent $C_{\alpha}^+$ can be equivalently characterized as the disk with center $\frac{1+\alpha^2}{1-\alpha^2}$ and radius $\frac{2\alpha}{1-\alpha^2}$, with $C_{\alpha}^-$ its image under a reflection across the imaginary axis. For varying $\alpha$, $\partial C_{\alpha}^+$ and $\partial C_{\alpha}^-$ define families of non-concentric circles, which collapse to the points $\pm 1$ as $\alpha \rightarrow 0$. Since this geometric picture will be important to have in mind, \cref{fig: circles} plots a handful of Apollonian circles. Throughout, we use $C_{\alpha}$ to denote the region $C_{\alpha}^+ \cup C_{\alpha}^-$. 

\subsubsection{Bounds for \text{\normalfont \bf IF-Newton} and \text{\normalfont \bf IF-Halley}}\label{subsub: convergence_1}
\indent Given their relationship to the points $\pm 1$, the circles of Apollonius are naturally equipped to describe convergence to the sign function. Indeed, the Newton iteration can be characterized by the following observation of Roberts \cite{SGN}.
\begin{prop}\label{prop: newton_cirlces}
    The function $f(z) = \frac{1}{2}(z+z^{-1})$ defining the Newton iteration maps $C_{\alpha}^{+}$ to $C_{\alpha^2}^{+}$ and $C_{\alpha}^{-}$ to $C_{\alpha^2}^{-}$.
\end{prop}
\indent Extending this to the Halley iteration is straightforward. \cref{lem: Halley_circle} captures the third-order convergence of  Halley's method for finding the roots of $z^2-1$. 
\begin{lem} \label{lem: Halley_circle}
The function $f(z) = zh(z^2) = z \frac{z^2+3}{3z^2+1}$ defining the Halley iteration maps $C_{\alpha}^{+}$ to $C_{\alpha^3}^+$ and  $C_{\alpha}^-$ to $C_{\alpha^3}^{-}$.
\end{lem}
\begin{proof}
    Applying the definition of $C_{\alpha}^{\pm}$, we have  
    \begin{equation}
        \frac{1 - f(z)}{1 + f(z)} = \frac{1 - \frac{z^3 + 3z}{3z^2 + 1}}{1 + \frac{z^3 + 3z}{3z^2 + 1}} = \frac{3z^2 + 1 - z^3 - 3z}{3z^2 + 1 + z^3 + 3z} = \frac{(1-z)^3}{(1+z)^3}.
        \label{eqn: Halley_circle_proof}
    \end{equation}
    The result follows immediately.
\end{proof}

\begin{remark}\label{rem: more_methods}
The proof of \cref{lem: Halley_circle} implies yet another strategy for deriving iterative methods for the sign function:\ work backwards from $\pm \frac{(1-z)^m}{(1+z)^m}$ given a desired order of convergence $m$. As an example, $-\frac{(1-z)^4}{(1+z)^4}$ can be written as 
\begin{equation}
    - \frac{(1-z)^4}{(1+z)^4} = \frac{1 - \frac{1 + 6z^2 + z^4}{4z + 4z^3}}{1 + \frac{1 + 6z^2 + z^4}{4z + 4z^3}}, 
    \label{eqn: fourth_order}
\end{equation}
which implies an iterative rational function
\begin{equation}
    f(z) = \frac{1 + 6z^2 + z^4}{4z + 4z^3} = \frac{z}{4} \left( \frac{1}{z^2} + \frac{z^2 + 5}{1 + z^2} \right),
    \label{eqn: fourth_order_function}
\end{equation}
where the latter expression indicates that this iteration can be implemented according to \cref{thm: matrix_relation_comp} with only three QR factorizations. 
\end{remark}


\indent We can translate results like \cref{lem: Halley_circle} into bounds for corresponding iterative methods with one more tool:\ the $\epsilon$-pseudospectrum. 

\begin{defn} \label{defn: matrix pseudospectrum}
For any $\epsilon > 0$, the \textit{$\epsilon$-pseudospectrum} of $A$ is
$$ \Lambda_{\epsilon}(A) = \left\{ z : \text{there exists a vector}\; u \neq 0 \; \;\text{s.t.} \;  (A + E)u = zu \; \text{for} \; \|E\|_2 \leq \epsilon \right\}. $$
\end{defn}

The following lemma of Banks et al. \cite[Lemma 4.3]{banks2020pseudospectral} implies that a bound on the pseudospectrum of a matrix, in terms of the circles of Apollonius, can be bootstrapped into a bound on the sign function. 
\begin{lem}[Banks et al.\ \cite{banks2020pseudospectral}]\label{lem: pseudospectral_convergence}
Suppose $\Lambda_{\epsilon}(A) \subset C_{\alpha}$ for some $\epsilon > 0$. Then,
$$ \|A - \text{\normalfont sign}(A) \|_2 \leq \frac{8 \alpha^2}{\epsilon(1+\alpha)(1-\alpha)^2}. $$
\end{lem}

To make use of this result, we need only characterize the way our iterations transform pseudospectral bounds. The following lemma provides a general picture (and in fact specifically generalizes another result of Banks et al.\ \cite[Lemma 4.4]{banks2020pseudospectral}).

\begin{lem}\label{lem: pseudo_evolution}
    Suppose the rational function $f$ has all of its poles on the imaginary axis and maps $C_{\alpha}^{\pm} \rightarrow C_{\alpha^m}^{\pm}$ with $ \partial C_{\alpha}^{\pm} \rightarrow \partial C_{\alpha^m}^{\pm}$ for any $\alpha \in (0,1)$. Let $f$ define an iteration for $\text{\normalfont sign}(X)$ according to 
    $$ X_{k+1} = f(X_k); \; \; \; \; X_0 = X.$$ 
    If $\Lambda_{\epsilon}(X_k) \subset C_{\alpha}$, then for any $\alpha' \in (\alpha^m, \alpha)$ we have $\Lambda_{\epsilon'}(X_{k+1}) \subset C_{\alpha'}$ with
    $$ \epsilon' = \frac{\epsilon (1 - \alpha^2)(\alpha'-\alpha^m)}{8 \alpha} .$$
\end{lem}
\begin{proof}
    Let $w$ be any point in the ``annulus" between $C_{\alpha}$ and $C_{\alpha'}$. Since $f$ maps $C_{\alpha}$ to $C_{\alpha'}$ and $w \notin C_{\alpha'}$, the rational function $\frac{1}{w - f(z)}$ is holomorphic on $C_{\alpha}$. Moreover, $C_{\alpha}$ contains $\Lambda(X_k)$, meaning we can bound $\|(wI - X_{k+1})^{-1}\|_2$ as 
    \begin{equation}
        \aligned 
        \|(wI - &X_{k+1})^{-1}\|_2 = \left| \left| \frac{1}{2 \pi i}\int_{\partial C_{\alpha}} \frac{(w - f(z))^{-1}}{zI - X_k} dz \right| \right|_2 \\
        &\leq \frac{1}{2 \pi} \int_{\partial C_{\alpha}^+} \frac{\| (zI - X_k)^{-1}\|_2}{|w - f(z)|} dz + \frac{1}{2 \pi} \int_{\partial C_{\alpha}^-} \frac{\|(zI - X_k)^{-1}\|_2}{|w -f(z)|} dz .
        \endaligned 
        \label{eqn: triangle_hfc}
    \end{equation}
   The first integral in this sum can be bounded by the length of the contour $\partial C_{\alpha}^+$ times the supremum of the integrand -- i.e., 
    \begin{equation}
        \int_{\partial C_{\alpha}^+} \frac{\|(zI-X_k)^{-1}\|_2}{|w-f(z)|} dz \leq \frac{4 \pi \alpha}{1 - \alpha^2} \sup_{z \in \partial C_{\alpha}^+} \frac{\|(zI-X_k)^{-1}\|_2}{|w-f(z)|}.
        \label{eqn: apply_ML_2}
    \end{equation}
    Now $\Lambda_{\epsilon}(X_k) \cap \partial C_{\alpha}^+ = \emptyset$, so $\|(zI - X_k)^{-1}\|_2 \leq \epsilon^{-1}$ for all $z \in \partial C_{\alpha}^+$. Using the fact that $f(z) \in C_{\alpha^m}^+$ if $z \in C_{\alpha}^+$, we therefore have
    \begin{equation}
        \aligned
        \int_{\partial C_{\alpha}^+} \frac{\|(zI-X_k)^{-1}\|_2}{|w-f(z)|} dz &\leq \frac{4 \pi \alpha}{\epsilon(1 - \alpha^2)} \sup_{y \in \partial C_{\alpha^m}^+} \frac{1}{|w-y|} \\
        &\leq \frac{8 \pi \alpha}{\epsilon (1 - \alpha^2)(\alpha' - \alpha^m)},
        \endaligned 
        \label{eqn: simplify_ML}
    \end{equation}
    where the last inequality follows from \cite[Lemma 4.5]{banks2020pseudospectral}. Since we obtain the same bound on the remaining term of \eqref{eqn: triangle_hfc}, we conclude 
    \begin{equation}
        \| (wI - X_{k+1}])^{-1}\|_2 \leq \frac{8 \alpha}{\epsilon (1 - \alpha^2)(\alpha' - \alpha^m)},
        \label{eqn: final_w_bound}
    \end{equation}
   and therefore $w \notin \Lambda_{\epsilon'}(X_{k+1})$ for $\epsilon' = \frac{\epsilon (1 - \alpha^2)(\alpha'-\alpha^m)}{8 \alpha}$. Since \eqref{eqn: final_w_bound} applies to any point $w$ between $C_{\alpha}$ and $C_{\alpha'}$ and $\Lambda(X_{k+1}) \subset C_{\alpha'}$, this suffices to show $\Lambda_{\epsilon'}(X_{k+1}) \subset C_{\alpha'}$.
 \end{proof}

 In exact arithmetic, the preceding lemmas immediately extend to the setting where $\text{sign}(B^{-1}A)$ is computed implicitly. Taking the place of $\Lambda_{\epsilon}(A)$ in the resulting error bound is the $\epsilon$-pseudospectrum of the pencil $(A,B)$, defined as follows.\footnote{There are actually several ways to define the pseudospectra of $(A,B)$.  \cref{defn:pencil_pseudospectrum} is originally due to Frayssé et al.\ \cite{Fraysse96spectralportraits}.}
 
  \begin{defn}\label{defn:pencil_pseudospectrum}
    For any $\epsilon > 0$, the $\epsilon$-pseudospectrum of $(A,B)$ is
    $$ \Lambda _{\epsilon }(A,B) =  \left\{ z: \begin{array}{c} (A + \Delta A)u = z(B + \Delta B)u \\ \text{\normalfont for} \; u \neq 0 \; \text{\normalfont and} \; \|\Delta A\|_2, \|\Delta B\|_2 \le \epsilon \\ \end{array} \right\}  .$$
\end{defn}

\cref{prop: general_error_bound} presents our main convergence result, which applies to both \textbf{IF-Newton} and \textbf{IF-Halley} (with $m = 2$ and $m = 3$, respectively) but is not specific to either. We draw a connection here to analogous results for \textbf{IRS} -- i.e., \cite[Theorem 1]{Bai:CSD-94-793} -- noting that a bound on $\|B_j^{-1}A_j - \text{sign}(B^{-1}A)\|_2$ can be bootstrapped into one for the projector $P_{R,\text{Re}(z)>0}$ as follows:

\begin{equation}\label{eqn: projctor_bound_from_sign}
  \left| \left| \frac{1}{2} B_j^{-1}(A_j+B_j) - P_{R,\text{Re}(z)>0} \right| \right|_2 = \frac{1}{2} \| B_j^{-1}A_j - \text{sign}(B^{-1}A) \|_2.
\end{equation}

\begin{prop}\label{prop: general_error_bound}
    Let $A,B \in {\mathbb C}^{n \times n}$ and let $f$ be a rational function satisfying the assumptions of \cref{lem: pseudo_evolution} for a corresponding value $m > 1$. Define the following inverse-free iteration for approximating $\text{\normalfont sign}(B^{-1}A)$:
    $$(B_{j+1} \backslash A_{j+1}) = f(B_j \backslash A_j); \; \; \;  (B_0 \backslash A_0) = (B \backslash A).$$
    If $\Lambda_{\epsilon}(A,B) \subset C_{\alpha}$ for some $\epsilon > 0$ and $\alpha \in (0,1)$, then for any $1 < c < \alpha^{-(m-1)}$ we have
    $$ \left| \left| B_j^{-1}A_j - \text{\normalfont sign}(B^{-1}A) \right| \right|_2 \leq  
    \left( c^{\frac{1}{m-1}} \alpha \right)^{m^j} \cdot  \frac{8 \alpha \|B\|_2}{\epsilon (c-1)^2} \cdot \left[ \frac{8c}{(1-\alpha^2)(c-1)} \right]^j.  $$
\end{prop}
\begin{proof}
    Since $\Lambda_{\epsilon}(A,B)$ is bounded $B$ is invertible.\footnote{In fact, $\Lambda_{\epsilon}(A,B)$ is bounded if and only if $\sigma_n(B) > \epsilon$.} Consider then the explicit iteration $X_{j+1} = f(X_j)$ for $X_0 = B^{-1}A$, which in exact arithmetic is equivalent to its inverse-free counterpart. Let $\alpha_j$ and $\epsilon_j$ be (decreasing) sequences of parameters defined recursively as follows:
    \begin{equation}\label{eqn: alpha_j}
        \begin{cases}
        \alpha_j = c \alpha_{j-1}^m,  & \alpha_0 = \alpha \\
        \epsilon_j = \frac{1}{8} \epsilon_{j-1} (1-\alpha^2)(c-1) \alpha_{j-1}^{m-1}, & \epsilon_0 = \frac{\epsilon}{\|B\|_2}.
        \end{cases}
    \end{equation}
    We show inductively that $\Lambda_{\epsilon_j}(X_j) \subset C_{\alpha_j}$. The base case $(j=0)$ is trivial and follows from the observation $\Lambda_{\epsilon/\|B\|_2}(B^{-1}A) \subseteq \Lambda_{\epsilon}(A,B)$. Consider now arbitrary $j$. Plugging our induction hypothesis into \cref{lem: pseudo_evolution}, and taking\footnote{Note that the restriction $c < \alpha^{-(m-1)}$ guarantees $\alpha_j \in (\alpha_{j-1}^m, \alpha_{j-1})$ for any $j$.} $\alpha' = \alpha_j$, we conclude that $\Lambda_{\epsilon_{j-1}}(X_{j-1}) \subset C_{\alpha_{j-1}}$ implies $\Lambda_{\epsilon'}(X_j) \subset C_{\alpha_j}$ for 
    \begin{equation} \label{eqn: inductive_epsilon}
        \aligned 
        \epsilon' = \frac{\epsilon_{j-1} (1 - \alpha_{j-1}^2)(\alpha_j - \alpha_{j-1}^m)}{8 \alpha_{j-1}} > \frac{\epsilon_{j-1} (1 - \alpha^2) (c-1)\alpha_{j-1}^{m-1}}{8}  = \epsilon_j.
        \endaligned 
    \end{equation}
    Hence, $\Lambda_{\epsilon_j}(X_j) \subset C_{\alpha_j}$. We can now bound error in the approximation via \cref{lem: pseudospectral_convergence}. Noting that $\text{sign}(B^{-1}A) = \text{sign}(X_j)$ since $X_j$ and $B^{-1}A$ have the same (right) eigenvectors, we have
    \begin{equation}\label{eqn: error_at_step_j}
        \aligned
        \|B_j^{-1}A_j - \text{sign}(B^{-1}A)\|_2 &= \|X_j - \text{sign}(X_j)\|_2 \\
        &\leq \frac{8 \alpha_j^2}{\epsilon_j (1 + \alpha_j)(1-\alpha_j)^2} \leq \frac{8 \alpha_j^2}{\epsilon_j (c-1)^2},
        \endaligned
    \end{equation}
    where the last inequality follows from the bounds $1 + \alpha_j > 1$ and $|1 - \alpha_j| > c-1$. We complete the proof by converting \eqref{eqn: alpha_j} into the non-recursive expressions $\alpha_j = c^{\frac{m^j - 1}{m-1}} \alpha^{m^j}$ and $\epsilon_j = \frac{\epsilon \alpha_j}{\|B\|_2 \alpha} \left[ \frac{(1-\alpha^2)(c-1)}{8c} \right]^j$ and applying them to \eqref{eqn: error_at_step_j}. To simplify, note that 
    \begin{equation} \label{eqn: simplify_alpha_j}
     c^{\frac{m^j-1}{m-1}} \alpha^{m^j} = \left[ c^{\frac{m^j - 1}{m^j(m-1)}} \alpha \right]^{m^j} \leq \left[ c^{\frac{1}{m-1}} \alpha \right]^{m^j}
    \end{equation}
    for any $j$.
\end{proof} 

In the setting where $B = I$ and $f$ is the rational function defining the standard Newton iteration (so that $m = 2$) this result reduces to \cite[Proposition 4.8]{banks2020pseudospectral}. As mentioned in its proof, \cref{prop: general_error_bound} implicitly assumes that $B$ is invertible via the pseudospectral bound $\Lambda_{\epsilon}(A,B) \subset C_{\alpha}$, though this is of course necessary for $\text{sign}(B^{-1}A)$ to be defined.

\subsubsection{Bounds for \text{\normalfont \bf IF-DWH}} \label{subsub: convergence_2}
\begin{figure}
    \centering
    \includegraphics[width=.8\linewidth]{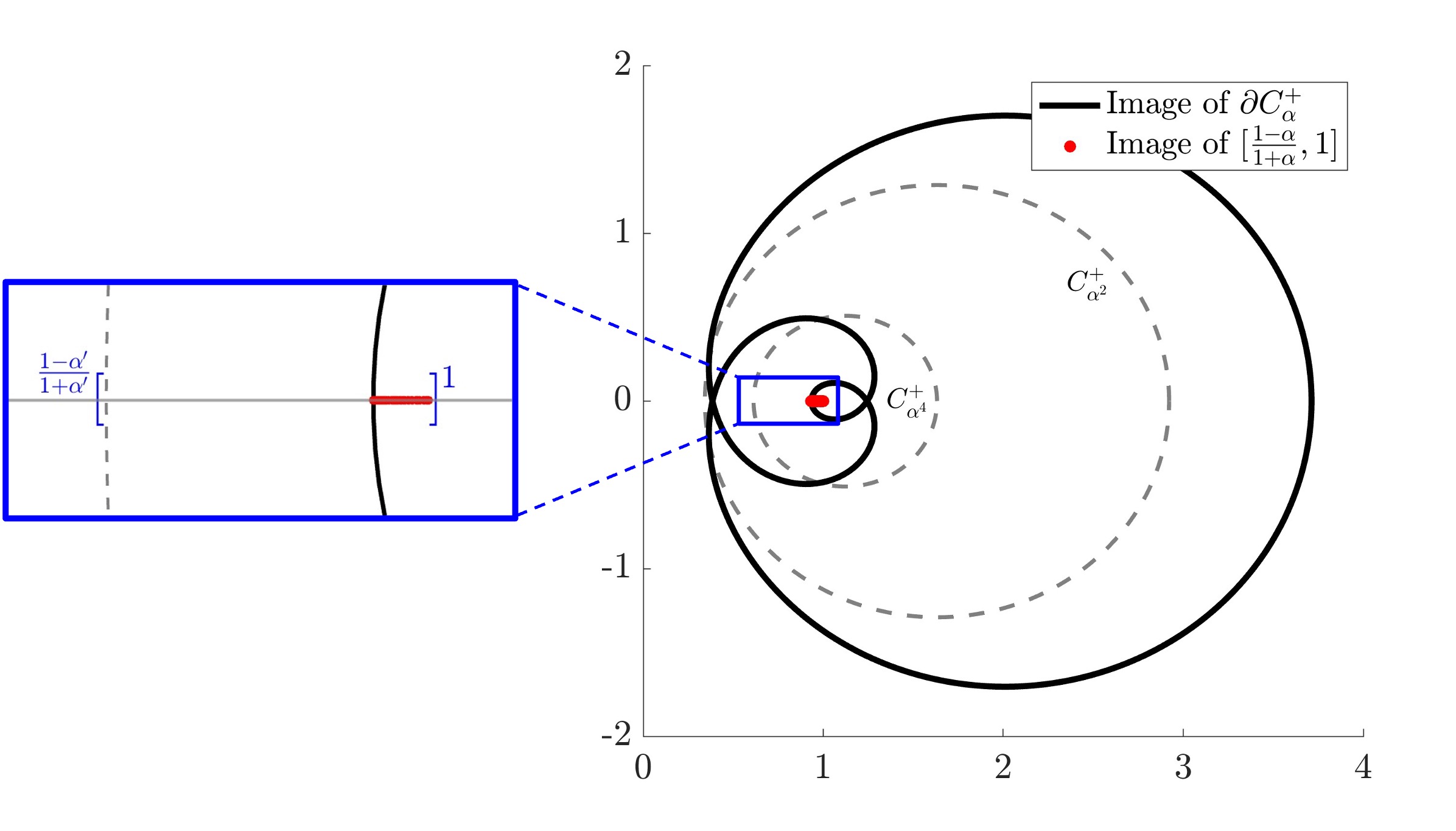}
    \caption{Result of applying $f_{\alpha}$ to portions of $C_{\alpha}^+$ (see \cref{defn: circles_of_Apollonius}), where $f_{\alpha}$ is the rational function defined according to \eqref{eqn: f_alpha} and \eqref{eqn: weighted_coefficients} with $l = \frac{1-\alpha}{1+\alpha}$ and $\alpha = 0.7$. A subplot focused on the real axis verifies \cref{prop: dwh_convergence_bound} for the corresponding value of $\alpha'$ (in this case $\frac{1-\alpha'}{1+\alpha'} \approx 0.6$). For reference we also mark the regions $C_{\alpha^2}^+$ and $C_{\alpha^4}^+$.}
    \label{fig: dwh_circle}
\end{figure}

\indent Unfortunately, the preceding results do not extend to \textbf{IF-DWH}, as the rational function applied at a given step of the weighted iteration is \textit{not} guaranteed to map one circle of Apollonius to another in general. See for example \cref{fig: dwh_circle}, which demonstrates how different portions of $C_{\alpha}^+$ are mapped by the rational function
\begin{equation}\label{eqn: f_alpha}
    f_{\alpha}(z) = z \frac{az^2+b}{(a+b-1)z^2+1}
\end{equation}
for $a$ and $b$ the weighted coefficients defined according to \eqref{eqn: weighted_coefficients} with $l = \frac{1-\alpha}{1+\alpha}$. The choice of $l$ here reflects the fact that
\begin{equation}\label{eqn: real_part_of_C}
    C_{\alpha}^+ \cap {\mathbb R} = \left[ \frac{1-\alpha}{1+\alpha}, \frac{1+\alpha}{1-\alpha} \right].
\end{equation}
\indent While $f_{\alpha}$ does not map $\partial C_{\alpha}^+$ to another circle of Apollonius, it does map $[\frac{1-\alpha}{1+\alpha}, 1]$ to another real interval bounded above by one, as guaranteed by our choice of weighted coefficients. Hence, we obtain a convergence result for the dynamic iteration -- one that is again described by the circles of Apollonius -- by showing that this interval is contained in $[\frac{1-\alpha'}{1+\alpha'}, 1]$ for some $\alpha' \in (0,\alpha)$. This is equivalent to showing $f_{\alpha}([\frac{1-\alpha}{1+\alpha},1]) \subseteq C_{\alpha'}^+$. \\
\indent Before stating such a result, we first prove a technical lemma.

\begin{lem}\label{lem: technical_alpha_lemm}
    For any $\alpha \in (0,1)$ we have
    $$ 0 < \frac{1+\alpha}{\sqrt{1-\alpha}} \cdot \frac{1 - \alpha^{\frac{1}{2} + \frac{1}{2} \sqrt{\frac{1}{1-\alpha}}}}{1+\alpha^{\frac{1}{2}+\frac{1}{2}\sqrt{\frac{1}{1-\alpha}}}} < 1.$$
\end{lem}
\begin{proof}
    Rewrite the upper bound as
    \begin{equation}\label{eqn: rewrite_upper_bound}
        \aligned 
        (&1+\alpha)\left(1-\alpha^{\frac{1}{2}+\frac{1}{2}\sqrt{\frac{1}{1-\alpha}}} \right) < \sqrt{1-\alpha} \left(1+\alpha^{\frac{1}{2}+\frac{1}{2}\sqrt{\frac{1}{1-\alpha}}} \right) \\
        \iff \; \; & 1 + \alpha - \alpha^{\frac{1}{2}+\frac{1}{2}\sqrt{\frac{1}{1-\alpha}}} - \alpha^{\frac{3}{2}+\frac{1}{2}\sqrt{\frac{1}{1-\alpha}}} < \sqrt{1-\alpha}+\sqrt{1-\alpha}\cdot \alpha^{\frac{1}{2}+\frac{1}{2}\sqrt{\frac{1}{1-\alpha}}} \\
        \iff \; \;  & 1 + \alpha - \sqrt{1-\alpha} < \alpha^{\frac{1}{2}+\frac{1}{2}\sqrt{\frac{1}{1-\alpha}}} \left( 1+\alpha + \sqrt{1-\alpha} \right)  \\
        \iff \; \; &\frac{1+\alpha-\sqrt{1-\alpha}}{1+\alpha+\sqrt{1-\alpha}} < \alpha^{\frac{1}{2}+\frac{1}{2}\sqrt{\frac{1}{1-\alpha}}}.
        \endaligned
    \end{equation}
    Letting $\alpha = 1-t$ for $0 < t < 1$, \eqref{eqn: rewrite_upper_bound} becomes 
    \begin{equation}\label{eqn: change_variables}
        (1-t)^{\frac{1}{2}+\frac{1}{2}\sqrt{\frac{1}{t}}} > \frac{2-(t+\sqrt{t})}{2+(\sqrt{t}-t)}.
    \end{equation}
    This now follows from Bernoulli's inequality, noting that $\sqrt{t}-t > 0$ since $t \in (0,1)$:
    \begin{equation}\label{eqn: apply_bernoulli}
        (1-t)^{\frac{1}{2}+\frac{1}{2}\sqrt{\frac{1}{t}}} \geq 1 - t\left( \frac{1}{2}+\frac{1}{2}\sqrt{\frac{1}{t}}\right) = \frac{2-(t+\sqrt{t})}{2} > \frac{2-(t+\sqrt{t})}{2+(\sqrt{t}-t)}.
    \end{equation}
    The corresponding lower bound is a consequence of the fact that $\alpha$ and $\alpha^{\frac{1}{2}+\frac{1}{2}\sqrt{\frac{1}{1-\alpha}}}$ are positive but strictly smaller than one. 
\end{proof}

We can now present our main convergence result for the dynamically weighted iteration.

\begin{prop} \label{prop: dwh_convergence_bound}
    Let $\alpha \in (0,1)$. If
    $f_{\alpha}(z) = z \frac{az^2 + b}{(a+b-1)z^2 + 1}$ for $a,b > 0$ defined according to \eqref{eqn: weighted_coefficients} with $l = \frac{1-\alpha}{1+\alpha}$, then 
    $$ f_{\alpha} \left(\left[ \frac{1-\alpha}{1+\alpha}, \; 1 \right] \right) \subseteq \left[ \frac{1 - \alpha'}{1 + \alpha'}, \; 1\right]$$
    for $\alpha' = \alpha^{2 + (1-\alpha)^{-1/2}}. $
\end{prop}
\begin{proof}
    It is sufficient to show $\left| \frac{1 - f_{\alpha}(z)}{1 + f_{\alpha}(z)} \right| \leq \alpha'$ for $z \in [ \frac{1-\alpha}{1+\alpha}, 1]$. Since the dynamic weights \eqref{eqn: weighted_coefficients} guarantee that $\left| \frac{1-f_{\alpha}(z)}{1 + f_{\alpha}(z)} \right|$ is maximized at $\frac{1-\alpha}{1+\alpha}$, we need only verify the inequality at the left endpoint of the interval. With this in mind, and recalling $a = \frac{1}{4}(b-1)^2$, we start by observing
    \begin{equation}\label{eqn: expand_and_factor}    \aligned
        \frac{1 - f_{\alpha}(z)}{1 + f_{\alpha}(z)} &= \frac{1 - \frac{az^3 + bz}{(a+b-1)z^2 + 1}}{1 + \frac{az^3 + bz}{(a+b-1)z^2 + 1}} \\
        & = \frac{(a+b-1)z^2 + 1 - az^3 - bz}{(a+b-1)^2z^2 + 1 + az^3 + bz} \\
        & = \frac{(1-z)(1 - (b-1)z + az^2)}{(1+z)(1 + (b-1)z + az^2)} \\
        & = \frac{(1-z)(1 - \frac{1}{2}(b-1)z)^2}{(1+z)(1 + \frac{1}{2} (b-1)z)^2}.
        \endaligned
    \end{equation}
    Set $\Theta(z) = \frac{1-z}{1+z}$ and $ \Psi(z) = \frac{(1-\frac{1}{2}(b-1)z)^2}{(1+\frac{1}{2}(b-1)z)^2}$ so that \eqref{eqn: expand_and_factor} becomes
    \begin{equation}\label{eqn: relabel}
        \frac{1 - f_{\alpha}(z)}{1 + f_{\alpha}(z)}  =  \Theta(z)\Psi(z).
    \end{equation}
    Since $\Theta \left( \frac{1-\alpha}{1+\alpha} \right) = \alpha$, we obtain the desired result provided $\Psi \left( \frac{1 - \alpha}{1+\alpha} \right) < \alpha^{1 + \sqrt{\frac{1}{1-\alpha}}}$. Noting $b -1 > 0$, it is easy to see
    \begin{equation}\label{eqn: first_condition}
        \aligned 
        &|\Psi(z)| > \alpha^{1 + \sqrt{\frac{1}{1-\alpha}}} \iff \left| \frac{ 1 - \frac{1}{2}(b-1)z}{1 + \frac{1}{2}(b-1)z} \right| > \alpha^{\frac{1}{2} + \frac{1}{2} \sqrt{\frac{1}{1-\alpha}}} \\
        & \iff z > \frac{2}{b-1} \cdot \frac{1 + \alpha^{\frac{1}{2} + \frac{1}{2} \sqrt{\frac{1}{1-\alpha}}}}{1 - \alpha^{\frac{1}{2} + \frac{1}{2} \sqrt{\frac{1}{1-\alpha}}}}  \; \; \; \text{or} \; \; \; z < \frac{2}{b-1} \cdot \frac{1 - \alpha^{\frac{1}{2} + \frac{1}{2}\sqrt{\frac{1}{1-\alpha}}}}{1 + \alpha^{\frac{1}{2} + \frac{1}{2}\sqrt{\frac{1}{1-\alpha}}}}.
        \endaligned 
    \end{equation}
    Hence, we complete the proof by showing
    \begin{equation}\label{eqn: inequality_one}
        \frac{2}{b-1} \cdot \frac{1-\alpha^{\frac{1}{2}+\frac{1}{1}\sqrt{\frac{1}{1-\alpha}}}}{1+\alpha^{\frac{1}{2} + \frac{1}{2}\sqrt{\frac{1}{1-\alpha}}}} < \frac{1-\alpha}{1+\alpha} < \frac{2}{b-1} \cdot \frac{1 + \alpha^{\frac{1}{2} + \frac{1}{2} \sqrt{\frac{1}{1-\alpha}}}}{1 - \alpha^{\frac{1}{2} + \frac{1}{2} \sqrt{\frac{1}{1-\alpha}}}}
    \end{equation}
    or equivalently
    \begin{equation}\label{eqn: inequality_two}
          \frac{1+\alpha}{1-\alpha} \cdot \frac{1-\alpha^{\frac{1}{2}+\frac{1}{2}\sqrt{\frac{1}{1-\alpha}}}}{1+\alpha^{\frac{1}{2} + \frac{1}{2}\sqrt{\frac{1}{1-\alpha}}}} <  \frac{b-1}{2} < \frac{1+\alpha}{1-\alpha} \cdot \frac{1 + \alpha^{\frac{1}{2} + \frac{1}{2} \sqrt{\frac{1}{1-\alpha}}}}{1 - \alpha^{\frac{1}{2} + \frac{1}{2} \sqrt{\frac{1}{1-\alpha}}}}.
    \end{equation}
   This is a consequence of the following key inequality:
   \begin{equation}\label{eqn: key_inequality}
       \frac{1}{\sqrt{1-\alpha}} < \frac{b-1}{2} < \frac{1+\alpha}{1-\alpha}.
   \end{equation}
   The right side of \eqref{eqn: key_inequality} follows from \cite[Remark 1]{optimizing_halley}\footnote{Note that the weights $a$ and $b$ in \eqref{eqn: weighted_coefficients} are labeled $b$ and $a$ in \cite{optimizing_halley}.}  and implies the upper bound of \eqref{eqn: inequality_two} via
   \begin{equation}\label{eqn: verify_upper}
        \frac{b-1}{2} < \frac{1+\alpha}{1-\alpha} < \frac{1+\alpha}{1-\alpha} \cdot \frac{1+\alpha^{\frac{1}{2} + \frac{1}{2} \sqrt{\frac{1}{1-\alpha}}}}{1 - \alpha^{\frac{1}{2}+ \frac{1}{2}\sqrt{\frac{1}{1-\alpha}}}}.
   \end{equation}
   The corresponding lower bound, meanwhile, can be obtained as
   \begin{equation}\label{eqn: verify_lower}
       \frac{b-1}{2} > \frac{1}{\sqrt{1-\alpha}} > \frac{1+\alpha}{1-\alpha} \cdot \frac{1-\alpha^{\frac{1}{2}+\frac{1}{2}\sqrt{\frac{1}{1-\alpha}}}}{1+\alpha^{\frac{1}{2} + \frac{1}{2}\sqrt{\frac{1}{1-\alpha}}}}
   \end{equation}
   via \cref{lem: technical_alpha_lemm}. The left side of \eqref{eqn: key_inequality} can be verified in Mathematica for all $\alpha \in (0,1)$. 
\end{proof}


As $\alpha \rightarrow 0$, \cref{prop: dwh_convergence_bound} implies asymptotically cubic convergence for the weighted Halley iteration, reflecting its eventual reduction to the standard version. Moreover, it guarantees fast convergence for \textbf{IF-DWH} on any problem with real eigenvalues, as captured by the following corollary. 

\begin{cor}\label{cor: log_log_convergence}
    Suppose $(A,B)$ has eigenvalues in $[-1,-l) \cup (l,1]$ for $l \ll 1$. Then \text{\normalfont \bf IF-DWH} requires at most $O\left(\log \left( \log \left( \frac{1+l}{2l} \right) \right) \right)$ iterations to drive the eigenvalues of $(A,B)$ into $\left[-1,-\frac{1}{2} \right) \cup \left(\frac{1}{2}, 1\right]$.
\end{cor}
\begin{proof}
    Write $l = \frac{1-\alpha}{1+\alpha}$ for some $\alpha \in (0,1)$. \cref{prop: dwh_convergence_bound} guarantees that after one iteration of \textbf{IF-DWH} the eigenvalues of $(A,B)$ lie in $[-1,-l') \cup (l',1]$ for $l' \geq \frac{1-\alpha'}{1+\alpha'}$ with
    \begin{equation}\label{eqn: bound_alpha_prime}
        \aligned 
        \alpha' &= \alpha^{2+\sqrt{\frac{1}{1-\alpha}}} \\
        &= (1-(1-\alpha))^2 (1-(1-\alpha))^{\sqrt{\frac{1}{1-\alpha}}} \\
        & \leq \frac{1-2(1-\alpha)+(1-\alpha)^2}{1+(1-\alpha)\sqrt{\frac{1}{1-\alpha}}} \\
        & = 1  - \sqrt{1-\alpha} \left(\frac{1 + 2\sqrt{1-\alpha} - (1-\alpha)^{3/2}}{1+\sqrt{1-\alpha}} \right) \\
        & \leq 1 - \sqrt{1-\alpha}.
        \endaligned
    \end{equation}
    Equivalently, $1-\alpha' \geq \sqrt{1 - \alpha}$. Applying this recursively,\footnote{After the first iteration, \textbf{IF-DWH} does not necessarily apply the rational function $f_{\alpha}$ from \eqref{eqn: f_alpha}. Nevertheless, an iterative method that \textit{does} apply $f_{\alpha}$ for the sequence of $\alpha$'s defined by \cref{prop: dwh_convergence_bound} shows strictly slower convergence than \textbf{IF-DWH}. This is equivalent to noting that the rational function selected by \textbf{IF-DWH} does ``better" than $f_{\alpha}$ beyond the first iteration, where ``better" means it maps the eigenvalues closer to one.} we conclude that after $j$ iterations eigenvalues lie in $[-1,-l_j) \cup (l_j,1]$ for $l_j \geq \frac{1-\alpha_j}{1+\alpha_j}$ and $\alpha_j \leq 1 - (1-\alpha)^{2^{-j}}$. Hence, to guarantee that the spectrum is contained in $[-1,-\frac{1}{2}) \cup (\frac{1}{2},1]$ after step $j$ it is sufficient to ensure $\alpha_j < \frac{1}{2}$, which can be accomplished by choosing $j$ so that
    \begin{equation}\label{eqn: bound_alpha_j}
        1 - (1-\alpha)^{2^{-j}} < \frac{1}{2}  \; \; \iff \; \; j > \frac{\log(\log((1-\alpha)^{-1})) - \log(\log(2))}{\log(2)}.
    \end{equation}
    We complete the proof by noting $\frac{1}{1-\alpha} = \frac{1+l}{2l}$.
\end{proof}

We consider only small values of $l$ here since the number of iterations required by \textbf{IF-DWH} increases as $l \rightarrow 0 $ (and of course $l$ could initially be greater than $\frac{1}{2}$ or large enough that a constant number of iterations is sufficient). Importantly, \cref{cor: log_log_convergence} implies that \textbf{IF-DWH} requires at most
\begin{equation}\label{eqn: if_dwh_steps}
    \left\lceil \frac{\log(\log(\frac{1+l}{2l})) - \log(\log(2))}{\log(2)} \right\rceil + \left\lceil \frac{\log(\log(\frac{1}{\delta})) - \log(\log(2))}{\log(3)} \right\rceil 
\end{equation}
iterations to drive eigenvalues into $[-1,-1+\delta) \cup (1-\delta,1]$ for \textit{any} $\delta,l \ll 1$. This can be shown in two steps:
\begin{enumerate}
    \item \cref{cor: log_log_convergence} guarantees that the first  $ \left\lceil \frac{\log(\log(\frac{1+l}{2l})) - \log(\log(2))}{\log(2)} \right\rceil$ iterations will map the initial spectral bound $[-1,-l) \cup (l,1]$ to $[-1,-\frac{1}{2}) \cup (\frac{1}{2},1]$.
    \item From there, the (at least) cubic convergence of the weighted Halley iteration ensures that at most $\left\lceil \frac{\log(\log(\frac{1}{\delta})) - \log(\log(2))}{\log(3)} \right\rceil$ additional iterations are needed to drive eigenvalues within $\delta$ of $\pm 1$.
\end{enumerate}
It is necessary to make this argument in two steps since the recurrence $1-\alpha' \geq \sqrt{1-\alpha}$ used in the proof of \cref{cor: log_log_convergence} slows dramatically as $\alpha \rightarrow 0$. \\
\indent Of course, \cref{cor: log_log_convergence} describes the convergence of eigenvalues to $\pm 1$, and in our general setting this doesn't fully characterize the challenge of sign function approximation. To provide a complete projector bound for \textbf{IF-DWH}, again for definite pencils, we can combine \cref{cor: log_log_convergence} with the following lemma. 

\begin{lem}\label{lem: shifted_dwh_bound}
    Suppose $(A,B)$ is a definite pencil with eigenvalues in $[-1,-l) \cup (l,1]$ and let $[A_j,B_j,l_j] = \text{\normalfont \bf IF-DWH}(A, B, j, l)$. Then
    $$ \|B_j^{-1}A_j - \text{\normalfont sign}(B^{-1}A)\|_2 \leq \frac{\|(A,B)\|_2}{\gamma(A,B)} (1-l_j).$$
\end{lem}
\begin{proof}
    Let $X$ be the invertible eigenvector matrix of $(A,B)$ satisfying 
    \begin{equation}
    (X^HAX,X^HBX) = (\Lambda_A,\Lambda_B) = (\text{diag}(\alpha_1, \ldots, \alpha_n), \text{diag}(\beta_1,\ldots,\beta_n))
    \label{eqn: standard_choice_X}
    \end{equation}
    for $\alpha_i,\beta_i \in {\mathbb R}$ with $\alpha_i^2 + \beta_i^2 = 1$, which exists since $(A,B)$ is definite (see \cite[Chapter VI]{stewart1990matrix}). Since $(A_j,B_j)$ has the same right eigenvectors as $(A,B)$, $X$ diagonalizes $B_j^{-1}A_j$. Writing $B_j^{-1}A_j = X\Lambda_j X^{-1}$ for $\Lambda_j$ diagonal and noting $B^{-1}A = X \Lambda_B^{-1} X^H X^{-H} \Lambda_A X^{-1} = X\Lambda_B^{-1}\Lambda_AX^{-1}$, we have
    \begin{equation}
        \aligned
        \|B_j^{-1}A_j - \text{sign}(B^{-1}A)\|_2 &= \|X\Lambda_jX^{-1} - X\text{sign}(\Lambda_B^{-1}\Lambda_A)X^{-1}\|_2 \\
        & \leq \kappa_2(X) \| \Lambda_j - \text{sign}(\Lambda_B^{-1} \Lambda_A) \|_2 \\
        & \leq \kappa_2(X)(1-l_j),
        \endaligned
        \label{eqn: shifted_eqn_one}
    \end{equation}
    where the last inequality follows from $\Lambda(A_p,B_p) \subseteq (-1,-l_j) \cup (l_j,1)$. We complete the proof by noting
    \begin{equation}\label{eqn: Crawford_number_bound}
        \kappa_2(X) \leq \frac{\|(A,B)\|_2}{\gamma(A,B)}.
    \end{equation}
    For the details see \cite[Proof of Theorem 2.3]{ELSNER1982341}.
\end{proof}

The results derived in this section imply, as we should expect, that convergence of sign-function-based methods is dependent on both the locations of the eigenvalues of $(A,B)$ -- relative to $\pm 1$ -- and the conditioning of its eigenvectors. The latter is present implicitly in \cref{prop: general_error_bound} and \cref{lem: shifted_dwh_bound} via the pseudospectrum and the Crawford number, respectively. \\
\indent Moreover, we can now describe the set of problems for which the methods considered in this section are fast, which -- as mentioned in \cref{section: intro} -- means they require at most $O(\log(\frac{n}{\delta})T_{\text{MM}}(n))$ operations to compute a projector to forward error $\delta$ in the spectral norm. Each of our algorithms is iterative, with one iteration requiring $O(T_{\text{MM}}(n))$ 
operations, so they are fast as long as forward error falls below $\delta$ after $O(\log( \frac{n}{\delta}))$ steps. \cref{prop: general_error_bound} implies that this is the case for sign-function-based methods like \textbf{IF-Newton} and \textbf{IF-Halley} (in exact arithmetic) when $1-\alpha$ and $\epsilon$ are at least polynomial in $\delta$ and $n^{-1}$. \\
\indent \textbf{IF-DWH} is even faster on definite problems. Together, \cref{cor: log_log_convergence} and \cref{lem: shifted_dwh_bound} imply that \textbf{IF-DWH} will produce a suitably accurate projector after at most $O(\log(\log(\frac{n}{\delta})))$ iterations, provided both the initial eigenvalue bound $l$ and the quantity $\frac{\gamma(A,B)}{\|(A,B)\|_2}$ are at least polynomial in $\delta$ and $n^{-1}$. When this occurs, recalling \eqref{eqn: Crawford_number_bound} and applying standard Bauer-Fike type results (e.g., \cite[Theorem 2.20]{arXiv}), we obtain the same kind of pseudospectral bound for $\Lambda_{\epsilon}(A,B)$ that guarantees fast convergence for \textbf{IF-Newton} and \textbf{IF-Halley}. \\
\indent Since $\Lambda_{\epsilon}(A,B) \cap \left\{ z : |z| = 1 \right\} = \emptyset$ implies $d_{(A,B)} \geq 2\epsilon$, \cite[Theorem 1]{Bai:CSD-94-793} tells a similar story for \textbf{IRS}. In total, these results suggests a general guideline for computing arbitrary spectral projectors $P_R$ and $P_L$ via the computational framework presented in this paper: the methods we develop provide fast approximations to $P_R$ and $P_L$ as long as $\Lambda_{\epsilon}(A,B)$ is well-separated from the boundary of $S$ for $\epsilon$ not too small (at least polynomial in $n^{-1}$ and the desired accuracy).

\begin{remark}
    We might ask whether it is possible to replicate the optimization of \textbf{IF-DWH} in general. That is, can we similarly modify the Newton/Halley iterations to achieve fast convergence on all of $C_{\alpha}^{\pm}$, not just their intersection with the real line? Unfortunately, the answer is essentially no; see \cref{appendix: Zolo}, which discusses the optimal rational approximation to the sign function on these disks. 
\end{remark}

\section{Numerical Examples}\label{section: examples}
In this section, we present a handful of numerical examples to explore the empirical performance of \textbf{IRS}, \textbf{IF-Newton}, \textbf{IF-Halley}, and \textbf{IF-DWH}. All results were obtained in Matlab R2025a. \\
\indent Our first test is a $500 \times 500$ pencil
\begin{equation}
    (A,B) = (X^H\Lambda X, X^HX),
    \label{eqn: example_1}
\end{equation}
where $X$ is invertible and $\Lambda$ is diagonal. We set 
\begin{equation}
    \Lambda = \begin{pmatrix} \Lambda_+ & 0 \\ 0 & \Lambda_- \end{pmatrix}
    \label{eqn: example_1_lambda}
\end{equation}
for $\Lambda_+, \Lambda_- \in {\mathbb R}^{250 \times 250}$ with diagonal entries sampled from ${\mathbb R}_{>0}$ and ${\mathbb R}_{<0}$, respectively.  By construction, the eigenvalues of $(A,B)$ belong to the diagonal of $\Lambda$, and the columns of $X^{-1}$ are right eigenvectors. We choose real eigenvalues here to test the potential improved convergence of \textbf{IF-DWH}; in this case, the pencil $(A,B)$ is definite. For a repeat of this experiment on an indefinite pencil, see \cref{appendix: examples}. \\
\indent The benefit of this construction is that it allows us to vary both the locations of the eigenvalues of $(A,B)$ as well as the conditioning of its eigenvectors. We consider in particular the following:
\begin{itemize}
    \item Well-separated eigenvalues: diagonal entries of $\Lambda$ are sampled uniformly from $(-4,-1) \cup (1,4)$. Well-separated here refers to distance from the imaginary axis and \textit{not} spacing between individual eigenvalues. 
    \item Poorly-separated eigenvalues: diagonal entries of $\Lambda_+$ and $\Lambda_{-}$ are $\pm 4/k$ for $1 \leq k \leq 250$.
    \item Well-conditioned eigenvectors: $X$ is a standard, complex Gaussian matrix, which is well-conditioned in the sense of \cite{RUDELSON2008600}. 
    \item Poorly-conditioned eigenvectors: $X$ is drawn randomly (again standard, complex Gaussian) and modified so that $\kappa_2(X) = 10^5$. This is done by computing a full SVD of $X$ and subtracting off a rank one matrix constructed from the singular vectors corresponding to its smallest singular value -- i.e., to ensure $\sigma_{\min}(X) = 10^{-5}\|X\|$ while leaving all other singular values unchanged.
\end{itemize} 

\indent Suppose we are interested in computing a unitary projector $P$ onto the right deflating subspace of $(A,B)$ corresponding to the right half plane (equivalently, corresponding to the eigenvalues in $\Lambda_+$). A stand-in for this projector can be obtained (without using matrix inversion) by computing an RQ factorization $X = RQ$ and taking $P = WW^H$ for $W$ containing the first 250 columns of $Q^H$. This follows from the fact that the right eigenvectors of $(A,B)$ are columns of $X^{-1} = Q^HR^{-1}$. Hence, we can compute the ``exact" projector $P$ without incurring significant error for either choice of $X$ and any set of eigenvalues.   \\
\indent In combination with the \textbf{GRURV} algorithm of Ballard et al.\ \cite{grurv}, each of the iterative methods considered in this paper can also produce an approximation of $P$. At a given iteration, each method yields a pencil $(A_j,B_j)$ with eigenvalues close to zero or one. An approximate projector can then be obtained as $\widetilde{P} = UU^H$ for $U$ a matrix containing the first 250 columns of the U-factor produced by \textbf{GRURV} when applied to $\frac{1}{2}B_j^{-1}(A_j+B_j)$ or -- in the case of \textbf{IRS} --  $(A_j+B_j)^{-1}A_j$. Note that in this approach, \textbf{IRS} must apply an initial Möbius transformation mapping the imaginary axis to the unit circle. 

\begin{figure}[t]
    \centering
    \includegraphics[width=\linewidth]{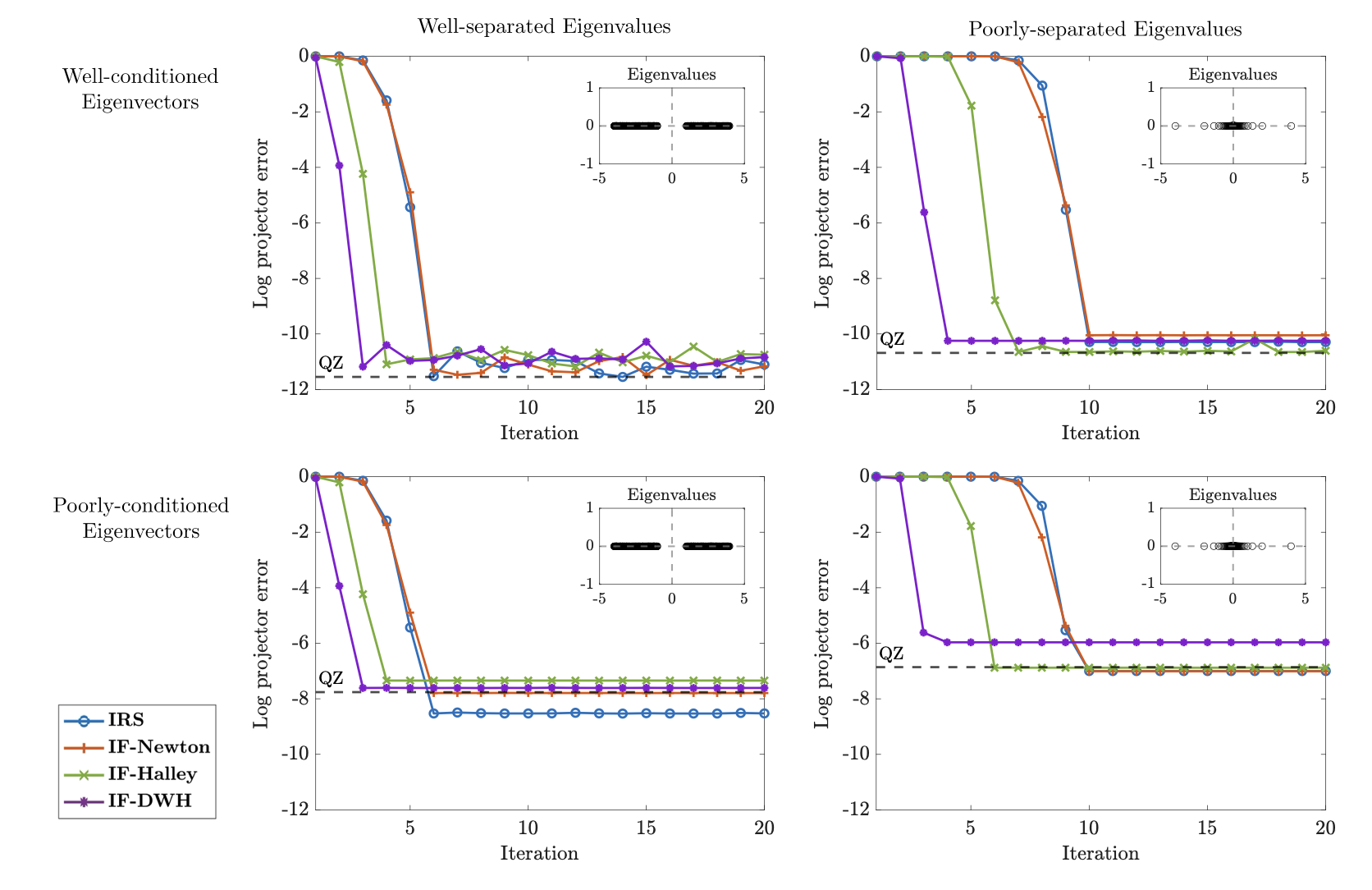}     
    \caption{Forward error for \textbf{IRS}, \textbf{IF-Newton}, \textbf{IF-Halley}, and \textbf{IF-DWH} when used to compute a spectral projector of a $500 \times 500$ pencil $(A,B)$ constructed according to \eqref{eqn: example_1} and \eqref{eqn: example_1_lambda}. Each plot corresponds to a different combination of eigenvalues and eigenvectors. For context, evaluating \eqref{eqn: if_dwh_steps} with $\delta = \varepsilon/\kappa_2(X)$ for $\varepsilon$ the error attained by QZ yields $5$ when eigenvalues are well separated from zero and $8$ when they are poorly separated. This provides an upper bound on the number of iterations required by \textbf{IF-DWH} to match QZ (in exact arithmetic).}
    \label{fig: proj_4x4}
\end{figure}

\indent \cref{fig: proj_4x4} presents the associated  forward projector error $\log_{10}(\|P-\widetilde{P}\|_2)$ for each method and all four combinations of eigenvalue placement and eigenvector conditioning. For comparison, we mark on each plot the accuracy of the standard QZ algorithm \cite{QZ}, which can approximate $P$ by computing first a full eigendecomposition (by way of Schur form) and then a QR factorization of corresponding eigenvectors. This provides a benchmark for the error produced by a backward-stable method; the gap between this error and machine precision can be explained by classical error bounds (see e.g., \cite{Stewart_subspace, Kagstrom_poromaa}). \\
\indent Two trends jump out from these plots:\ (1) eigenvalues near the imaginary axis drive up the number of iterations required by all of the methods to converge and (2) eigenvector conditioning impacts attainable forward error but not convergence. With the exception of \textbf{IF-DWH}, each of the iterative methods is capable of matching the error produced by QZ, which recall requires an expensive full Schur decomposition.  

\begin{figure*}[t!]
    \centering 
    \subfigure[\centering $l_0 = 2.53 \times 10^{-4}$ / 14 iterations.]{\includegraphics[width=.45\linewidth]{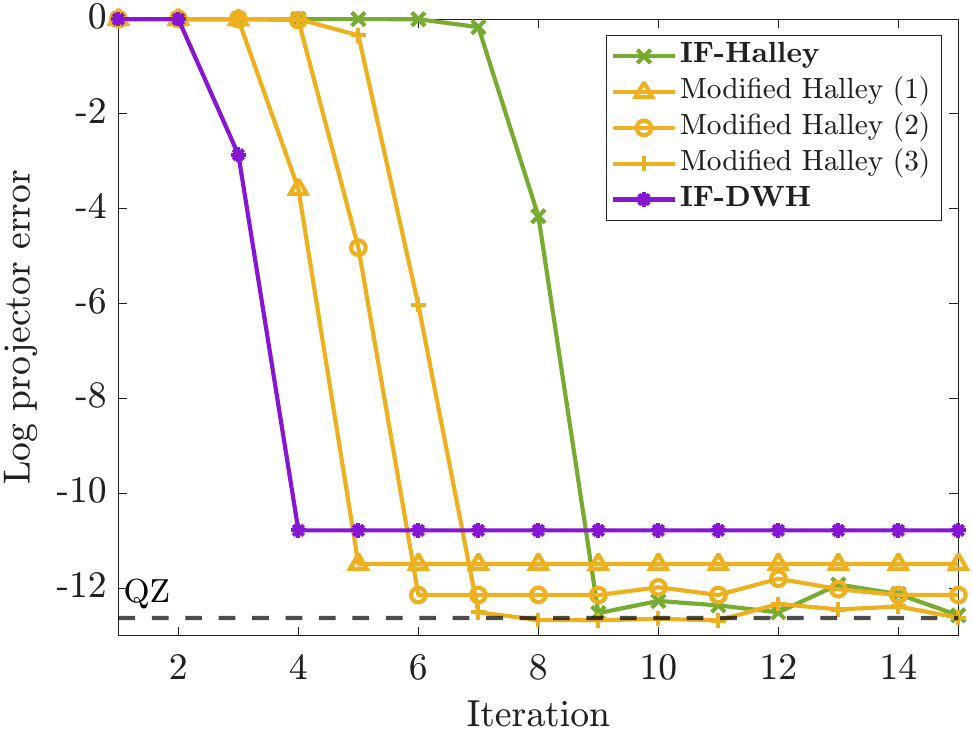}} \hspace{5mm}
    \subfigure[\centering $l_0 = 1.6 \times 10^{-5}$/ 18 iterations.]{\includegraphics[width=.45\linewidth]{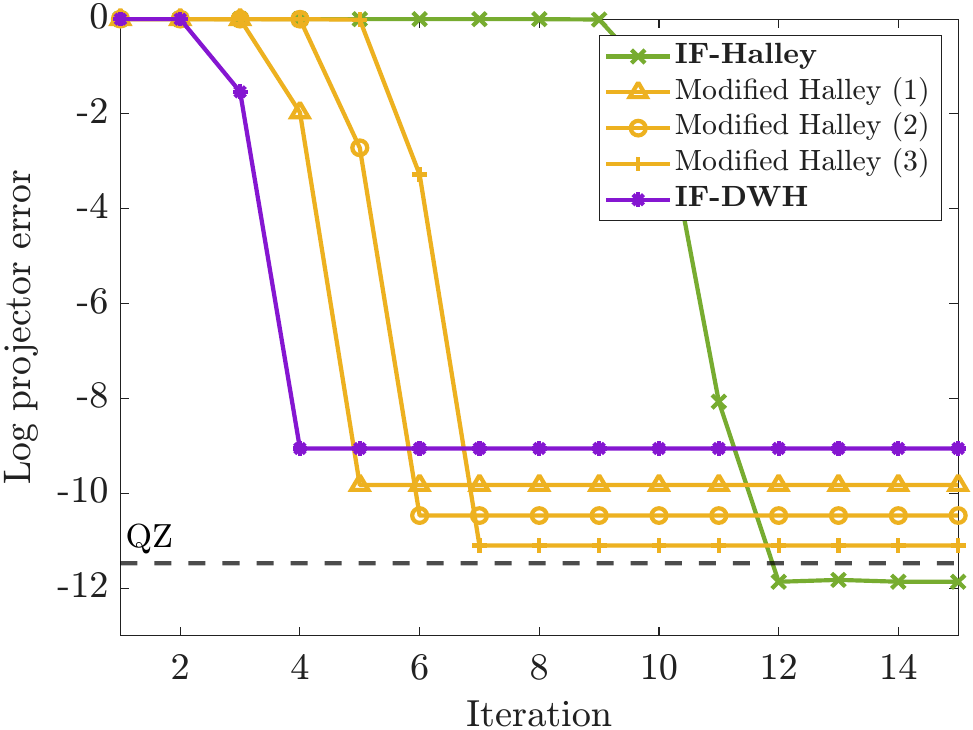}}
    \caption{Forward projector error for methods based on the Halley iteration. The input pencil $(A,B)$ is again $500 \times 500$ and constructed according to \eqref{eqn: example_1}, this time with a Haar unitary eigenvector matrix $X$ and eigenvalues $\pm 4/k^{3/2}$ (plot (a)) and $\pm 4/k^2$ (plot (b)) for $1 \leq k \leq 250$, which corresponds to the listed values of $l_0$. For each modified iteration, we note in parentheses the number of standard Halley steps applied before switching to the weighted version. For context, we also label the plots with the number of iterations required by \textbf{IRS} and \textbf{IF-Newton}; to be more efficient overall, a Halley-based method must converge in fewer than half as many.}\label{fig: modified_projector_comp}
\end{figure*}

\indent \textbf{IF-DWH} converges not only faster than \textbf{IF-Halley} but fast enough to be more efficient than \textbf{IRS} or \textbf{IF-Newton}, though it is also the least accurate. While all of the methods stagnate on the harder problems -- a sign that the variance in iterations beyond a certain point is masked by the larger forward error -- \textbf{IF-DWH} does so without reaching the QZ benchmark, particularly when eigenvalues are close to the imaginary axis. \\
\indent This is likely due to inherent instability in the weighted iteration. When $l_0$, or more generally $l_j$, is small (as is the case in the first few steps of the poorly-separated plots of \cref{fig: proj_4x4}), the weights in \eqref{eqn: weighted_scalar_Halley} are unbalanced in the sense that $a_j,b_j,$ and $c_j$ are large while $d_j = 1$. As a result, the input to the second QR factorization computed by \textbf{IF-DWH} (line 11 in \cref{alg:IF_DWH}) is potentially ill-conditioned, specifically if $Q_{12}^HA_j$ is itself ill-conditioned. When this occurs, significant error is incurred after only one iteration, enough to move \textbf{IF-DWH} away from the error achieved by the other methods. Interestingly, this phenomenon was not reported in the original version of the dynamically weighted Halley iteration \cite{optimizing_halley} and is possibly specific to our general, inverse-free implementation. \\
\indent To overcome this drawback, we suggest using a modified Halley iteration:\ if $l_0$ is initially small (recall that $l_0$ is an input for \textbf{IF-DWH}) run a few standard Halley iterations -- i.e., \textbf{IF-Halley} -- before switching to \textbf{IF-DWH}. The result is a method that converges somewhere in between \textbf{IF-DWH} and \textbf{IF-Halley}; when the gap between the two is significant, this modified iteration can both avoid the error stagnation of \textbf{IF-DWH} while converging fast enough to beat \textbf{IRS} and \textbf{IF-Newton}. We verify this in \cref{fig: modified_projector_comp}.

\indent More testing is necessary to provide rigorous guidelines for choosing the number of standard Halley iterations to run before defaulting to \textbf{IF-DWH}. It may also be worth initially running one of the second order methods instead of \textbf{IF-Halley} to further save on computational costs.

\section{Conclusion}\label{section: conclusion}
\indent This paper presented a high-level, inverse-free framework for computing the spectral projectors of any regular matrix pencil $(A,B)$. We used this framework to state a number of new iterative methods for the problem, providing theoretical and empirical evidence for their efficacy. In doing so, we have also reinforced the practical value of the high-level framework and its accompanying Indicator Approximation Problem. For practitioners looking to handle this fundamental problem in numerical linear algebra without relying on inversion, our work offers a wealth of new options.

\section{Acknowledgements}\label{section: acknowledgments} 
This project was supported by Graduate Fellowships for STEM Diversity (GFSD) and the following National Science Foundation (NSF) grants:\ DMS-2154099, FRG award DMS-1952786, and MSPRF 2402027. Special thanks to Yuji Nakatsukasa for suggesting the indicator-approximation perspective and to Heather Wilber for pointing out the references used  in \cref{appendix: Zolo}.  Thanks also to Volker Mehrmann and two anonymous referees for helpful correspondence and feedback regarding an earlier draft of this paper.  

\bibliographystyle{abbrv}{}
\bibliography{bib}

\renewcommand{\arraystretch}{1}
\appendix
\section{Fast Generalized Schur Decomposition}\label{section: fast_schur}
In this appendix, we verify that the pseudospectral divide-and-conquer approach of \cite{arXiv} can be adapted to produce the generalized Schur decomposition of any pencil $(A,B)$. The result is an inverse-free and randomized  algorithm that can produce an approximate Schur decomposition, of arbitrarily small backward error, in nearly matrix multiplication time. \\
\indent We start by recalling the necessary theory from linear algebra. The generalized Schur decomposition of $(A,B)$, as introduced by Stewart \cite{Stewart_Schur}, takes the form 
\begin{equation}\label{eqn: schur_form}
    (A,B) = Q_L(T_A,T_B)Q_R^H
\end{equation}
for $Q_L,Q_R$ unitary and $T_A,T_B$ upper triangular. It is easy to see from \eqref{eqn: schur_form} that any set of leading columns of $Q_R$ span a right deflating subspace of $(A,B)$. Moreover, leading columns of $Q_L$ span the corresponding left deflating subspaces. When $(A,B)$ is regular, its eigenvalues are encoded as ratios of the diagonal entries of $T_A$ and $T_B$, with zeros on the diagonal of $T_B$ corresponding to eigenvalues at infinity. Eigenvectors, meanwhile, can be computed from \eqref{eqn: schur_form} via cheap triangular solves. \\
\indent To set the stage for a fast Schur decomposition, we first describe \textbf{EIG},  the main divide-and-conquer routine used in \cite{arXiv}. This algorithm divides the input pencil recursively, at each step searching for a line that splits the spectrum and subsequently computing bases for the corresponding pairs of deflating subspaces (the matrices $U_R^{(1)}, U_L^{(1)}, U_R^{(2)}$ and $U_L^{(2)}$ from \cref{section: intro}). Efficiency is derived from the fact that a split separating at least a fifth of the eigenvalues can be found at each step, which is a consequence of assuming that the input pencil is suitably well-behaved -- i.e., it satisfies a certain guarantee of \textit{pseudospectral shattering} as originally defined by Banks et al.\ \cite{banks2020pseudospectral}. As we discuss in more detail below, and was proved rigorously in \cite[Section 3]{arXiv}, such a guarantee can be obtained for arbitrary $(A,B)$ by applying small random perturbations to the matrices, allowing us to use \textbf{EIG} to approximately diagonalize any pencil (in the backward-error sense). \\
\indent Intuitively, we can adapt this approach to obtain the Schur form of $(A,B)$ by block upper triangularizing the pencil at each step. This requires only half the work of the diagonalization  case, as we need only compute $U_R^{(1)}, U_L^{(1)} \in {\mathbb C}^{n \times k}$. Indeed, if $U_R,U_L \in {\mathbb C}^{n \times n}$ are unitary matrices such that 
\begin{equation} \label{eqn: left_right_blocks}
    U_R = \begin{pmatrix} U_R^{(1)} & W_R \end{pmatrix}, \; \; U_L = \begin{pmatrix} U_L^{(1)} & W_L \end{pmatrix} 
\end{equation}
for some $W_R,W_L \in {\mathbb C}^{(n-k) \times n}$, then 
\begin{equation}\label{eqn: block_schur}
U_L^HAU_R = \begin{pmatrix} A_{11} & A_{12} \\
0 & A_{22} \end{pmatrix}, \; \; U_L^HBU_R = \begin{pmatrix} B_{11} & B_{12} \\ 0 & B_{22} \end{pmatrix}.
\end{equation}
This follows from the observation that the lower left block of $U_L^HAU_R$ is $W_L^HAU_R^{(1)}$; since $AU_R^{(1)}$ belongs to the left deflating subspace of $(A,B)$ spanned by the columns of $U_L^{(1)}$ -- and moreover $\text{range}(W_L)$ is the orthogonal complement of this deflating subspace since $U_L$ is unitary -- we have $W_L^HAU_R^{(1)} = 0$. The same argument implies $W_L^HBU_R^{(1)} = 0$. \\
\indent Before stating a divide-and-conquer Schur algorithm that leverages this observation, we pause to recap pseudospectral shattering. Let the grid $g = \text{grid}(z_0, \omega, s_1, s_2)$ be the boundary of the $s_1 \times s_2$ lattice in the complex plane consisting of ($\omega \times \omega$)-sized squares with lower left corner $z_0 \in {\mathbb C}$ (and grid lines parallel to the real/complex axis). Recalling \cref{defn:pencil_pseudospectrum}, we say that $\Lambda_{\epsilon}(A,B)$ is \textit{shattered} with respect to $g$ if the following conditions hold:\ (1) each eigenvalue of $(A,B)$ belongs to a unique grid box of $g$ and (2) $\Lambda_{\epsilon}(A,B) \cap g = \emptyset$. As mentioned above, this definition is originally due to Banks et al. \cite{banks2020pseudospectral}. \\
\indent The key insight of \cite{arXiv} -- which generalizes \cite{banks2020pseudospectral} -- is that small, Gaussian perturbations to $A$ and $B$ guarantee shattering with high probability for any input pencil and a certain choice of $\epsilon$ and $g$ (see \cite[Theorem 3.12]{arXiv}). Shattering positions us perfectly for divide-and-conquer; the grid provides a collection of viable splitting lines, each guaranteed to be minimally well-separated (in the sense of \cref{section: SIGN}) from a certain pseudospectrum of the perturbed problem, and we are also guaranteed that one of these lines splits off at least a fifth of the eigenvalues. Running divide-and-conquer over this shattering grid efficiently produces a diagonalization (or in our case a Schur decomposition) of the perturbed pencil, which can approximate the original problem $(A,B)$ provided the initial perturbation is small. 

\begin{algorithm}
\caption{Divide-and-Conquer Generalized Schur Decomposition (\textbf{SCHUR}) \\
\textbf{Input:} $n \in {\mathbb N}_{+}$, $A,B \in {\mathbb C}^{m \times m}$, $\epsilon > 0$, $C > 0$, $g \subset  \left\{ z : |\text{Re}(z)|, |\text{Im}(z)| < r \right\}$ an $s_1 \times s_2$ grid of box size $\omega$, $\eta > 0$ a desired accuracy, and $\theta \in (0,1)$ a failure probability.\\
\textbf{Requires:} $m \leq n$, $||A||_2, ||B||_2 \leq C$, and $\Lambda_{\epsilon}(A,B)$ shattered with respect to $g$. \\
\textbf{Output:} Unitary $Q_R,Q_L$ and an upper triangular pencil $(T_A,T_B)$ such that $||A - Q_LT_AQ_R^H||_2 \leq \eta$ and $||B - Q_LT_BQ_R^H||_2 \leq \eta $ with probability at least $1 - \theta$ (see \cref{thm: Schur_succeeds}).}\label{alg: Schur}
\begin{algorithmic}[1]
\If{$m=1$}
    \State $Q_R = Q_L = 1$; $T_A = A$; $T_B = B$
\Else
    \State $\zeta = 2 \left( \lfloor \log_2(\max \left\{s_1, s_2 \right\})  + 1 \rfloor\right)  $
    \State $\delta = \min \left\{ \frac{\theta}{2(\theta + 10n^6 \zeta)}, \sqrt{ \frac{\theta}{5}} \frac{1}{n^3C^2} \min \left\{ \frac{\epsilon^2}{800}, \frac{\eta^2}{96} \right\} \right\} $
    \State $p = \left\lceil \max \left\{ 7, \; -2 \log_2 \left( -\frac{1}{2} \log_2 \left( 1 - \frac{\epsilon}{r(1+r)C} \right) \right), \; 1 + \log_2 \left[ \frac{ \log_2 \left( \frac{\delta \pi \epsilon}{4 m \omega C + \delta \pi \epsilon} \right)}{\log_2 \left( 1 - \frac{\epsilon}{r(1+r)C}\right)} \right] \right\} \right\rceil$
    \State Choose a grid line $\text{Re}(z) = h$ of $g$
    \State $({\mathcal A}, {\mathcal B}) = (A - (h-1)B, A - (h+1)B)$
    \State $[A_p,B_p] = \textbf{IRS}({\mathcal A}, {\mathcal B}, p)$
    \State $[U,R_1,R_2, V] = \textbf{GRURV}(2, A_p + B_p, A_p, -1, 1)$
    \State $k = \# \left\{ i : \left| \frac{R_2(i,i)}{R_1(i,i)} \right| \geq \sqrt{\frac{\theta}{10 \zeta}}\frac{1 - \delta}{n^3} \right\}  $
    \If{$k < \frac{1}{5}m$ or $k > \frac{4}{5}m$} 
        \State Return to step 7 and choose a new grid line, executing a binary search if necessary. If this fails, \indent  \; \; search over horizontal grid lines $\text{Im}(z) = h$, this time setting $({\mathcal A}, {\mathcal B}) = (A - i(h-1)B, A - i(h+1)B)$.
    \Else 
        \State $U_R = \textbf{GRURV}(2,A_p+B_p,A_p,-1,1)$
        \State $[A_p,B_p] = \textbf{IRS}({\mathcal A}^H,{\mathcal B}^H,p)$
        \State $U_L = \textbf{GRURV}(2,A_p^H, (A_p+B_p)^H,1,-1)$.
        \vspace{2mm}
        \State $U_L^HAU_R = \begin{pmatrix} A_{11}& A_{12}\\ E & A_{22} \end{pmatrix}, \; \; U_L^HBU_R = \begin{pmatrix} B_{11} & B_{12} \\ F & B_{22} \end{pmatrix} $\Comment{$A_{11},B_{11} \in {\mathbb C}^{k \times k}$}
        \vspace{2mm}
        \State $g_R = \left\{ z \in g : \text{Re}(z) > h \right\}$ (or $g_R = \left\{ z \in g : \text{Im}(z) > h \right\}$)
        \State $[\widehat{Q}_R, \widehat{Q}_L, \widehat{T}_A, \widehat{T}_B] = \textbf{SCHUR}(n, A_{11}, B_{11}, \frac{4}{5} \epsilon, C, g_R, \frac{1}{\sqrt{3}} \eta, \theta)$
        \State $g_L = \left\{ z \in g : \text{Re}(z) < h \right\}$ (or $g_L = \left\{ z \in g : \text{Im}(z) < h \right\}$)
        \State $[\widetilde{Q}_R,  \widetilde{Q}_L, \widetilde{T}_A, \widetilde{T}_B] = \textbf{SCHUR}(n, A_{22}, B_{22}, \frac{4}{5} \epsilon, C, g_L, \frac{1}{\sqrt{3}} \eta, \theta)$
        \vspace{2mm}
        \State $Q_R = U_R \begin{pmatrix} \widehat{Q}_R & 0 \\ 0& \widetilde{Q}_R \end{pmatrix}, \; \; Q_L = U_L \begin{pmatrix} \widehat{Q}_L & 0 \\ 0 & \widetilde{Q}_L \end{pmatrix} $
        \vspace{2mm}
        \State $T_A = \begin{pmatrix} \widehat{T}_A & \widehat{Q}_{L}^HA_{12}\widetilde{Q}_R \\ 0 & \widetilde{T}_A \end{pmatrix}, \; \; T_B = \begin{pmatrix} \widehat{T}_B & \widehat{Q}_L^HB_{12}\widetilde{Q}_R \\ 0 & \widetilde{T}_B \end{pmatrix} $
        \vspace{2mm}
\EndIf
\EndIf
\State \Return $Q_R,Q_L,T_A,T_B$
\end{algorithmic}
\end{algorithm}

\indent  Of course, establishing a shattering grid first requires localizing the eigenvalues in a certain disk around the origin. For the generalized eigenvalue problem in particular, this requires some care. While in the single-matrix case we can simply scale the input matrix appropriately, a pencil $(A,B)$ can have arbitrarily large eigenvalues regardless of $||A||_2$ and $||B||_2$. Hence, the proof of shattering in  \cite{arXiv} requires not only Gaussian perturbations to $A$ and $B$ but a scaling factor that is polynomial in the size of the pencil (and is applied to only one of the matrices); that is, shattering is proved for the pencil $(\widetilde{A}, n^{\alpha} \widetilde{B})$, where $\alpha > 1$ is a constant and $\widetilde{A}$ and $\widetilde{B}$ are obtained from $A$ and $B$, respectively, by applying Gaussian perturbations. Accordingly, \textbf{EIG} assumes that the input matrices satisfy norm bounds in terms of $n^{\alpha}$. As discussed in \cite[Section 6]{arXiv}, this scaling factor is primarily a theoretical necessity and can likely be dropped in practice. \\
\indent  We are now ready to present \textbf{SCHUR} (\cref{alg: Schur}), a version of \textbf{EIG} that produces a Schur decomposition instead of a full diagonalization. As in the diagonalization case, \textbf{SCHUR} assumes that the input pencil comes with a guarantee of pseudospectral shattering for a corresponding $\epsilon$ and grid $g$. Nevertheless, this routine is somewhat more general. In particular,  we no longer build in the $n^{\alpha}$ scaling, instead taking an input $C > 0$ such that $||A||_2,||B||_2 \leq C$. While we will eventually set $C = 3n^{\alpha}$ to obtain a decomposition for arbitrary inputs, this allows the parameters of \textbf{SCHUR} to be relaxed if a better bound is available (e.g., for problems where shattering can be established without the polynomial scaling).  \\
\indent Throughout, \textbf{SCHUR} leverages the \textbf{GRURV} algorithm of Ballard et al.\ to (implicitly) compute rank-revealing factorizations of spectral projectors approximated by  \textbf{IRS}.\footnote{Repeated squaring could of course be replaced by any other inverse-free method considered in this paper, provided line 8 of \cref{alg: Schur} is adjusted accordingly.} A version of divide-and-conquer for Schur form built on these subroutines was originally proposed as \textbf{RGNEP} by Ballard, Demmel, and Dumitriu \cite{Ballard2010MinimizingCF}, though without the key step of pseudospectral shattering. Note that performance guarantees for \textbf{GRURV} can be found in \cite{grurv}. 

\begin{thm}\label{thm: Schur_succeeds}
Let $(A,B)$ and $g$ be a pencil and grid satisfying the requirements of \cref{alg: Schur} for corresponding $\epsilon > 0$ and $C > 0$. For any $\theta \in (0,1)$ and $\eta > 0$, suppose
$$[Q_R,Q_L,T_A,T_B] = \text{\normalfont \bf SCHUR}(n, A, B,  \epsilon, C, g, \eta, \theta) $$ 
in exact arithmetic. With probability at least $ 1- \theta$ this call is successful -- meaning its recursive procedure converges -- and moreover the outputs satisfy
$$ ||A - Q_LT_AQ_R^H||_2 \leq \eta \; \; \; \text{and} \; \; \; \; ||B - Q_LT_BQ_R^H||_2 \leq \eta. $$ 
\end{thm}
\begin{proof}
    We follow closely the proof of \cite[Theorem 4.8]{arXiv}, which exhibits success for \textbf{EIG}. Note that, in terms of our more general parameters, \textbf{EIG} takes $C = 3n^{\alpha}$ and $r = 5$, where $r$ bounds (in magnitude) the imaginary and real parts of the grid points of $g$. In repeating the analysis of \cite{arXiv}, we tighten a few inequalities where convenient and discard \textbf{EIG}'s second constraint on $p$ -- i.e., that $p \geq \log_2(\frac{105n^{\alpha}}{\epsilon} - 1)$ -- as it is redundant. \\
    \indent The first few steps of the proof carry over directly:\ pseudospectral shattering guarantees that the pencil $({\mathcal A}, {\mathcal B})$ in line 8 is minimally well-conditioned for repeated squaring and consequently that the choice of $p$ in line 6 ensures that \textbf{IRS} approximates projectors to accuracy $\delta$ in the spectral norm. Similarly, we know that a good eigenvalue split always exists, and -- by requiring that $k$ is computed correctly for every grid line checked -- \textbf{SCHUR} finds such a split with probability at least $ 1- \frac{\theta}{10n^4}$, provided $\delta \leq \frac{\theta}{2(\theta + 10n^6 \zeta)}$. \\
    \indent The first deviation from the diagonalization case occurs in lines 15-17, where \textbf{SCHUR} computes bases for one pair of deflating subspaces instead of two. Intuitively, we can account for this by applying the guarantees of \textbf{DEFLATE} (i.e., \cite[Algorithm 4]{arXiv}) once instead of twice. While this changes the parameters by (at most) a constant factor, we include the details for completeness. \\
    \indent Let $U_R = \begin{pmatrix} U_R^{(1)} & W_R \end{pmatrix}$ and $U_L = \begin{pmatrix} U_L^{(1)} & W_L \end{pmatrix}$ for $U_R^{(1)}, U_L^{(1)} \in {\mathbb C}^{m \times k}$. Taking $\nu = \sqrt{\frac{\theta}{5n^4}}$ in \cite[Theorem 4.7]{arXiv}, we are guaranteed that with probability at least $ 1 - \frac{2 \theta}{5n^4}$ there exist $\overline{U}_R, \overline{U}_L \in {\mathbb C}^{n \times k}$ whose columns span the corresponding true right/left deflating subspaces of $({\mathcal A}, {\mathcal B})$ and 
    \begin{equation}\label{eqn: grurv_error}
        \max \left\{||U_R^{(1)} - \overline{U}_R||_2, \; ||U_L^{(1)} - \overline{U}_L||_2 \right\} \leq \sqrt{8 n^2 \delta \sqrt{\frac{5k(m-k)}{\theta}}} \leq \sqrt{\sqrt{\frac{5}{\theta}} 8n^3 \delta}.
    \end{equation}
   The second constraint on $\delta$ therefore implies
    \begin{equation}\label{eqn: deflating_subspace_accuracy}
        \max \left\{ ||U_R^{(1)} - \overline{U}_R||_2, \; ||U_L^{(1)} - \overline{U}_L||_2 \right\} \leq \min \left\{\frac{\epsilon}{10C}, \frac{\eta}{2\sqrt{3}C} \right\},
    \end{equation}
    again with probability at least $ 1- \frac{2 \theta}{5n^4}$. As we will see, these bounds guarantee, respectively, that shattering is preserved as we recur and moreover that the error in the resulting factorization can be controlled. A simple union bound implies that at any step of the recursive procedure a good dividing line is found, $k$ is computed correctly, and \eqref{eqn: deflating_subspace_accuracy} holds with probability at least $ 1 - \frac{\theta}{2n^4}$. \\
    \indent Before moving on, we pause to consider how bounds like $||U_R^{(1)} - \overline{U}_R||_2 \leq \mu$ and  $||U_L^{(1)} - \overline{U}_L||_2 \leq \mu$ can be extended to cover $||E||_2$ and $||F||_2$ in line 18. The heuristic here is fairly straightforward:\ if the deflating subspaces are approximated accurately enough, $E$ and $F$ should be nearly zero. With this in mind, let $U_R^{(1)} = \overline{U}_R + \Delta_1$ for some $\Delta_1 \in {\mathbb C}^{m \times k}$ with $||\Delta_1||_2 \leq \mu $. In this case,
    \begin{equation}\label{eqn: bound_E_1}
    E = W_L^HAU_R^{(1)} = W_L^HA(\overline{U}_R + \Delta_1) = W_L^HA\overline{U}_R + W_L^HA\Delta_1.
    \end{equation}
    Now $\overline{U}_R$ contains a basis for the right deflating subspace of $(A,B)$ corresponding to $\overline{U}_L$, so $A\overline{U}_R = \overline{U}_LX$ for some $X \in {\mathbb C}^{k \times k}$. Applying this to \eqref{eqn: bound_E_1} and further letting $\overline{U}_L = U_L^{(1)} + \Delta_2$ for another $\Delta_2 \in {\mathbb C}^{m \times k}$ satisfying $||\Delta_2||_2 \leq \mu$, we have
    \begin{equation}\label{eqn: bound_E_2}
        \aligned
        E &= W_L^H\overline{U}_LX + W_L^HA\Delta_1 = W_L^H(U_L^{(1)} + \Delta_2)X +W_L^HA \Delta_1= W_L^H\Delta_2 X + W_L^H A \Delta_1,
        \endaligned
    \end{equation}
    where the last equality follows from $W_L^H U_L^{(1)} = 0$. Hence, we conclude
    \begin{equation}
        ||E||_2 \leq ||W_L^H \Delta_2 X||_2 + || W_L^H A \Delta_1||_2 \leq \mu (||X||_2 + ||A||_2).
    \end{equation}
    We can remove $||X||_2$ from this bound by noting $X = \overline{U}_L^HA\overline{U}_R$ and therefore $||X||_2 \leq ||A||_2$. Thus, $||E||_2 \leq 2 \mu ||A||_2$. A similar argument yields $||F||_2 \leq 2 \mu ||B||_2$. \\
    \indent We now show that the first bound in \eqref{eqn: deflating_subspace_accuracy} implies that the recursive calls in lines 20 and 22 are valid, meaning the inputs satisfy the listed requirements. Since multiplying by blocks of unitary matrices cannot violate the norm constraints, this reduces to showing that $\Lambda_{4 \epsilon /5}(A_{11}, B_{11})$ and $\Lambda_{4 \epsilon /5}(A_{22}, B_{22})$ are shattered with respect to $g_R$ and $g_L$, respectively. For $(A_{11}, B_{11})$ we can simply repeat the argument used in the proof of \cite[Theorem 4.8]{arXiv}; that is, $||A_{11} - \overline{U}_L^HA\overline{U}_R||_2$ and  $||B_{11} - \overline{U}_L^HB\overline{U}_R||_2$ can both be bounded by  $\frac{\epsilon}{5}$, in which case $\Lambda_{4\epsilon/5}(A_{11}, B_{11})$ is shattered with respect to $g$ with each eigenvalue of $(A_{11}, B_{11})$ sharing a unique grid box with an eigenvalue of $(\overline{U}_L^HA\overline{U}_R, \overline{U}_L^HB \overline{U}_R)$. Hence, we also know that $\Lambda(A_{11}, B_{11}) \subset g_R$.\\
    \indent Unfortunately, the same approach cannot be used for $(A_{22}, B_{22})$, as the columns of $W_L$ and $W_R$ do \textit{not} approximately span deflating subspaces in general. Instead, we'll need to establish the following two items separately:
    \begin{enumerate}
        \item $\Lambda_{4 \epsilon /5}(A_{22}, B_{22})$ is shattered with respect to $g$.
        \item $(A_{22}, B_{22})$ has all of its eigenvalues in $g_L$.
    \end{enumerate}
    
    \indent We'll tackle shattering first. Suppose $z \in \Lambda_{4\epsilon /5}(A_{22}, B_{22})$. In this case,  $\sigma_{m-k}(A_{22} - zB_{22}) \leq \frac{4}{5}\epsilon (1 + |z|)$ (see e.g., \cite{Fraysse96spectralportraits}). Since $A_{22} - z B_{22}$ is the lower right $(m-k) \times (m-k)$ block of $U_L^HAU_R - z U_L^HBU_R$, where $U_L$ and $U_R$ are unitary, it is easy to see that
    \begin{equation} \label{eqn: extend_sv_bound}
        \sigma_m(A - zB) = \sigma_m(U_L^HAU_R - z U_L^HBU_R) \leq \sqrt{\sigma_{m-k}(A_{22} - z B_{22})^2 + ||E - zF||_2^2}.
    \end{equation}
    Extending \eqref{eqn: deflating_subspace_accuracy} as outlined above, recalling $||A||_2, ||B||_2 \leq C$, we can bound $||E - zF||_2$ as follows:
    \begin{equation}\label{eqn: bound_E-zF}
        ||E - zF||_2 \leq ||E||_2 + |z|||F||_2 \leq \frac{\epsilon}{5}(1+|z|).
    \end{equation}
    Applying this to \eqref{eqn: extend_sv_bound} alongside our bound for $\sigma_{m-k}(A_{22} - zB_{22})$ yields
    \begin{equation}\label{eqn: simplify_sv_bound}
        \sigma_m(A-zB) \leq \sqrt{\frac{17}{25}} \epsilon (1 + |z|) \leq \epsilon (1 + |z|),
    \end{equation}
    which implies $z \in \Lambda_{\epsilon}(A,B)$. Thus, $\Lambda_{4 \epsilon / 5}(A_{22},B_{22}) \subseteq \Lambda_{\epsilon}(A,B)$; since the latter is assumed to be shattered with respect to $g$ on input, this proves shattering for $\Lambda_{4\epsilon/ 5}(A_{22}, B_{22})$. \\
    \indent It remains to show that the eigenvalues of $(A_{22}, B_{22})$ lie in $g_L$. To do this, consider the $m \times m$ matrices
    \begin{equation}\label{eqn: zero_out_EF}
        M_1 = \begin{pmatrix} A_{11} & A_{12} \\
        0 & A_{22} \end{pmatrix}, \; \; M_2 = \begin{pmatrix} B_{11} & B_{12} \\
        0 & B_{22} \end{pmatrix}, 
    \end{equation}
    which are obtained from $U_L^HAU_R$ and $U_L^HBU_R$, respectively, by zeroing out $E$ and $F$. Since $(M_1, M_2)$ is block upper triangular, its spectrum is $\Lambda(A_{11}, B_{11}) \cup \Lambda(A_{22},B_{22})$. Moreover, we have
    \begin{equation} \label{eqn: bound_error_in_removing_E}
        U_L^HAU_R - M_1 = \begin{pmatrix} 0 & 0 \\ E & 0 \end{pmatrix} \; \; \Longrightarrow \; \; ||U_L^HAU_R - M_1||_2 \leq ||E||_2 \leq \frac{\epsilon}{5}
    \end{equation}
    and similarly $||U_L^HBU_R - M_2||_2 \leq \frac{\epsilon}{5}$. Now $\Lambda_{\epsilon}(U_L^HAU_R, U_L^HBU_R)$ is shattered with respect to $g$, so by \cite[Lemma 3.14]{arXiv} we know that each eigenvalue of $(M_1,M_2)$ shares a grid box with an eigenvalue of $(U_L^HAU_R, U_L^HBU_R)$.\footnote{This again follows from the fact that $\Lambda_{\epsilon}(A,B)$ is shattered with respect to $g$, as shattering is both unchanged when $A$ and $B$ are multiplied by unitary matrices and additionally implies that $(A,B)$ -- and therefore also $(U_L^HAU_R, U_L^HBU_R)$ -- is regular.} But $\Lambda(U_L^HAU_R, U_L^HBU_R) = \Lambda(A,B)$ and $(A,B)$ has exactly $k$ eigenvalues in $g_R$ and $m-k$ in $g_L$ (recall that we are assuming $k$ has been computed accurately in line 11). Since $\Lambda(M_1,M_2) = \Lambda(A_{11}, B_{11}) \cup \Lambda(A_{22},B_{22})$ and we showed above that $(A_{11}, B_{11})$ has all of its eigenvalues in $g_R$, this implies that $(A_{22},B_{22})$ must have all of its eigenvalues in $g_L$.  \\
    \indent Combining points one and two above, we conclude that $\Lambda_{4 \epsilon / 5}(A_{22},B_{22})$ is shattered with respect to $g_L$ and therefore that the call to \textbf{SCHUR} in line 22 is valid. As a result, we have proved success for one step of the recursive procedure under the conditions that (1) a good dividing line is found, (2) the corresponding value of $k$ is accurate, and (3) the upper bound \eqref{eqn: deflating_subspace_accuracy} holds. Since these occur simultaneously (at any step) with probability at least $1 - \frac{\theta}{2n^4}$, and moreover the recursive tree of \textbf{SCHUR} has depth at most $\log_{5/4}(n)$, the algorithm fails to converge with probability at most $2 \cdot 2^{\log_{5/4}(n)} \frac{\theta}{2n^4} \leq \theta$. \\
    \indent To complete the proof, we derive the error bounds $||A - Q_LT_AQ_R^H||_2 \leq \eta$ and $||B - Q_LT_BQ_R^H||_2 \leq \eta$ under the same assumptions. We do this inductively. The base case is $m=1$, in which case the decomposition produced by \textbf{SCHUR} is exact. Consider now any other step of the procedure. 
    We have
    \begin{equation}\label{eqn: Schur_error_1}
        \aligned 
        ||A - Q_LT_AQ_R^H ||_2 &= ||U_L^HAU_R - U_L^HQ_LT_AQ_R^HU_R||_2 \\
        & = \left| \left| \begin{pmatrix} A_{11} & A_{12} \\
        E & A_{22} \end{pmatrix} - \begin{pmatrix} \widehat{Q}_L & 0 \\ 0 & \widetilde{Q}_L \end{pmatrix} \begin{pmatrix} \widehat{T}_A & \widehat{Q}_{L}^HA_{12}\widetilde{Q}_R \\ 0 & \widetilde{T}_A \end{pmatrix} \begin{pmatrix} \widehat{Q}_R^H & 0 \\ 0 & \widetilde{Q}_R^H \end{pmatrix} \right| \right|_2 \\
        & = \left| \left| \begin{pmatrix} A_{11} - \widehat{Q}_L \widehat{T}_A \widehat{Q}_R^H & 0 \\
        E & A_{22} - \widetilde{Q}_L \widetilde{T}_A \widetilde{Q}_R^H \end{pmatrix} \right| \right|_2 \\
        & \leq \left( ||A_{11} - \widehat{Q}_L \widehat{T}_A \widehat{Q}_R^H||_2^2 + ||A_{22} - \widetilde{Q}_L \widetilde{T}_A \widetilde{Q}_R^H||_2^2 + ||E||_2^2 \right)^{1/2}.
        \endaligned
    \end{equation}
    By our induction hypothesis we can assume that $||A_{11} - \widehat{Q}_L \widehat{T}_A \widehat{Q}_R^H||_2$ and $||A_{22} - \widetilde{Q}_L\widetilde{T}_A \widetilde{Q}_R||_2$ are both bounded by $\frac{\eta}{\sqrt{3}}$. Since the second piece of \eqref{eqn: deflating_subspace_accuracy} implies the same bound for $||E||_2$, we therefore conclude 
    \begin{equation}\label{eqn: Schur_error_2}
        ||A - Q_LT_AQ_R^H||_2 \leq \sqrt{ \left( \frac{\eta}{\sqrt{3}}\right)^2 + \left( \frac{\eta}{\sqrt{3}}\right)^2 + \left( \frac{\eta}{\sqrt{3}}\right)^2} = \eta.
    \end{equation}
    Repeating this argument yields the corresponding bound $||B - Q_LT_BQ_R^H||_2 \leq \eta$.
\end{proof}

As a corollary, we obtain an algorithm that can compute an accurate Schur decomposition of any pencil in nearly matrix multiplication time (with high probability). Once again, note that -- like all of the methods considered in this paper -- this algorithm is entirely inverse-free. 

\begin{cor}\label{cor: Schur_in_nmmt}
    Let $(A,B) \in {\mathbb C}^{n \times n} \times {\mathbb C}^{n \times n}$ be any pencil with $||A||_2,||B||_2 \leq 1$ and let $ \varepsilon < 1$. There exists an exact-arithmetic, randomized algorithm that takes as inputs $(A,B)$ and $\varepsilon$ and -- in $O(\log^2(\frac{n}{\varepsilon})T_{\text{\normalfont MM}}(n))$ operations -- produces $Q_R,Q_L,T_A,T_B \in {\mathbb C}^{n \times n}$ such that $Q_R$ and $Q_L$ are unitary, $T_A$ and $T_B$ are upper triangular, and 
    $$ ||A - Q_LT_AQ_R^H||_2 \leq \varepsilon \; \; \; \text{and} \; \; \; \; ||B - Q_LT_BQ_R^H||_2 \leq \varepsilon $$
    with probability at least $1 - O(\frac{1}{n})$.
\end{cor}
\begin{proof}
    Consider the following algorithm, whose only inputs are $(A,B)$ and $\varepsilon$.
    \begin{enumerate}
        \item Set parameters: $\gamma = \frac{\varepsilon}{8}$; $\alpha = \frac{\lceil 2 \log_n(1/\gamma)+3 \rceil}{2}$; $\epsilon = \frac{\gamma^5}{64n^{\frac{11 \alpha + 25}{3}} + \gamma^5}$; $\omega = \frac{\gamma^4}{4}n^{- \frac{8\alpha + 13}{3}}$.
        \item Draw two independent Ginibre matrices $G_1,G_2 \in {\mathbb C}^{n \times n}$ and let $(\widetilde{A}, \widetilde{B}) = (A + \gamma G_1, B + \gamma G_2)$. 
        \item Draw $z$ uniformly at random from the box of side length $\omega$ with bottom left corner $-4-4i$ and let $g = \text{grid}(z, \omega, \lceil 8/\omega \rceil, \lceil 8/ \omega \rceil)$.
        \item $[Q_R, Q_L, T_A, T_B] = \textbf{SCHUR}(n, \widetilde{A}, n^{\alpha} \widetilde{B}, \epsilon, 3n^{\alpha}, g, \frac{\varepsilon}{2}, \frac{1}{n})$.
        \item Rescale: $T_B = n^{-\alpha} T_B$.
        \item Output $Q_R,Q_L,T_A$ and $T_B$.
    \end{enumerate}
    For the listed values of $\alpha$, $\epsilon$, and $\omega$, \cite[Theorem 3.12]{arXiv} guarantees that $||\widetilde{A}||_2 \leq 3$, $||n^{\alpha} \widetilde{B}||_2 \leq 3n^{\alpha}$, and $\Lambda_{\epsilon}(\widetilde{A}, n^{\alpha} \widetilde{B})$ is shattered with respect to $g$ with probability at least $1 - O(\frac{1}{n})$. In this case the call to \textbf{SCHUR} is valid, and moreover (by \cref{thm: Schur_succeeds}) it succeeds with probability at least $ 1- \frac{1}{n}$, achieving 
    \begin{equation}\label{eqn: accuracy_on_perturbed}
        ||\widetilde{A} - Q_LT_AQ_R^H||_2 \leq \frac{\varepsilon}{2} \; \; \text{and} \; \; ||n^{\alpha} \widetilde{B} - Q_LT_BQ_R^H||_2 \leq \frac{\varepsilon}{2}
    \end{equation}
    on exit (i.e., before rescaling $T_B$). Consequently,
    \begin{equation}\label{eqn: Schur_error_A}
        ||A - Q_LT_AQ_R^H||_2 \leq ||A - \widetilde{A}||_2 + ||\widetilde{A} - Q_LT_AQ_R^H||_2 \leq \frac{\varepsilon}{2} + \frac{\varepsilon}{2} \leq \varepsilon
    \end{equation}
    and similarly, again for the initial $T_B$,
    \begin{equation}\label{eqn: schur_error_B}
        ||B - Q_L(n^{-\alpha}T_B)Q_R^H||_2 \leq ||B - \widetilde{B}||_2 + || \widetilde{B} -Q_L(n^{-\alpha} T_B)Q_R^H||_2 \leq \frac{\varepsilon}{2} + \frac{\varepsilon}{2n^{\alpha}} = \varepsilon.
    \end{equation}
    Here, we make use of the fact that $||A - \widetilde{A}||_2 = \gamma||G_1||_2$ and $||B - \widetilde{B}||_2 = \gamma ||G_2||_2$, where we can assume $||G_1||_2, ||G_2||_2 \leq 4$.\footnote{The probability that one of $||G||_1$ and $||G||_2$ is larger than four is included in the failure probability of shattering.} Hence, a final union bound implies that \textbf{SCHUR} succeeds, and we achieve the desired accuracy, with probability at least $1 - O(\frac{1}{n})$. Since this algorithm does less work than its diagonalization counterpart, the complexity $O(\log^2(\frac{n}{\varepsilon})T_{\text{MM}}(n))$ follows from \cite[Proposition 5.3]{arXiv}. 
\end{proof}

\indent In the proof of \cref{cor: Schur_in_nmmt}, \textbf{SCHUR} is applied to a pencil $(A,B)$ where $A$ and $B$ are minimally well-conditioned with high probability. For such inputs, we can simplify \cref{alg: Schur} by replacing lines 16 and 17 with one of the  following:
\begin{equation}\label{eqn: simplify_Schur}
    [U_L,R] = \textbf{QR}(BU_R) \; \; \; \text{or} \; \; \; [U_L,R] = \textbf{QR}(AU_R).
\end{equation}
In other words, we can leverage the fact that $A$ and $B$ map right deflating subspaces to the corresponding left ones, trading an expensive call to \textbf{IRS} for a simple QR factorization. This results in significant savings, as QR can be done in matrix multiplication time while \textbf{IRS} requires $O(\log ( \frac{n}{\varepsilon}) T_{\text{MM}}(n))$ operations for the parameters chosen in the proof of \cref{cor: Schur_in_nmmt}. Moreover, when a bound on $\kappa_2(B)$ or $\kappa_2(A)$ is available, we can easily describe how error in $U_R$ propagates to $U_L$. We provide a sketch for $[U_L,R] = \textbf{QR}(BU_R)$ below. \\
\indent As in the proof of \cref{thm: Schur_succeeds}, let $U_R = \begin{pmatrix} U_R^{(1)} & W_R \end{pmatrix}$ and $U_L = \begin{pmatrix} U_L^{(1)} & W_L \end{pmatrix}$ for $U_R^{(1)}, U_L^{(1)} \in {\mathbb C}^{m \times k}$ and suppose that $U_R^{(1)}$ has been computed to some accuracy $\mu$ -- i.e., there exists a true $\overline{U}_R \in {\mathbb C}^{m \times k}$ (whose columns span the associated true right deflating subspace) such that $||U_R^{(1)} - \overline{U}_R||_2 \leq \mu$. Further let $R_{11}$ be the upper left $k \times k$ block of $R$. In this case, $BU_R^{(1)} = U_L^{(1)} R_{11}$ and  $||BU_R^{(1)} - B\overline{U}_R||_2 \leq \mu ||B||_2$. Hence, classical perturbation theory for reduced QR (e.g., \cite{10.1007/BF01931293}) implies that for $\mu$ sufficiently small there exists a (unique) reduced QR factorization 
\begin{equation}\label{eqn: nearby_QR}
    B\overline{U}_R = \overline{U}_L(R_{11} + \Delta R)
\end{equation}
such that\footnote{The $\sqrt{k}$ dependence here can actually be improved to $\log(k)$ \cite[Chapter 3]{My_thesis}.}
\begin{equation}\label{eqn: nearby_QR2}
    ||U_L^{(1)} - \overline{U}_L||_2 \leq O (\sqrt{k}) \cdot \kappa_2(BU_R^{(1)}) \cdot \frac{\mu ||B||_2}{||BU_R^{(1)}||_2} \leq O(\sqrt{k}) \kappa_2(B) \mu
\end{equation}
Clearly, the columns of $\overline{U}_L$ span the left deflating subspace corresponding to $\overline{U}_R$. Thus, we can simply tighten $\mu$ by $O(\sqrt{k}) \kappa_2(B)$ to obtain a similar accuracy guarantee for the left deflating subspace. In terms of the parameters of \textbf{SCHUR}, this would correspond to a slightly more restrictive choice of $\delta$ in line 5. \\
\indent Of course, this requires a bound on $\kappa_2(B)$ (or similarly $\kappa_2(A)$). While a guarantee of pseudospectral shattering does imply that $B$ is nonsingular, it does not necessarily guarantee that $B$ is particularly well-conditioned; if it is not, obtaining $U_L$ from $U_R$ may actually be more costly than repeating \textbf{IRS} and \textbf{GRURV}. In an effort to further maintain the generality of \cref{alg: Schur}, we therefore compute $U_L$ via \textbf{IRS}. Nevertheless, \eqref{eqn: simplify_Schur} is likely to be the cheaper option in practice.
\begin{remark}
    While the parameters chosen in \cref{alg: Schur} and the proof of \cref{cor: Schur_in_nmmt} seem restrictive, they are in fact much nicer than the diagonalization case, as we might expect. Compare in particular the dependence of $\delta$ on $\eta$ in \textbf{SCHUR} versus $\beta$ in \textbf{EIG}. Since most of the work in both algorithms is spent locating good dividing lines, this is the real advantage -- in terms of efficiency -- of going to Schur form instead of complete diagonalization. Note that recent work of Shah \cite{fast_herm_diag} opens the possibility of further relaxing parameters in both algorithms, though not enough to change their final asymptotic complexities. 
\end{remark}

\section{Proof of Theorem 2.7}\label{section: appendix}
In this appendix, we provide a full proof of \cref{thm: main_IRS_bound}. 
\subsection{Technical Lemmas}
We begin by stating a handful of technical lemmas. The first two are straightforward consequences of the black-box, finite-precision assumptions from \cref{section: IRS}. Here, $\textbf{MM}(\cdot )$ and $\textbf{QR}(\cdot)$ are again the floating-point algorithms covered by \cref{MM assumption} and \cref{QR assumption}, respectively. 
\begin{lem}\label{lem: finite_addition}
Let $A,B \in {\mathbb C}^{n \times n}$. If $C$ is the floating-point sum of $A$ and $B$ then 
$$||C - (A+B)||_2 \leq \sqrt{n} {\bf u} ||A + B||_2.$$
\end{lem}
\begin{proof}
    By \eqref{eqn: floating point} each entry $C_{ij}$ of $C$ satisfies $C_{ij} = (A+B)_{ij}(1 + \Delta_{ij})$ for some $|\Delta_{ij}| \leq {\bf u}$. Consequently,
    \begin{equation}
        [C - (A+B)]_{ij} = (A+B)_{ij} \Delta_{ij}
        \label{eqn: entry_wise_diff}
    \end{equation}
    and therefore
    \begin{equation}
        ||C - (A+B)||_F^2 \leq \sum_{ij} |(A+B)_{ij}\Delta_{ij}|^2 \leq {\bf u}^2 \sum_{ij} |(A+B)_{ij}|^2 = {\bf u}^2 ||A+B||_F^2.
        \label{eqn: addition_error_bound}
    \end{equation}
    We complete the proof by noting $||C - (A+B)||_2 \leq ||C - (A+B)||_F$ and $||A+B||_F \leq \sqrt{n} ||A+B||_2$.
\end{proof}

\begin{lem}\label{lem: QR_backward_R} 
Let $[Q,R] = \normalfont{\textbf{QR}(A)}$ for 
$$ A = \begin{pmatrix} A_1 \\ A_2 \end{pmatrix}, \; \; Q = \begin{pmatrix} Q_{11} & Q_{12} \\
Q_{21} & Q_{22} \end{pmatrix}, \; \text{and} \; \; R = \begin{pmatrix} R' \\0 \end{pmatrix}, $$
where all blocks are $n \times n$. Define the following matrices:
\begin{enumerate}
    \item $\widetilde{R}$ -- the floating-point sum of $\normalfont{\textbf{MM}}(Q_{11}^H,A_1)$ and $\normalfont{\textbf{MM}}(Q_{21}^H,A_2)$.
    \item $E$ -- the floating-point sum of $\normalfont{\textbf{MM}}(Q_{12}^H,A_1)$ and $\normalfont{\textbf{MM}}(Q_{22}^H, A_2)$.
\end{enumerate}
If $\mu_{\text{\normalfont QR}}(n) {\bf u}, \mu_{\text{\normalfont MM}}(n) {\bf u} \leq 1$, then
$$ ||\widetilde{R} - R||_2, ||E||_2 \leq 4\left(\mu_{\text{\normalfont QR}}(n) + \mu_{\text{\normalfont MM}}(n) + 2 \sqrt{n} \right){\bf u}||A||_2.$$
\end{lem}
\begin{proof}
    By \cref{QR assumption}, there exist matrices $\widehat{A} \in {\mathbb C}^{2n \times n}$ and $\widehat{Q} \in {\mathbb C}^{2n \times 2n}$ such that $\widehat{Q}$ is truly unitary, $\widehat{A} = \widehat{Q}R$, $||Q - \widehat{Q}||_2 \leq \mu_{\text{QR}}(n) {\bf u}$, and $||A - \widehat{A}||_2 \leq \mu_{\text{QR}}(n) {\bf u} ||A||_2$. Let
    \begin{equation}
        \widehat{A} = \begin{pmatrix} \widehat{A}_1 \\ \widehat{A}_2 \end{pmatrix} \; \text{and} \; \; \widehat{Q} = \begin{pmatrix} \widehat{Q}_{11} & \widehat{Q}_{12} \\ \widehat{Q}_{21} & \widehat{Q}_{22} \end{pmatrix},
        \label{eqn: block_A_and_Q}
    \end{equation}
    where again all blocks are $n \times n$. Consider first $\textbf{MM}(Q_{11}^H, A_1)$. Applying \cref{MM assumption} we observe
    \begin{equation}
        \aligned 
        || \textbf{MM}(Q_{11}^H, A_1) - \widehat{Q}_{11}^H \widehat{A}_1||_2 &\leq || \textbf{MM}(Q_{11}^H, A_1) - Q_{11}^H A_1||_2 + ||Q_{11}^HA_1 - \widehat{Q}_{11}^H \widehat{A}_1 ||_2 \\
        & \leq \mu_{\text{MM}}(n) {\bf u} ||Q_{11}||_2 ||A_1||_2 + ||Q_{11}^HA_1 - \widehat{Q}_{11}^HA_1||_2 + ||\widehat{Q}_{11}^H A_1 - \widehat{Q}_{11}^H \widehat{A}_1||_2 \\
        & \leq \mu_{\text{MM}}(n) {\bf u} ||Q||_2||A||_2 + ||Q - \widehat{Q}||_2 ||A||_2 + ||\widehat{Q}||_2 ||A - \widehat{A}||_2.
        \endaligned 
        \label{eqn: expand_with_MM_guarantee}
    \end{equation}
    Since $||\widehat{Q}||_2 = 1$ and therefore $||Q||_2 \leq ||\widehat{Q}||_2 + \mu_{\text{QR}}(n) {\bf u} = 1 + \mu_{\text{QR}}(n) {\bf u}$, \eqref{eqn: expand_with_MM_guarantee} implies 
    \begin{equation}
        \aligned 
        ||\textbf{MM}(Q_{11}^H,A_1) - \widehat{Q}_{11}^H \widehat{A}_1 ||_2 & \leq \left[ 2 \mu_{\text{QR}}(n) + \mu_{\text{MM}}(n)(1 + \mu_{\text{QR}}(n) {\bf u}) \right] {\bf u} ||A||_2\\
        &\leq 2 (\mu_{\text{QR}}(n) + \mu_{\text{MM}}(n)) {\bf u}||A||_2.
        \endaligned 
        \label{eqn: simplified_block_bound}
    \end{equation}
    Repeating this argument, swapping blocks accordingly, we obtain the same result for $||\textbf{MM}(Q_{21}^H,A_2) - \widehat{Q}_{21}^H \widehat{A}_2||_2$, $||\textbf{MM}(Q_{12}^H,A_1) - \widehat{Q}_{12}^H\widehat{A}_1||_2$, and $||\textbf{MM}(Q_{22}^H,A_2) - \widehat{Q}_{22}^H \widehat{A}_2||_2$. To now bound $||\widetilde{R} - R'||_2$, note that 
    \begin{equation}
        R' = \widehat{Q}_{11}^H \widehat{A}_1 + \widehat{Q}_{21}^H \widehat{A}_2
        \label{eqn: R_prime}
    \end{equation}
    since $\widehat{A} = \widehat{Q}R$. Consequently, 
    \begin{equation}
        \aligned 
        ||\widetilde{R} - R'||_2 \leq ||\widetilde{R} - &(\textbf{MM}(Q_{11}^H, A_1) + \textbf{MM}(Q_{21}^H,A_2))||_2 \\
         & + ||\textbf{MM}(Q_{11}^H,A_1) - \widehat{Q}_{11}^H \widehat{A}_1||_2 + ||\textbf{MM}(Q_{21}^H,A_2) - \widehat{Q}_{21}^H \widehat{A}_2||_2.
        \endaligned 
        \label{eqn: R_bound}
    \end{equation}
    By \cref{lem: finite_addition}, we can bound the first term by $\sqrt{n} {\bf u} ||\textbf{MM}(Q_{11}^H,A_1) + \textbf{MM}(Q_{21}^H,A_2)||_2$, where
    \begin{equation}
        \aligned
        ||\textbf{MM}(Q_{11}^H,A_1) + \textbf{MM}(Q_{21}^H,A_2)||_2 &\leq ||\textbf{MM}(Q_{11}^H,A_1)||_2 + ||\textbf{MM}(Q_{21}^H,A_2)||_2 \\
        & \leq ||Q_{11}^HA_1||_2 + \mu_{\text{MM}}(n) {\bf u}||Q_{11}||_2||A_1||_2  \\
        & \; \; \; \; \; \; \;  + ||Q_{21}^HA_2||_2 + \mu_{\text{MM}}(n) {\bf u} ||Q_{21}||_2||A_2||_2 \\
        & \leq 2(1 + \mu_{\text{MM}}(n) {\bf u})(1 + \mu_{\text{QR}}(n) {\bf u}) ||A||_2 \\
        & \leq 8 ||A||_2.
        \endaligned
        \label{eqn: norm_bounds}
    \end{equation}
    Applying this to \eqref{eqn: R_bound} alongside \eqref{eqn: simplified_block_bound} yields
    \begin{equation}
        ||\widetilde{R} - R'||_2 \leq 4 \left(\mu_{\text{QR}}(n) + \mu_{\text{MM}}(n) + 2 \sqrt{n}\right){\bf u}||A||_2.
        \label{eqn: final_R_bound}
    \end{equation}
    We obtain the same bound for $||E||_2$ by repeating this argument with $\textbf{MM}(Q_{12}^H,A_1)$ and $\textbf{MM}(Q_{22}^H,A_2)$ and noting that $\widehat{Q}_{12}^H\widehat{A}_1 + \widehat{Q}_{22}^H \widehat{A}_2 = 0$.
\end{proof}
\indent Next, we state a pair of lemmas due to Malyshev \cite[Lemmas 4.1 and 4.2]{Malyshev1}.
\begin{lem}[Malyshev \cite{Malyshev1}]\label{lem: Malyshev_one}
Suppose $R \in {\mathbb C}^{m \times m}$ is nonsingular and $E \in {\mathbb C}^{n \times m}$ for $m \geq n$. There exists a matrix $S \in {\mathbb C}^{(m+n) \times (m+n)}$ such that 
\begin{enumerate}[itemsep = 0em]
    \item $(I+S) \binom{R}{E} = \binom{R'}{0}$.
    \item $(I+S)^H(I+S) = I$
    \item $||S||_2 \leq ||E R^{-1}||_2 \leq ||E||_2 ||R^{-1}||_2$.
\end{enumerate}
\end{lem}
\begin{proof}
    We take the opportunity to correct a small error in Malyshev's proof. Define
    \begin{equation}
        \widetilde{S} = \begin{pmatrix} 0 & (ER^{-1})^H \\
        -ER^{-1} & 0 \end{pmatrix}.
        \label{eqn: malyshev_proof_1}
    \end{equation}
    Then the matrix
    \begin{equation}
        S = \begin{pmatrix} [I + (ER^{-1})^HER^{-1}]^{-1/2} & 0 \\
        0 & [I+ER^{-1}(ER^{-1})^H]^{-1/2} \end{pmatrix} (I + \widetilde{S}) - I
        \label{eqn: malyshev_proof_2}
    \end{equation}
    satisfies the listed requirements.
\end{proof}

\begin{lem}[Malyshev \cite{Malyshev1}]\label{lem: Malyshev_two}
Let $A \in {\mathbb C}^{m \times n}$ be full rank and suppose 
$$ A = Q_1 \begin{pmatrix} K_1 & L_1 \\ 0 & M_1 \end{pmatrix} = Q_2 \begin{pmatrix} K_2 & L_2 \\ 0 & M_2 \end{pmatrix} $$
for $Q_1, Q_2 \in {\mathbb C}^{m \times m}$ unitary, $K_1,K_2 \in {\mathbb C}^{k \times k}$ nonsingular, and $M_1,M_2 \in {\mathbb C}^{(m-k) \times (n-k)}$ full rank. Then there exist unitary matrices $P \in {\mathbb C}^{k \times k}$ and $Q \in {\mathbb C}^{(n-k) \times (n-k)}$ such that $K_2 = PK_1$, $L_2 = PL_1$, and $M_2 = QM_1$. 
\end{lem}
\subsection{Detailed Proof}
We are now ready to bound error in \textbf{IRS}. To simplify the analysis, we will not track the individual polynomials $\mu_{\text{MM}}(n)$ and $\mu_{\text{QR}}(n)$, instead working with a ``general polynomial"
\begin{equation}
    \mu(n) = \max \left\{ \mu_{\text{MM}}(n), \mu_{\text{QR}}(n), \sqrt{n} \right\}
    \label{eqn: generic_poly}
\end{equation}
and the associated quantity $\tau = \mu(n) {\bf u}$. Throughout, we can think of $\tau$ as small, corresponding to a choice ${\bf u} < \mu(n)^{-1}$. \\
\indent Before proceeding with our main argument we state one final lemma, which bounds norm growth in repeated squaring. Because finite-precision \textbf{IRS} repeatedly multiplies the inputs by pieces of nearly unitary matrices, we expect that norms should grow by (at most) small constants.  Here, $\widetilde{A}_j$ and $\widetilde{B}_j$ are the outputs of $j$ steps of finite-precision \textbf{IRS}, beginning with the input matrices $\widetilde{A}_0 = A$ and $\widetilde{B}_0 = B$. 
\begin{lem} \label{lem: norms grow at step j}
    At any step $j$, $|| \binom{\widetilde{A}_{j+1}}{\widetilde{B}_{j+1}} ||_2 \leq (1 + 2 \tau)^2 || \binom{\widetilde{A}_j}{\widetilde{B}_j} ||_2$.
\end{lem}
\begin{proof}
    Recall that $\widetilde{A}_{j+1} = \textbf{MM}(\widetilde{Q}_{12}^H,\widetilde{A}_j)$ and $\widetilde{B}_{j+1} = \textbf{MM}(\widetilde{Q}_{22}^H, \widetilde{B}_j)$ for $\widetilde{Q}_{12}$ and $\widetilde{Q}_{22}$ blocks of a nearly unitary $\widetilde{Q}$ obtained by computing a finite-precision, full QR factorization of $\binom{\widetilde{B}_j}{-\widetilde{A}_j}$. With this in mind, write
    \begin{equation}
        \left| \left| \binom{\widetilde{A}_{j+1}}{\widetilde{B}_{j+1}}\right| \right|_2 \leq \left| \left| \binom{\widetilde{A}_{j+1} - \widetilde{Q}_{12}^H\widetilde{A}_j}{\widetilde{B}_{j+1} - \widetilde{Q}_{22}^H\widetilde{B}_j} \right| \right|_2 + \left| \left| \binom{\widetilde{Q}_{12}^H\widetilde{A}_j}{\widetilde{Q}_{22}^H\widetilde{B}_j} \right| \right|_2 .
        \label{eqn: triangle to bound norm}
    \end{equation}
    By \cref{MM assumption}, $||\widetilde{A}_{j+1} - \widetilde{Q}_{12}^H\widetilde{A}_j||_2 \leq \tau ||\widetilde{Q}_{12}||_2||\widetilde{A}_j||_2 $ and $||\widetilde{B}_{j+1} - \widetilde{Q}_{22}^H\widetilde{B}_j||_2 \leq \tau ||\widetilde{Q}_{22}||_2||\widetilde{B}_j||_2$, so
    \begin{equation}
        \left| \left| \binom{\widetilde{A}_{j+1} - \widetilde{Q}_{12}^H\widetilde{A}_j}{\widetilde{B}_{j+1} - \widetilde{Q}_{22}^H\widetilde{B}_j}\right| \right|_2 \leq \sqrt{2}\tau (1+\tau) \left| \left| \binom{\widetilde{A}_j}{\widetilde{B}_j} \right| \right|_2 
        \label{eqn: term one in triangle}
    \end{equation}
    since $\widetilde{Q}_{12}$ and $\widetilde{Q}_{22}$ satisfy $||\widetilde{Q}_{12}||_2, ||\widetilde{Q}_{22}||_2 \leq ||\widetilde{Q}||_2 \leq 1+\tau$ and $||\widetilde{A}_j||_2, ||\widetilde{B}_j||_2 \leq || \binom{\widetilde{A}_j}{\widetilde{B}_j} ||_2$. Similarly,
    \begin{equation}
        \left| \left| \binom{\widetilde{Q}_{12}^H \widetilde{A}_j}{\widetilde{Q}_{22}^H \widetilde{B}_j} \right| \right|_2 \leq \left| \left| \begin{pmatrix} \widetilde{Q}_{12}^H & 0 \\ 0 & \widetilde{Q}_{22}^H \end{pmatrix} \right| \right|_2 \left| \left| \binom{\widetilde{A}_j}{\widetilde{B}_j} \right| \right|_2  \leq (1+\tau) \left| \left| \binom{\widetilde{A}_j}{\widetilde{B}_j} \right| \right|_2 .
        \label{eqn: term two in triangle}
    \end{equation}
    We obtain the final inequality by combining \eqref{eqn: term one in triangle} and \eqref{eqn: term two in triangle} and using the loose\footnote{We use this bound for convenience to simplify constants. As we will see, it does not significantly impact the final result.} upper bound $(\sqrt{2}\tau+1)(1+\tau) \leq (1 + 2 \tau)^2$.
\end{proof}

\indent We now derive the main bound of \cref{thm: main_IRS_bound}. The high-level strategy of its proof, which again is due to Malyshev \cite{Malyshev1}, can be summarized as follows. Consider the $2^pn \times (2^p+1)n$ block matrix
\begin{equation}
    M = \begin{pmatrix} -A & B & & &  \\
    & - A & B & &  \\
    & & \ddots & \ddots &  \\
    & & & -A & B \\
    \end{pmatrix}.
    \label{eqn: block_M}
\end{equation}
As we demonstrate below, the floating-point matrices used to obtain $\widetilde{A}_p$ and  $\widetilde{B}_p$ via \textbf{IRS} can be built into an approximate block QR factorization of $M$ (containing $\widetilde{A}_p$ and $\widetilde{B}_p$). Our goal will be to derive a nearby, \textit{exact} block QR factorization of $M$, which will contain exact outputs $\mathring{A}_p$ and $\mathring{B}_p$ with bounds on $||\widetilde{A}_p - \mathring{A}_p||_2$ and $|| \widetilde{B}_p - \mathring{B}_p||_2$ available. Since it will be relevant later on, note that the middle $2^pn \times (2^p-1)n$ block of $M$ is the matrix $D_{(A,B)}^p$ from \cref{defn: irs_condition_num}.  \\
\indent To improve readability, we break this proof into several steps, each of which is labeled in bold. \\

\noindent \textbf{Step One: What happens if we apply the output of \textbf{IRS} to $M$ in blocks?} 
\begin{tcolorbox}[blanker, breakable, width=\textwidth, left=0cm, top=2pt, bottom=0pt, before skip=2pt]
\noindent Consider the first iteration of \textbf{IRS}, which computes 
\begin{equation}
    \left[ \widetilde{Q}_1, \binom{R_1}{0} \right] = \textbf{QR} \left( \begin{bmatrix} B \\ -A \end{bmatrix} \right)
    \label{eqn: irs_step_one}
\end{equation}
for nearly unitary $\widetilde{Q}_1 \in {\mathbb C}^{2n \times 2n}$ and upper triangular $R_1 \in {\mathbb C}^{n \times n}$. Let $\widetilde{P}_1$ be the matrix
\begin{equation}
    \widetilde{P}_1 = \begin{pmatrix} \widetilde{Q}_1 & & \\
    & \ddots & \\
    & & \widetilde{Q}_1 \end{pmatrix} \in {\mathbb C}^{2^pn \times 2^pn}
    \label{eqn: first_finite_projector}
\end{equation}
containing $2^{p-1}$ copies of $\widetilde{Q}_1$ on its diagonal. Further, let $\widetilde{M}_1$ be a floating-point approximation of $\widetilde{P}_1^H M$ obtained by applying \textbf{MM} (and finite-precision matrix addition) in $n \times n$ blocks. It is easy to see that $\widetilde{M}_1$ has block structure
\begin{equation}
    \widetilde{M}_1 = \begin{pmatrix} * & \widetilde{R}_1 & * & & & & & &\\
    - \widetilde{A}_1 & \widetilde{E}_1 & \widetilde{B}_1 & & & & & & \\
    & & * & \widetilde{R}_1 & * & & & & \\
    & & - \widetilde{A}_1 & \widetilde{E}_1 & \widetilde{B}_1 & & & & \\
    & & & & \ddots & \ddots & & &   \\
    & & & & & * & \widetilde{R}_1 & * &  \\
    & & & & & - \widetilde{A}_1 & \widetilde{E}_1 & \widetilde{B}_1 &  \\
    \end{pmatrix},
    \label{eqn: irs_after_one}
\end{equation}
where $*$ blocks are arbitrary. We use finite-arithmetic block matrix multiplication here -- as opposed to a separate black-box algorithm for large, non-square matrices -- to guarantee that $\widetilde{A}_1$ and $\widetilde{B}_1$ appear in \eqref{eqn: irs_after_one}. Note that the zero blocks of $\widetilde{M}_1$ are computed exactly by \textbf{MM}. Moreover, $\widetilde{R}_1$ and $\widetilde{E}_1$ are covered by \cref{lem: QR_backward_R} -- i.e., 
\begin{equation}
    ||\widetilde{R}_1 - R_1||_2, ||\widetilde{E}_1||_2 \leq 16 \tau \left| \left| \binom{A}{B} \right| \right|_2.
    \label{eqn: apply_technical_lem}
\end{equation} 
\end{tcolorbox} 

\pagebreak

\noindent \textbf{Step Two: Repeat this argument for the next iteration}.
\begin{tcolorbox}[blanker, breakable, width=\textwidth, left=0cm, top=2pt, bottom=5pt, before skip=2pt]
\noindent The second step of \textbf{IRS} computes
\begin{equation}
    \left[ \widetilde{U}, \binom{R_2}{0} \right] = \textbf{QR} \left( \begin{bmatrix} \widetilde{B}_1 \\ - \widetilde{A}_1 \end{bmatrix} \right),
    \label{eqn: irs_step_two}
\end{equation}
for $\widetilde{U} \in {\mathbb C}^{2n \times 2n}$ and $R_2 \in {\mathbb C}^{n \times n}$. Breaking $\widetilde{U}$ into $n \times n$ blocks $\widetilde{U} = \begin{pmatrix} \widetilde{U}_{11} & \widetilde{U}_{12} \\ \widetilde{U}_{21} & \widetilde{U}_{22} \end{pmatrix}$ and constructing
\begin{equation}
    \widetilde{Q}_2 = \begin{pmatrix} I_n & & & \\ & \widetilde{U}_{11} & & \widetilde{U}_{12} \\ & & I_n & \\ & \widetilde{U}_{21} & & \widetilde{U}_{22} \end{pmatrix} \in {\mathbb C}^{4n \times 4n},
    \label{eqn: Q_2}
\end{equation}
let $\widetilde{P}_2$ be the matrix containing $2^{p-2}$ copies of $\widetilde{Q}_2$ on its diagonal -- that is,
\begin{equation}
    \widetilde{P}_2 = \begin{pmatrix} \widetilde{Q}_2 & & \\ & \ddots & \\ & & \widetilde{Q}_2 \end{pmatrix} \in {\mathbb C}^{2^pn \times 2^pn}.
    \label{eqn: second_finite projector}
\end{equation}
Again applying \textbf{MM} in $n \times n$ blocks to left-multiply $\widetilde{M}_1$ by $\widetilde{P}_2$,  we obtain $\widetilde{M}_2$ -- a finite-precision version of $\widetilde{P}_2^H \widetilde{M}_1$ consisting of $2^{p-2}$ blocks of the form
\begin{equation}
\begin{pmatrix}
        * & \widetilde{R}_1 & * & &  \\
        * & * & \widetilde{R}_2 & * & *  \\
        & & * & \widetilde{R}_1 & * \\
        - \widetilde{A}_2 & \textbf{MM}(\widetilde{U}_{12}^H, \widetilde{E}_1) & \widetilde{E}_2 & \textbf{MM}(\widetilde{U}_{22}^H, \widetilde{E}_1) & \widetilde{B}_2 \\
    \end{pmatrix}.
    \label{eqn: irs_after_two}
\end{equation} 
\end{tcolorbox}

\vspace{3mm}

\noindent \textbf{Step Three: Generalize to an arbitrary step of IRS.} 
\begin{tcolorbox}[blanker, breakable, width=\textwidth, left=0cm, top=2pt, bottom=5pt, before skip=2pt]
\indent The process outlined above yields a sequence of $2^pn \times (2^p+1)n$ matrices $\widetilde{M}_1, \ldots, \widetilde{M}_p$. Each $\widetilde{M}_i$ is an approximation of the exact product $\widehat{M}_i = \widetilde{P}_i^H \widetilde{P}_{i-1}^H \cdots \widetilde{P}_1^H M$ for a corresponding set of nearly unitary matrices $\widetilde{P}_1, \ldots, \widetilde{P}_p$, each of which is constructed from the blocks of a $2n \times 2n$ nearly unitary matrix as in \eqref{eqn: Q_2}. Moreover, $\widetilde{M}_i$ consists of $2^{p-i}$ blocks with structure
\begin{equation}
    \begin{pmatrix} * & * & * \\
    - \widetilde{A}_i & \widetilde{\Delta}_i & \widetilde{B}_i \\
    \end{pmatrix}
    \label{eqn: general_block_form}
\end{equation}
for $\widetilde{\Delta}_i$ a small $n \times (2^i-1) n$ matrix. Indeed, the center $n \times n$ block $\widetilde{E}_i$ of $\widetilde{\Delta}_i$ satisfies
\begin{equation}
    ||\widetilde{E}_i||_2 \leq 16 \tau \left| \left| \binom{\widetilde{B}_{i-1}}{-\widetilde{A}_{i-1}} \right| \right|_2 \leq 16 \tau (1 + 2 \tau)^{2i-2} \left| \left| \binom{A}{B} \right| \right|.
    \label{eqn: delta_block_bound}
\end{equation} 
\end{tcolorbox}

\vspace{3mm}

\noindent \textbf{Step Four: Construct a corresponding set of exact-arithmetic block matrices.} 
\begin{tcolorbox}[blanker, breakable, width=\textwidth, top=2pt, left = 0pt, bottom=5pt, before skip=2pt,parbox=false]
\noindent Suppose that $\widetilde{P}_i$ is constructed from the blocks of the nearly unitary matrix $\widetilde{Q} \in {\mathbb C}^{2n \times 2n}$. Since we use \textbf{QR} to obtain $\widetilde{Q}$, as described above, we know by \cref{QR assumption} that there exists a truly unitary matrix $Q \in {\mathbb C}^{2n \times 2n}$ such that $||\widetilde{Q} - Q||_2 \leq \mu_{\text{QR}}(n) {\bf u} \leq \tau$. With this in mind, let $P_i$ be the truly unitary matrix that has the same block structure as $\widetilde{P}_i$ but swaps the blocks of $\widetilde{Q}$ for the corresponding blocks of $Q$ and define the $2^pn \times (2^p + 1)n$ matrices $M_i = P_i^H \cdots P_1^H M$. \\
\indent We now have two sets of exact-arithmetic matrices to work with: $M_i$ and $\widehat{M}_i$. The former can be thought of as an exact-arithmetic counterpart of $\widetilde{M}_i$ while $\widehat{M}_i$ is an intermediate matrix, obtained via exact multiplication with the nearly unitary $\widetilde{P}_i$. Since $||P_i - \widetilde{P}_i||_2 \leq \tau$ by construction, we can easily bound $||\widehat{M}_i - M_i||_2$ recursively:
\begin{equation}
    \aligned 
    ||\widehat{M}_i - M_i||_2 &= || \widehat{M}_i - P_i^H\widehat{M}_{i-1} + P_i^H \widehat{M}_{i-1} - M_i||_2 \\
    & \leq ||\widehat{M}_i - P_i^H \widehat{M}_{i-1}||_2 + ||P_i^H\widehat{M}_{i-1} - M_i||_2 \\
    & \leq ||\widetilde{P}_i - P_i||_2||\widehat{M}_{i-1}||_2 + ||P_i||_2 ||\widehat{M}_{i-1} - M_{i-1}||_2 \\
    & \leq \tau (1+\tau)^{i-1}||M||_2 + ||\widehat{M}_{i-1} - M_{i-1}||_2
    \endaligned 
    \label{eqn: hat_minus_exact}
\end{equation}
The base case here is $||\widehat{M}_1 - M_1||_2 \leq ||\widetilde{P}_1 - P_1||_2||M||_2 \leq \tau ||M||_2$, so by induction we obtain
\begin{equation}
    ||\widehat{M}_i - M_i||_2 \leq \left( \sum_{j=0}^{i-1} (1+\tau)^j \right) \tau ||M||_2 = \left[ (1 + \tau)^i - 1 \right] ||M||_2.
    \label{eqn: inductive_bound}
\end{equation}
\indent Note that $\widehat{M}_i$ and $M_i$ have the same block structure as $\widetilde{M}_i$, a consequence of the fact that each $P_i$ has the same block structure as $\widetilde{P}_i$. Following \eqref{eqn: general_block_form}, label the blocks of $\widehat{M}_i$ and $M_i$ as 
\begin{equation}
    \begin{pmatrix} * & * & * \\ - \widehat{A}_i & \widehat{\Delta}_i & \widehat{B}_i \\ \end{pmatrix} \; \; \; \;  \text{and} \; \; \; \; \begin{pmatrix} * & * & * \\ - A_i & \Delta_i & B_i \\ \end{pmatrix}, 
    \label{eqn: extend_block_notation}
\end{equation}
respectively, where again $*$ blocks are arbitrary and $\widehat{\Delta}_i, \Delta_i \in {\mathbb C}^{n \times (2^i - 1)n}$.
\end{tcolorbox}

\vspace{3mm}

\noindent \textbf{Step Five: Bound $||\widetilde{A}_i - A_i||_2$ and $||\widetilde{B}_i - B_i||_2$.}
\begin{tcolorbox}[blanker, breakable, width=\textwidth, top=2pt, left = 0pt, bottom=5pt, before skip=2pt,parbox=false]
The matrices $A_i$ and $B_i$ in \eqref{eqn: extend_block_notation} are not necessarily the result of applying $i$ steps of exact-arithmetic repeated squaring to $A$ and $B$, as the matrices $P_i$ used to obtain $A_i$ and $B_i$ from $M$ do not correspond to true QR factorizations of exact outputs. Rather, \cref{QR assumption} implies that each $P_i$ is constructed from the Q-factor of a matrix nearby $\binom{\widetilde{B}_j}{-\widetilde{A}_j}$. With this in mind, we next bound $||\widetilde{A}_i - A_i||_2$  and $||\widetilde{B}_i - B_i||_2$. \\
\indent Consider first $||\widetilde{A}_i - A_i||_2$. Since 
\begin{equation}
    || \widetilde{A}_i - A_i||_2 \leq || \widetilde{A}_i - \widehat{A}_i||_2 + || \widehat{A}_i - A_i||_2 \leq || \widetilde{A}_i - \widehat{A}_i||_2 + || \widehat{M}_i - M_i||_2 ,
    \label{eqn: A_i_triangle}
\end{equation} 
and given \eqref{eqn: inductive_bound}, we can bound $|| \widetilde{A}_i - A_i||_2$ via \eqref{eqn: A_i_triangle} by bounding $|| \widetilde{A}_i - \widehat{A}_i||_2$, which records the error due to finite-precision, block matrix multiplication. With this in mind, suppose $\widetilde{A}_i = \textbf{MM}(\widetilde{Q}_{12}^H, \widetilde{A}_{i-1})$ for $\widetilde{Q}_{12}$ an $n \times n$ block of a nearly unitary $2n \times 2n$ matrix. In this case, $\widehat{A}_i = \widetilde{Q}_{12}^H \widehat{A}_{i-1}$ and we have
\begin{equation}
    \aligned
    ||\widetilde{A}_i - \widehat{A}_i||_2 & = ||\textbf{MM}(\widetilde{Q}_{12}^H, \widetilde{A}_{i-1}) - \widetilde{Q}_{12}^H \widehat{A}_{i-1}||_2 \\
    & = || \textbf{MM}(\widetilde{Q}_{12}^H, \widetilde{A}_{i-1}) - \widetilde{Q}_{12}^H \widetilde{A}_{i-1} + \widetilde{Q}_{12}^H \widetilde{A}_{i-1} - \widetilde{Q}_{12}^H \widehat{A}_{i-1}||_2 \\
    & \leq || \textbf{MM}(\widetilde{Q}_{12}^H, \widetilde{A}_{i-1}) - \widetilde{Q}_{12}^H \widetilde{A}_{i-1}||_2 + ||\widetilde{Q}_{12}^H \widetilde{A}_{i-1} - \widetilde{Q}_{12}^H \widehat{A}_{i-1}||_2 .
    \endaligned 
    \label{eqn: error_in_block_MM1}
\end{equation}
Applying our black-box assumptions and \cref{lem: norms grow at step j}, we have
\begin{equation}
    \aligned 
    || \widetilde{A}_i - \widehat{A}_i||_2 & \leq \tau ||\widetilde{Q}_{12}||_2 || \widetilde{A}_{i-1}||_2 + ||\widetilde{Q}_{12}||_2|| \widetilde{A}_{i-1} - \widehat{A}_{i-1}||_2 \\
    & \leq \tau (1 + \tau) \left| \left| \binom{\widetilde{A}_{i-1}}{\widetilde{B}_{i-1}} \right| \right|_2 + (1 + \tau) || \widetilde{A}_{i-1} - \widehat{A}_{i-1} ||_2 \\
    & \leq \tau (1 + \tau) (1 + 2 \tau)^{2i-2} \left| \left| \binom{A}{B} \right| \right|_2 + (1 + \tau) || \widetilde{A}_{i-1} - \widehat{A}_{i-1}||_2.
    \endaligned
    \label{eqn: error_in_block_MM2}
\end{equation}
Once again we obtain a recursive bound. In this notation $\widetilde{A}_0 = \widehat{A}_0 = A$, so the base case here is simply the error in one finite-precision $n \times n$ matrix multiplication -- i.e.,  
\begin{equation}
    ||\widetilde{A}_1 - \widehat{A}_1||_2 \leq \tau (1+\tau) || A||_2 \leq \tau (1 +\tau) \left| \left| \binom{A}{B} \right| \right|_2.
    \label{eqn: A_tilde_base_case}
\end{equation}
Thus, we conclude inductively
\begin{equation}
    \aligned 
    || \widetilde{A}_i -\widehat{A}_i||_2 &\leq \tau \left( \sum_{j=1}^i ( 1+\tau)^{i-j+1} (1 + 2 \tau)^{2j-2} \right) \left| \left| \binom{A}{B} \right| \right|_2 \\
    & = \tau (1 + 2 \tau)^{2i} \left(\sum_{j=1}^i \left[ \frac{1 + \tau}{(1 + 2 \tau)^2} \right]^{i-j+1}  \right) \left| \left| \binom{A}{B} \right| \right|_2 \\
    & = \tau (1 + 2\tau)^{2i} \cdot \frac{1 + \tau}{(1+2\tau)^2} \cdot \frac{1 - \left( \frac{1 + \tau}{(1+2\tau)^2} \right)^i}{1 - \frac{1 + \tau}{(1+2\tau)^2}} \cdot \left| \left| \binom{A}{B} \right| \right|_2 \\
    & =  \frac{1 + \tau}{3 + 4 \tau}\left[ (1 + 2 \tau)^{2i} - (1+\tau)^i \right] \left| \left| \binom{A}{B} \right| \right|_2.
    \endaligned 
    \label{eqn: inductive_A_error}
\end{equation} 
Combining this with \eqref{eqn: inductive_bound} and $\frac{1 + \tau}{3 + 4 \tau} < 1$, and noting $||M||_2 \leq ||A||_2 + ||B||_2 \leq 2 || \binom{A}{B}||_2$, we have
\begin{equation}
    \aligned 
    || \widetilde{A}_i - A_i||_2 &\leq \left[ (1 + 2 \tau)^{2i} - (1+\tau)^i \right] \left| \left| \binom{A}{B} \right| \right|_2 + 2 \left[(1+\tau)^i - 1 \right] \left| \left| \binom{A}{B} \right| \right|_2 \\
    & = \left[(1+2 \tau)^{2i} + (1+\tau)^i - 2 \right] \left| \left| \binom{A}{B} \right| \right|_2.
    \endaligned
    \label{eqn: final_A_bound}
\end{equation}
Repeating this argument implies the same bound for $||\widetilde{B}_i - B_i||_2$.
\end{tcolorbox}

\vspace{3mm}

\noindent \textbf{Step Six: Show that $||\Delta_i||_2$ is small.}
\begin{tcolorbox}[blanker, breakable, width=\textwidth, top=2pt, left = 0pt, bottom=5pt, before skip=2pt,parbox=false]
If $M_i$ \textit{was} obtained from $M$ via exact-arithmetic repeated squaring, we would have $\Delta_i = 0$. Hence, the norm of $\Delta_i$ is an indication of how far $A_i$ and $B_i$ are from exact outputs. With this in mind, we next derive a bound on $||\Delta_i||_2$. \\
\indent We start by bounding $||\widehat{\Delta}_i||_2$, beginning with its middle $n \times n$ block $\widehat{E}_i$, which corresponds to $\widetilde{E}_i$ in $\widetilde{M}_i$. If we again assume that $\widetilde{P}_i$ is built from the $2n \times 2n$ nearly unitary matrix $\widetilde{Q} = \begin{pmatrix} \widetilde{Q}_{11} & \widetilde{Q}_{12} \\ \widetilde{Q}_{21} & \widetilde{Q}_{22} \end{pmatrix}$, we know that $\widetilde{E}_i$ is the finite-arithmetic sum of $\textbf{MM}(\widetilde{Q}_{12}^H, \widetilde{B}_{i-1})$ and $\textbf{MM}(\widetilde{Q}_{22}^H, -\widetilde{A}_{i-1})$ while $\widehat{E}_i = \widetilde{Q}_{12}^H \widehat{B}_{i-1} - \widetilde{Q}_{22}^H \widehat{A}_{i-1}$. Hence, we have
\begin{equation}
    \aligned 
    ||\widetilde{E}_i - \widehat{E}_i||_2 \leq& \; ||\widetilde{E}_i - (\textbf{MM}(\widetilde{Q}_{12}^H, \widetilde{B}_{i-1}) + \textbf{MM}(\widetilde{Q}_{22}^H, - \widetilde{A}_{i-1}))||_2 \\
    &+ ||\textbf{MM}(\widetilde{Q}_{12}^H, \widetilde{B}_{i-1}) - \widetilde{Q}_{12}^H \widehat{B}_{i-1}||_2 + ||\textbf{MM}(\widetilde{Q}_{22}^H, \widetilde{A}_{i-1}) - \widetilde{Q}_{22}^H \widehat{A}_{i-1} ||_2.
    \endaligned 
    \label{eqn: triangle_for_E}
\end{equation}
By \cref{lem: finite_addition}, the first term in this expression can be bounded by 
\begin{equation}
    \aligned 
    \tau || \textbf{MM}(&\widetilde{Q}_{12}^H, \widetilde{B}_{i-1}) + \textbf{MM}(\widetilde{Q}_{22}^H, - \widetilde{A}_{i-1})||_2 \\
    &\leq \tau \left[ || \textbf{MM}(\widetilde{Q}_{12}^H, \widetilde{B}_{i-1})||_2 + ||\textbf{MM}(\widetilde{Q}_{22}^H, - \widetilde{A}_{i-1})||_2 \right] \\
    & \leq \tau \left[||\widetilde{Q}_{12}^H \widetilde{B}_{i-1}||_2 + \tau || \widetilde{Q}_{12}||_2||\widetilde{B}_{i-1}||_2 + || \widetilde{Q}_{22}^H\widetilde{A}_{i-1}||_2 + \tau ||\widetilde{Q}_{22}||_2 || \widetilde{A}_{i-1}||_2 \right] \\
    & \leq \tau (1 + \tau)^2 \left[ || \widetilde{B}_{i-1}||_2 + || \widetilde{A}_{i-1}||_2 \right] \\
    & \leq 2 \tau ( 1 + \tau)^2 \left| \left| \binom{\widetilde{A}_{i-1}}{\widetilde{B}_{i-1}} \right| \right|_2 \leq 2\tau(1+\tau)^2 (1 + 2\tau)^{2i-2} \left| \left| \binom{A}{B} \right| \right|_2, \\
    \endaligned
    \label{eqn: E_term_one}
\end{equation}
where the last inequality follows from \cref{lem: norms grow at step j}. Using \eqref{eqn: inductive_A_error}, the remaining terms of \eqref{eqn: triangle_for_E} satisfy the following: 
\begin{equation}
    \aligned 
    || \textbf{MM}(\widetilde{Q}_{22}^H, \; &\widetilde{A}_{i-1}) - \widetilde{Q}_{22}^H \widehat{A}_{i-1} ||_2 \\
    &\leq || \textbf{MM}(\widetilde{Q}_{22}^H, \widetilde{A}_{i-1}) - \widetilde{Q}_{22}^H \widetilde{A}_{i-1}||_2 + ||\widetilde{Q}_{22}^H \widetilde{A}_{i-1} - \widetilde{Q}_{22}^H \widehat{A}_{i-1}||_2 \\
    & \leq \tau ||\widetilde{Q}_{22}||_2 ||\widetilde{A}_{i-1}||_2 + ||\widetilde{Q}_{22}||_2 || \widetilde{A}_{i-1} - \widehat{A}_{i-1} ||_2 \\
    & \leq \tau ( 1 + \tau) || \widetilde{A}_{i-1}||_2 + (1 + \tau) || \widetilde{A}_{i-1} - \widehat{A}_{i-1}||_2 \\
    & \leq \left[(1 + \tau)^2(1 + 2\tau)^{2i-2} - (1 + \tau)^i \right] \left| \left| \binom{A}{B} \right| \right|_2.
    \endaligned
    \label{eqn: E_term_two}
\end{equation}
Putting everything together, we obtain
\begin{equation}
    \aligned
    ||\widetilde{E}_i - \widehat{E}_i||_2 \leq 2 \left[ (1+\tau)^3(1+2\tau)^{2i-2} - (1+\tau)^i \right] \left| \left| \binom{A}{B} \right| \right|_2.
    \endaligned 
    \label{eqn: final_overline_E_bound}
\end{equation}
To extend this bound to all of $\widehat{\Delta}_i$, note that
\begin{equation}
    \widehat{\Delta}_i = \begin{pmatrix} \widetilde{Q}_{12}^H \widehat{\Delta}_{i-1} & \widehat{E}_i & \widetilde{Q}_{22}^H \widehat{\Delta}_{i-1} \end{pmatrix}
    \label{eqn: overline_delta_blocks}
\end{equation}
for the same $\widetilde{Q}_{12}$ and $\widetilde{Q}_{22}$ used above. Hence, applying both \eqref{eqn: delta_block_bound} and \eqref{eqn: final_overline_E_bound}, we have
\begin{equation}
    \aligned
    || \widehat{\Delta}_i||_2 &\leq \left| \left| \begin{bmatrix} \widetilde{Q}_{12}^H\widehat{\Delta}_{i-1} & \widetilde{Q}_{22}^H \widehat{\Delta}_{i-1} \end{bmatrix} \right| \right|_2 + || \widehat{E}_i||_2 \\
    & \leq \left| \left| \begin{bmatrix} \widetilde{Q}_{12}^H & \widetilde{Q}_{22}^H \end{bmatrix} \right| \right|_2 || \widehat{\Delta}_{i-1}||_2 + ||\widetilde{E}_i||_2 + ||\widetilde{E}_i - \widehat{E}_i||_2 \\
    & \leq (1 + \tau)|| \widehat{\Delta}_{i-1}||_2 + 2\left[ (8\tau + (1+\tau)^3)(1+2\tau)^{2i-2} - (1+\tau)^i \right] \left| \left| \binom{A}{B} \right| \right|_2 \\
    & \leq (1+\tau)|| \widehat{\Delta}_{i-1}||_2 + 2\left[ (1+15\tau)(1+2\tau)^{2i-2} - (1+\tau)^{i} \right] \left| \left| \binom{A}{B} \right| \right|_2,
    \endaligned 
    \label{eqn: overline_delta_norm_bound}
\end{equation}
where we obtain the final inequality via $8\tau + (1+\tau)^3 \leq 1 + 15 \tau$. Observing $|| \widehat{\Delta}_1||_2 = ||\widehat{E}_1||_2 \leq 28 \tau || \binom{A}{B}||_2$, \eqref{eqn: overline_delta_norm_bound} implies inductively
\begin{equation}
    || \widehat{\Delta}_i||_2 \leq 2(1+\tau)^i\left[ \frac{14 \tau}{1+\tau}+ (1+15 \tau) \sum_{j=1}^{i-1} \frac{(1+2\tau)^{2j}}{(1+\tau)^{j+1}} - (i-1) \right] \left| \left| \binom{A}{B} \right| \right|_2.
    \label{eqn: inductive_overline_delta_bound}
\end{equation}
We therefore conclude,
\begin{equation}
    || \widehat{\Delta}_i||_2 \leq 2 \left( 14 \tau(1+\tau)^{i-1} + (i-1) \left[ (1+15 \tau)(1+2\tau)^{2i} - (1+\tau)^{i} \right] \right) \left| \left| \binom{A}{B} \right| \right|_2,
    \label{eqn: final_overline_delta_norm_bound}
\end{equation}
which we obtain by bounding the sum in \eqref{eqn: inductive_overline_delta_bound} as 
\begin{equation}
    (1+\tau)^i\sum_{j=1}^{i-1} \frac{(1+2\tau)^{2j}}{(1+\tau)^{j+1}} = \frac{(1+2\tau)^{2i}}{1+\tau}\sum_{j=1}^{i-1} \left[ \frac{1+\tau}{(1+2\tau)^2} \right]^{i-j} \leq (i-1)(1+2\tau)^{2i},
    \label{eqn: delta_sum_bound}
\end{equation}
noting $\frac{1+\tau}{(1+2\tau)^2} < 1$. Combining \eqref{eqn: final_overline_delta_norm_bound} with \eqref{eqn: inductive_bound} we have a final bound
\begin{equation}
    \aligned
    ||\Delta_i||_2 &\leq ||\widehat{\Delta}_i||_2 + || \widehat{\Delta}_i - \Delta_i||_2 \\
    &\leq || \widehat{\Delta}_i||_2 + || \widehat{M}_i - M_i||_2 \\
    & \leq 2 \left( 14\tau(1 + \tau)^{i-1} + (i-1)\left[(1+15 \tau)(1+2\tau)^{2i} - (1+\tau)^{i} \right] + \; (1+\tau)^{i} - 1 \right) \left| \left| \binom{A}{B} \right| \right|_2.
    \endaligned
    \label{eqn: bound_Delta_p_norm}
\end{equation}
\end{tcolorbox}

\vspace{3mm}

\noindent \textbf{Step Seven: Obtain $\mathring{A}_p, \mathring{B}_p$ by transforming $M_p$ to block upper triangular.}
\begin{tcolorbox}[blanker, breakable, width=\textwidth, top=2pt, left = 0pt, bottom=5pt, before skip=2pt,parbox=false]
When $i = p$, the matrix $M_i$ consists of only one block of the form \eqref{eqn: extend_block_notation}. Hence, we have shown so far
\begin{equation}
    \widehat{P}_p^H \widehat{P}_{p-1}^H \cdots \widehat{P}_1^H M = \begin{pmatrix} * & * & * \\ -A_p & \Delta_p & B_p \end{pmatrix},
    \label{eqn: i_equals_p_case}
\end{equation}
where each $\widehat{P}_i$ is unitary, $\Delta_p \in {\mathbb C}^{n \times (2^p-1)n}$ is small, and $A_p$ and $B_p$ are close to our finite-precision outputs $\widetilde{A}_p$ and $\widetilde{B}_p$. Letting $\Pi \in {\mathbb C}^{(2^p+1)n \times (2^p+1)n}$ be the permutation matrix that swaps the blocks of \eqref{eqn: i_equals_p_case} containing $-A_p$ and $\Delta_p$, we have constructed an exact, almost-block-QR factorization 
\begin{equation}
    \widehat{P}_p^H \widehat{P}_{p-1} \cdots \widehat{P}_1^HM \Pi = \begin{pmatrix} * & * & * \\ \Delta_p & -A_p & B_p \end{pmatrix}.
    \label{eqn: almost_block_QR}
\end{equation}
Equivalently, recalling that the middle $2^pn \times (2^p-1)n$ block of $M$ is $D_{(A,B)}^p$, we have found an exact factorization $\widehat{P}_p^H \cdots \widehat{P}_1^H D_{(A,B)}^p = \binom{*}{\Delta_p}$. \\
\indent Label the $*$ block of this matrix as $F \in {\mathbb C}^{(2^p-1)n \times (2^p - 1)n}$. By \cref{lem: Malyshev_one}, there exists $S \in {\mathbb C}^{2^pn \times 2^pn}$ such that $I + S$ is unitary,  $(I+S) \binom{F}{\Delta_p} = \binom{F'}{0}$, and
\begin{equation}
    ||S||_2 \leq ||\Delta_p||_2 ||F^{-1}||_2 \leq \frac{||\Delta_p||_2}{\sigma_{\min}(D_{(A,B)}^p) - ||\Delta_p||_2},
    \label{eqn: S_norm_bound}
\end{equation}
assuming $\sigma_{\min}(D_{(A,B)}^p) > ||\Delta_p||_2$. Supposing this is the case, let 
\begin{equation}
    (I+S)\widehat{P}_p \widehat{P}_{p-1} \cdots \widehat{P}_1^H M \Pi = \begin{pmatrix} * & * & * \\
    0 & -\mathring{A}_{p} & \mathring{B}_p \end{pmatrix} 
    \label{eqn: block_QR_2}
\end{equation}
and note
\begin{equation}
    ||\mathring{A}_p - A_p||_2, ||\mathring{B}_p - B_p||_2 \leq ||S||_2 ||M||_2 \leq \frac{2 ||\Delta_p||_2}{\sigma_{\min}(D_{(A,B)}^p) - ||\Delta_p||_2} \left| \left| \binom{A}{B} \right| \right|_2.
    \label{eqn: dot_bound}
\end{equation}
Combining \eqref{eqn: dot_bound} with \eqref{eqn: final_A_bound}, we obtain a final bound
\begin{equation}
    \aligned
    ||\widetilde{A}_p - \mathring{A}_p||_2 &\leq || \widetilde{A}_p - A_p||_2 + ||A_p - \mathring{A}_p||_2  \\
    & \leq \left[ (1+2\tau)^{2p} + (1+\tau)^{p} - 2 + \frac{2 ||\Delta_p||_2}{\sigma_{\min}(D_{(A,B)}^p) - || \Delta_p||_2} \right] \left| \left| \binom{A}{B} \right| \right|_2,
    \endaligned
    \label{eqn: final_bound}
\end{equation}
which also applies to $||\widetilde{B}_p - \mathring{B}_p||_2$. 
\end{tcolorbox}

\vspace{3mm}

\noindent \textbf{Step Eight: Bound $\tau$ by enforcing $||\widetilde{A}_p - \mathring{A}_p||_2, ||\widetilde{B}_p - \mathring{B}_p||_2 \leq \delta || \binom{A}{B} ||_2$.}
\begin{tcolorbox}[blanker, breakable, width=\textwidth, top=2pt, left = 0pt, bottom=5pt, before skip=2pt,parbox=false]
Given \eqref{eqn: final_bound}, we obtain the desired bound on $||\widetilde{A}_p - \mathring{A}_p||_2$ and $||\widetilde{B} - \mathring{B}_p||_2$ provided each of $(1+2\tau)^{2p} - 1$, $(1+\tau)^p - 1$, and $\frac{2 ||\Delta_p||_2}{\sigma_{\min}(D_{(A,B)}^p) - ||\Delta_p||_2}$ is at most $\frac{\delta}{3}$. We focus on the latter, since it is the largest. Here, we note that taking $||\Delta_p||_2 \leq \frac{\delta}{9} \sigma_{\min}(D_{(A,B)}^p)$ guarantees not only that the bound \eqref{eqn: S_norm_bound} holds but also
\begin{equation}
    \frac{2 ||\Delta_p||_2}{\sigma_{\min}(D_{(A,B)}^p) - || \Delta_p||_2} \leq \frac{2\delta}{9 - \delta} < \frac{\delta}{3},
    \label{eqn: Delta_criteria_equiv }
\end{equation}
as desired. Appealing to \eqref{eqn: bound_Delta_p_norm} and \cref{defn: irs_condition_num}, we obtain $||\Delta_p||_2 \leq \frac{\delta}{9} \sigma_{\min}(D_{(A,B)}^p)$ by requiring that each of $14\tau(1+\tau)^{p-1}$, $(p-1)\left[(1+15\tau)(1+2\tau)^{2p} - (1+\tau)^{p}\right]$, and $(1+\tau)^{p}-1$ is bounded by $\frac{\delta}{54 \kappa_{\text{IRS}}(A,B,p)}$. Once again, we focus on the largest of these terms, which in this case is $X = (p-1)\left[(1+15\tau)(1+2\tau)^{2p} - (1+\tau)^{p}\right]$, assuming $p > 1$. \\
\indent We begin by rewriting $X$ as follows:
\begin{equation}
    \aligned
    X  = (p-1)(1+2\tau)^{2p} \left[ 1 + 15 \tau - \left[ 1 - \tau \left( \frac{3 + 4\tau}{(1+2\tau)^2} \right) \right]^p \right].
    \endaligned 
    \label{eqn: rewrite_X}
\end{equation}
Since $\tau \left( \frac{3 + 4 \tau}{(1+2\tau)^2} \right) \leq  1$, we can bound $X$ from above via Bernoulli's inequality
\begin{equation}
    \aligned 
    X &\leq (p-1)(1+2\tau)^{2p} \left[1+15\tau - \left[1 - p\tau \left( \frac{3 + 4 \tau}{(1+2\tau)^2}\right)\right] \right] \\
    & = (p-1)(1+2\tau)^{2p} \tau \left[ 15 + p \left( \frac{3 + 4\tau}{(1+2\tau)^2} \right) \right] \\
    & \leq 3\tau (p-1)(p+5)(1+2\tau)^{2p},
    \endaligned
    \label{eqn: bound_X}
\end{equation}
where the last inequality follows by loosely bounding $\frac{3+4\tau}{(1+2\tau)^2} \leq 3$. Finally assuming $(1+2\tau)^{2p} \leq 2$, we obtain a final bound
\begin{equation}
    X \leq 6\tau(p-1)(p+5) = 6\tau (p^2 + 4p - 5),
    \label{eqn: final_X_bound}
\end{equation}
which implies a criterion on $\tau$:
\begin{equation}
    \tau \leq \frac{\delta}{324 \kappa_{\text{IRS}}(A,B,p) (p^2 + 4p - 5)}.
    \label{eqn: tau_criteria_one}
\end{equation}
Note that if $p = 1$, and therefore $X = 0$, we require instead $15\tau \leq \frac{\delta}{18 \kappa_{\text{IRS}}(A,B,p)}$ above, which is clearly satisfied by \eqref{eqn: tau_criteria_one}. It is similarly not hard to show that this requirement on $\tau$ guarantees the remaining bounds and therefore yields $|| \widetilde{A}_p - \mathring{A}_p||_2, || \widetilde{B}_p - \mathring{B}_p||_2 \leq \delta||\binom{A}{B}||_2$. 
\end{tcolorbox}

\vspace{3mm}

\indent It remains to show that $\mathring{A}_p$ and $\mathring{B}_p$ can be obtained via exact-arithmetic repeated squaring. This follows from \cref{lem: Malyshev_two}; exact-arithmetic repeated squaring implies an alternative block-QR factorization of $M\Pi$, which is equivalent to \eqref{eqn: block_QR_2} up to a rotation/reflection. Since such a rotation/reflection can be baked into the final QR factorization computed by exact-arithmetic repeated squaring (which is agnostic to the specific QR factorizations used), $\mathring{A}_p$ and $\mathring{B}_p$ are indeed exact outputs of repeated squaring satisfying $\mathring{A}_p^{-1} \mathring{B}_p = (A^{-1}B)^{2^p}$.

\section{Zolotarev's Problem on the Circles of Apollonius}\label{appendix: Zolo}
This appendix summarizes existing results on optimal rational approximation to the sign function, specifically with respect to the disks of Apollonius considered in \cref{section: SIGN}. We include these results to justify the following claim:\ the standard Halley or Newton iterations are essentially optimal (in terms of fast convergence to the sign function) in the general setting -- i.e., for arbitrary points $z  \in C_{\alpha}^{\pm}$, which in particular are not known to be real.

\begin{defn}
    ${\mathcal R}_n$ is the set of rational functions $r(z) = p(z)/q(z)$, where $p$ and $q$ are polynomials of degree at most $n$.
\end{defn}

The search for an optimal approximations to the sign function in ${\mathcal R}_{n}$ goes back to work of Zolotarev \cite{Zolo}. We state the relevant problem, tailored to our context, below. 

\begin{mdframed}
    \textbf{Zolotarev's Problem on the Circles of Apollonius:} Given a fixed positive integer $n$, find the rational function $r_n(z)$ satisfying 
    $$ r_n = \argmin_{r \in {\mathcal R}_n} \max_{z \in C_{\alpha}} |r(z)-\text{sign}(z)|.$$
\end{mdframed}

Importantly, a solution to this problem is already known. 

\begin{thm}\label{thm: Zolo_soln}
    The solution $r_n(z)$ to Zolotarev's problem on the circles of Apollonius is 
    $$ r_n(z) = -\frac{1-\alpha^{2n}}{1+\alpha^{2n}} \cdot \frac{(z-1)^n - (z+1)^n}{(z-1)^n + (z+1)^n} $$
    and satisfies
    $$ \max_{z \in C_{\alpha}}|r_n(z) - \text{\normalfont sign(z)}| \leq \frac{2 \alpha^n}{1+\alpha^{2n}}.$$
\end{thm}
\begin{proof}
    This follows from work of Starke \cite{STARKE1992115} and Istace and Thiran \cite{3_4_equivalence} (see also \cite[Section 2]{Bounding_zolo}).
\end{proof}

To put this result in context, recall that $C_{\alpha}^+$ is the disk with center $\frac{1+\alpha^2}{1-\alpha^2}$ and radius $\frac{2\alpha}{1-\alpha^2}$. Hence, any $z \in C_{\alpha}^{+}$ satisfies
\begin{equation}\label{eqn: circle_to_sign}
    |z-\text{sign}(z)| \leq \left| z - \frac{1+\alpha^2}{1-\alpha^2}  + \frac{1+\alpha^2}{1-\alpha^2} - 1 \right| \leq \frac{2\alpha}{1-\alpha^2} + \frac{2\alpha^2}{1-\alpha^2}. 
\end{equation}
Since the rational functions that define the Newton and Halley iterations map $C_{\alpha}^{\pm}$ to $C_{\alpha^2}^{\pm}$ and $C_{\alpha^3}^{\pm}$, respectively, \eqref{eqn: circle_to_sign} and \cref{thm: Zolo_soln} imply that they are optimal in ${\mathcal R}_2$ and ${\mathcal R}_3$ up to an additive factor of, at most, $2\alpha^2/(1-\alpha^2)$. Indeed, when $n = 3$ the optimal rational function $r_3(z)$ is simply the reciprocal of the rational function corresponding to the Halley iteration rescaled by factor of $\frac{1-\alpha^6}{1+\alpha^6}$:
\begin{equation}
    r_3(z) = \frac{1 - \alpha^6}{1+\alpha^6} \cdot \frac{3z^2+1}{z^3+3z}.
\end{equation}
Moreover, this argument implies that the rational functions obtained by iteratively composing the Halley/Newton rationals are similarly (nearly) optimal. As a result, we cannot hope to replicate the fast convergence of \textbf{IF-DWH} in general.





\section{Additional Numerical Examples}\label{appendix: examples}
\begin{figure}[t]
    \centering
    \includegraphics[width=\linewidth]{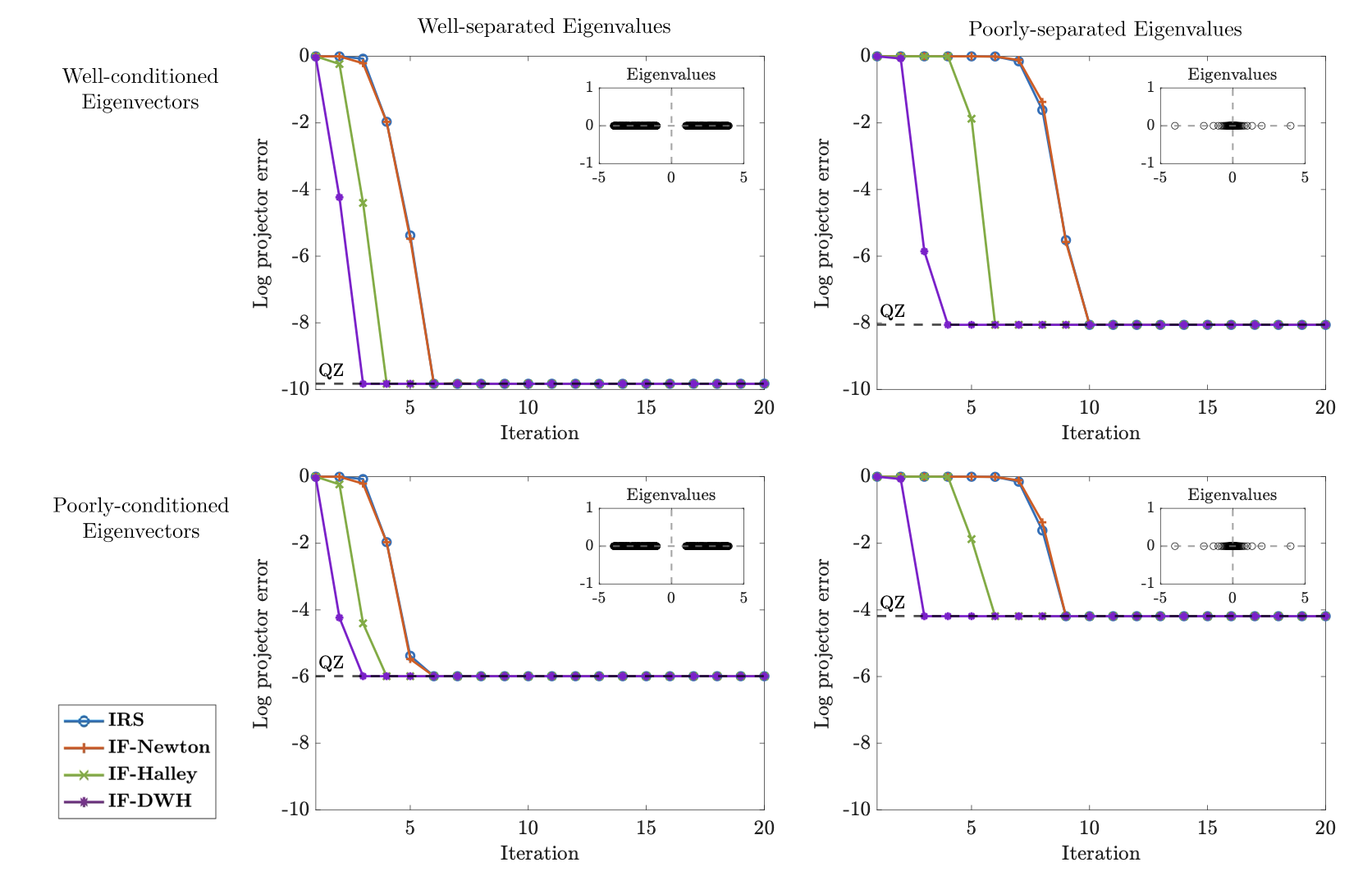}     
    \caption{A repeat of \cref{fig: proj_4x4} for an indefinite pencil $(A,B)$ constructed according to \eqref{eqn: indefinite_example} for the same choices of $\Lambda $ and $X$ as in \cref{section: examples}.}
    \label{fig: indef_proj_4x4}
\end{figure}

This appendix repeats the main numerical example from \cref{section: examples}, which applied the algorithms discussed in \cref{section: IRS,section: SIGN} to a $500 \times 500$ definite pencil $(A,B)$ constructed from the diagonalization \eqref{eqn: example_1} (with various eigenvalue placements and eigenvector conditionings). Here, we modify this experiment by making the input pencil indefinite; for the same construction of $\Lambda$ and $X$, we take $B$ to be random complex Gaussian and set 
\begin{equation}\label{eqn: indefinite_example}
    A = BX\Lambda X^{-1}.
\end{equation}
In this case, the eigenvalues of $(A,B)$ are the diagonal entries of $\Lambda$ (again chosen to be real to verify the efficacy of the dynamically weighted Halley iteration) and the columns of $X$ are right eigenvectors. The latter implies that an ``exact" projector onto the right deflating subspace of $(A,B)$ corresponding to eigenvalues in the right half plane can be obtained from a QR factorization of $X$. \\
\indent Results of this test are shown in \cref{fig: indef_proj_4x4} for the four possible combinations of $\Lambda$ and $X$. The trends here match exactly what was observed in \cref{section: examples} with one exception:\ we can no longer identify subtle differences between the methods once they are converged. This is likely a consequence of error derived from the construction of $A$, which requires matrix inversion. 

\end{document}